\documentclass[graybox]{svmult}

\usepackage{type1cm}        
\usepackage{makeidx}         
\usepackage{graphicx}        
\usepackage{multicol}        
\usepackage[bottom]{footmisc}
\usepackage{url}
\usepackage{amsmath,amsthm,amssymb,enumerate}
\usepackage{euscript,mathrsfs}
\usepackage{dsfont}
\usepackage{bm}
\usepackage{CJKutf8} \newcommand{\chinese}[1]{\begin{CJK*}{UTF8}{gbsn}\! #1 \end{CJK*}}

\usepackage[labelformat=simple]{subcaption}

\usepackage{multirow}

\newcommand{\Td}{\mathbb{T}^d}

\newcommand{\intTd}[1]{\int_{\mathbb{T}^d} #1 \ \dx}

\newcommand{\br}{\\}
\newcommand{\intO}[1]{\int_{\Td} #1 \ \dx}
\newcommand{\R}{\Bbb R}

\newcommand{\pd}{\partial}
\newcommand{\vr}{\varrho}
\newcommand{\vc}[1]{{\bm #1}}
\newcommand{\vu}{\vc{u}}
\newcommand{\vm}{\vc{m}}

\newcommand{\tvr}{\tilde \vr}
\newcommand{\tvu}{{\tilde \vu}}
\newcommand{\tvt}{\tilde \vt}

\newcommand{\vt}{\vartheta}
\newcommand{\Div}{{\rm div}_x}
\newcommand{\Grad}{\nabla_x}

\newcommand{\dx}{\,{\rm d} {x}}
\newcommand{\dt}{\,{\rm d} t }

\newcommand{\tvm}{\widetilde{\vc{m}}}
\newcommand{\tvS}{\widetilde{S}}

\newcommand{\bfphi}{\boldsymbol{\varphi}}
\newcommand{\bfvarphi}{\boldsymbol{\varphi}}

\newcommand{\abs}[1]{\left|#1\right|}

\usepackage[linkcolor=blue,colorlinks=true]{hyperref}

\usepackage{newtxtext}       %
\usepackage[varvw]{newtxmath}       

\def\softd{{\leavevmode\setbox1=\hbox{d}%
          \hbox to 1.05\wd1{d\kern-0.4ex{\char039}\hss}}}

\makeindex             

\definecolor{Cgrey}{rgb}{0.85,0.85,0.85}
\definecolor{Cblue}{rgb}{0.50,0.85,0.85}
\definecolor{Cred}{rgb}{1,0,0}
\definecolor{fancy}{rgb}{0.10,0.85,0.10}
\definecolor{forestgreen}{rgb}{0.13, 0.55, 0.13}
\definecolor{pinegreen}{rgb}{0.0, 0.47, 0.44}

\begin{document}

\title*{What is a limit of structure-preserving numerical methods for compressible flows?}
\titlerunning{Limit of structure-preserving numerical methods}
\author{M\'aria Luk\'a\v{c}ov\'a-Medvi{\softd}ov\'a
\orcidID{0000-0002-4351-0161},
\\ Bangwei She\orcidID{0000-0002-5025-0070}
and\\ Yuhuan Yuan\orcidID{0000-0001-6392-9202} }
\authorrunning{Luk\'a\v{c}ov\'a-Medvi{\softd}ov\'a, She and Yuan}
\institute{Luk\'a\v{c}ov\'a-Medvi{\softd}ov\'a \at
Institute of Mathematics, University of Mainz, Staudingerweg 9, Germany. \email{lukacova@uni-mainz.de}
\and Bangwei She \at Academy for Multidisciplinary studies, Capital Normal University, West 3rd Ring North Road 105, 100048 Beijing, P. R. China. \email{bangweishe@cnu.edu.cn}
\and Yuhuan Yuan \at School of Mathematics, Nanjing University of Aeronautics and Astronautics, Jiangjun Avenue No. 29, 211106 Nanjing, P. R. China. \email{yuhuanyuan@nuaa.edu.cn}
}
%
%
\maketitle

\abstract{
We present an overview of recent developments on the convergence analysis of numerical methods for inviscid multidimensional compressible flows that preserve underlying physical structures. We introduce the concept of generalized solutions, the so-called dissipative solutions, and explain their relationship to other commonly used solution concepts. In numerical experiments we apply $\mathcal{K}-$ convergence of numerical solutions and  approximate turbulent solutions together with the Reynolds stress defect and the energy defect.
}

\section{Mathematical model}
\label{sec:1}
Hyperbolic conservation laws are fundamental for many applications in science, engineering or socioeconomic. An iconic example is represented by  the
Euler equations of gas dynamics expressing the conservation of mass, momentum and energy
\begin{eqnarray}
\label{euler}
\pd_t \vr  + \Div\vm  &=& 0,
\nonumber \\ \nonumber
\pd_t \vm   + \Div(\vm \otimes \vu  ) + \Grad p &=& 0,
\\ \label{Euler}
\pd_t E  +  \Div( (E + p) \vu ) &=& 0,
\end{eqnarray}
where $\vr, \vm$ and $E$ are the conservative variables, representing the density, momentum and the total energy, respectively.
Further,  $p$ and $ \vu = \vm / \vr$ stand for the pressure and velocity. The total energy  $E = \frac 12\frac{\abs{\vm}^2}{ \vr} + \vr e$ consists of the kinetic energy and the internal energy $e$.
Denoting $s$ the specific entropy and $S$ the total entropy, system \eqref{euler} is closed by the pressure law of perfect gas
$$
p(\vr, S) = \vr^\gamma \exp\left(\frac{S}{c_V \vr}\right), \quad e = c_V \vr^{\gamma -1} \exp\left(\frac{S}{c_V \vr}\right) \quad \mbox{ for }\vr >0,
$$
where $\gamma$ is the adiabatic coefficient and $c_V = \frac{1}{\gamma -1}$ the specific heat at constant volume.
In what follows we assume that $\gamma > 1$ and note that physically reasonable range for gases is $1 < \gamma \leq 5/3.$
In addition to conservation  law \eqref{euler} we also consider the second law of thermodynamics that is expressed by the fact that entropy is nondecreasing in time
\begin{equation}
\label{entropy}
\pd_t S + \Div (S \vu) \geq 0.
\end{equation}

Solutions of \eqref{euler}, \eqref{entropy} are known to develop finite-time discontinuities even for infinitely smooth initial data. Since a classical solution may not exists at all,  \eqref{euler}, \eqref{entropy} are considered in the weak (distributional) sense.
However, for the multidimensional compressible Euler equations De Lellis and Sz\'ekelyhidi constructed in their pioneering work~\cite{dlsz2} infinitely many weak entropy solutions, see also Chiodaroli et al.~\cite{kreml}, Feireisl et al.~\cite{feireisl}. Taking into account ill-posedness of multidimensional Euler equations in the class of weak entropy solutions, there is a need to propose new selection criteria to obtain a physically reasonable solution concept. We note that these questions still remain open.

Another research direction looks for a generalized solution concept that is better compatible  with the  vanishing viscosity approach and convergence analysis of well-known numerical methods.

In the last decades we can observe large progress in the development and experimental validation of numerical methods for hyperbolic conservation laws, and in particular the Euler equations. We refer the interested reader to monographs Dolej\v{s}\'{\i} and Feistauer \cite{feist-dolejsi}, Eymard et al. \cite{herbin_FVM}, Kr\"oner~\cite{kroner}, Kuzmin et al.~\cite{kuzmin}, LeVeque \cite{leveque} or Toro~\cite{Toro09} as well as to the handbook by Abgrall and Shu \cite{Abgrall_Shu_1,Abgrall_Shu_2} and further references therein. These methods typically share the common features of being highly accurate, efficient and robust. In addition, they  are typically entropy stable and positivity preserving.
However their numerical analysis, concerning the convergence and error estimates is missing, in general.

Main aim of this paper is to review our recent results on the convergence and error estimates that have been obtained recently for various numerical methods ranging from the finite volume, finite element to discontinuous Galerkin or residual distribution methods, see Feireisl, Luk\'a\v{c}ov\'a-Medvi{\softd}ov\'a et al. \cite{book}, Kuzmin et al.~\cite{kuzmin1}, Luk\'a\v{c}ov\'a-Medvi{\softd}ov\'a and \"Offner \cite{po} and Abgrall et al.~\cite{abgrall}, respectively.

A common feature of all of them is their \emph{structure-preserving property}. This means that numerical solutions satisfy important physical properties, such as the entropy stability and  positivity preservation of density and internal energy.

\section{Dissipative weak solutions}
\label{sec:2}

In view of the ill–posedness of the Euler equations in the class of weak entropy solutions, the relevance of the Euler system to describe the behaviour of
fluids in higher space dimensions may be questionable. In fact, the Euler equations should be seen as an inviscid (vanishing viscosity) limit of real viscous fluids. The low viscosity regime is typical for turbulent flows,  where the solutions may develop oscillatory behaviour. As it follows from \cite{book, FLSS} weak limit of weak solutions to compressible Navier-Stokes equations may not be a weak solution  of the Euler equations. Instead, it is  a generalized,  dissipative weak solution, cf.~\cite{book}.

To be specific,  let us consider \eqref{euler}, \eqref{entropy} on a space-time cylinder $[0,T) \times \Omega,$ $\Omega \subset \Bbb R^d,$ $d=2,3$ and $T>0.$ The system of Euler equations is accompanied with
initial data
$$
\vr(t=0) = \vr_0,\ \ \vm(t=0) = \vm_0,\ E(t=0) = E_0
$$
and periodic boundary conditions. As usual, we identify  domain $\Omega$ with a flat torus  $\Td=\left( [0,1]\Big|_{\{0,1\}} \right)^d.$

\begin{definition}[{\bf Dissipative weak (DW) solution}] \label{DE3}
Let the initial data satisfy
	\begin{eqnarray*}
&& \vr_0 \in L^1(\Td), \ \vm_0 \in L^1(\Td; \R^d),
\\
&& E_0 = E(\vr_0, \vm_0, S_0)\, \mbox{ and }
\intTd{ E(\vr_0, \vm_0, S_0) } < \infty.
	\end{eqnarray*}
	We say that $(\vr, \vm, S)$ is \emph{dissipative weak  solution} to the Euler system
in $[0,T) \times \Td$,	$0 < T \leq \infty$, if the following holds:
	\begin{itemize}
		
		\item {\bf Regularity.}
		The solution $(\vr, \vm, S)$ belongs to the class
		\begin{align*}
			\vr &\in C_{\rm weak, loc}([0,T); L^\gamma(\Td)),\
            \vm \in C_{\rm weak, loc}([0,T); L^{\frac{2\gamma}{\gamma + 1}}(\Td; \R^d)),\br
            S &\in L^\infty(0,T; L^\gamma(\Td)) \cap BV_{\rm weak}([0,T); L^\gamma(\Td)) \br
			&\intTd{ E(\vr, \vm, S)(t, \cdot) } \leq \intTd{ E(\vr_0, \vm_0, S_0) } \ \mbox{ for any }\ 0 \leq t < T.
			\nonumber
		\end{align*}
		\item {\bf Equation of continuity.}
		The integral identity
		\[
		\int_0^T \intTd{ \Big[ \vr \partial_t \varphi + \vm \cdot \Grad \varphi \Big] } \dt =  - \intTd{ \vr_0 \varphi(0, \cdot)}
		\]
		holds for any $\varphi \in C^1_c([0,T) \times \Td)$.
		
		\item {\bf Momentum equation.}  The integral identity
		\begin{align}
			\int_0^T &\intTd{ \left[ \vm \cdot \partial_t \bfphi + \mathds{1}_{\vr > 0} \frac{\vm \otimes \vm}{\vr} : \Grad \bfphi +
				p(\vr, S) \Div\bfphi \right] } \dt \nonumber  \\ &=
			- \int_0^T \intTd{ \Grad \bfphi : \D \mathfrak{R} (t) }  - \intTd{ \vm_0 \cdot \bfphi(0, \cdot) }
			\label{E19}
		\end{align}
		holds for any $\bfphi \in C^1_c([0,T) \times \Td; \R^d)$, where the Reynolds defect stress reads as
		\begin{equation} \label{E20}
			\mathfrak{R} \in  L^\infty(0, T; \mathcal{M}^+(\Td; \R^{d \times d}_{\rm sym})).
		\end{equation}

\item{\bf Entropy inequality.}
\begin{align} \label{DS10}
& \left[ \intO{ S \varphi } \right]_{t = \tau_1-}^{t = \tau_2+} \geq
\int_{\tau_1}^{\tau_2} \intO{ \left[ S \partial_t \varphi + \left< {\mathcal{V}_{t,x}}; 1_{\tvr > 0} \left( \tvS \tvu \right) \right> \cdot \Grad \varphi \right] } \dt, \\
& \hspace{8cm} S(0-, \cdot) = S_0, \nonumber
\end{align}
for any $0 \leq \tau_1 \leq \tau_2 < T$, any $\varphi \in C^1_c([0,T) \times \Td)$, $\varphi \geq 0$,
where
$
\{ {\mathcal{V}}_{t,x} \}_{ (t,x) \in (0,T) \times \Td }
$
is a parametrized probability (Young) measure,
\[
{\mathcal{V}_{t,x}} \in L^\infty((0,T) \times \Td); {\mathcal{P}}(\Bbb R^{d + 2})),\ \Bbb R^{d+2} = \left\{ \tvr \in \Bbb R, \tvm \in \Bbb R^d, \tvS \in \Bbb R \right\};
\]
\begin{equation} \label{DS10a}
\left< {\mathcal{V}}; \tvr \right> = \vr, \ \left< {\mathcal{V}}; \tvm \right> = \vm,\
\left< {\mathcal{V}}; \tvS \right> = S.
\end{equation}
		
		\item{\bf Compatibility of the energy and Reynolds stress defects.} There exists a non--increasing function $\mathcal{E}:
		[0,T) \to [0, \infty)$ satisfying
		\begin{align}
			\mathcal{E}(0-) &= \intTd{ E(\vr_0, \vm_0, S_0) }, \nonumber \br
			\mathcal{E}(\tau+) & =  \intTd{ E(\vr, \vm, S)(\tau, \cdot) } + \mathfrak{E} 	\quad \mbox{ for any } 0 \leq \tau < T,
			\label{E21}	
		\end{align}
		where $\mathfrak{E} \in   L^\infty(0, T; \mathcal{M}^+(\Td)) $ is the energy defect satisfying
		\begin{align}
		 \min \left\{ 2, d(\gamma - 1) \right\} \mathfrak{E}  \leq  {\rm trace} [\mathfrak{R}] \leq  \max \left\{ 2, d(\gamma - 1) \right\}  \mathfrak{E}.
			\nonumber	
		\end{align}

	\end{itemize}
	
\end{definition}

Unlike the weak entropy solutions, the DW solutions are known to exist globally in time. A suitable way to show their global-in-time existence is via convergence analysis of structure-preserving numerical methods as we will discuss in the next section. Moreover, for DW solutions the following properties hold.
\begin{itemize}
\item \textbf{Weak-strong uniqueness}\\
If a strong solution to the Euler equation \eqref{euler} exists, then any DW solution emanated from the same initial data coincides with the strong solution on its lifespan. This result is proved using the relative energy, see B\v{r}ezina and Feireisl~\cite{brezina}
\begin{equation} \label{relative_energy}
\begin{split}
{E} &\left( \vr, \vc{m}, S \Big| \tvr, \tvu, \tvS \right)
= \frac{1}{2} \vr \left| \frac{\vm}{\vr}
- \tvu \right|^2
\\&+  \vr e(\vr, S)
- \frac{\partial (\tvr e(\tvr, \tvS))}{\partial \vr} (\vr - \tvr)
- \frac{\partial (\tvr e(\tvr, \tvS))}{\partial S} (S - \widetilde{S} ) -  \tvr e(\tvr, \tvS),
\end{split}
\end{equation}
that measures a ``distance'' between a DW solution $(\vr, \vc{m}, S)$  and a strong solution $(\tvr, \tvu, \tvS),$
$\tvu = \widetilde \vm / \tvr,$ belonging to the class
\begin{equation} \label{strong_solution}
\begin{split}
\tvr,\ \tvS &\in W^{1, \infty}((0,T) \times \Td),\ \tvu \in W^{1, \infty}((0,T) \times \Td; \R^d) \\
\tvr &\geq \underline{\vr} > 0,\
\tvt \equiv \frac{1}{\gamma - 1}  \frac{\partial p(\tvr, \tvS)}{\partial S}  \geq
\underline{\vt} > 0 \ \mbox{in} \ (0,T) \times \Td.
\end{split}
\end{equation}

\item \textbf{Compatibility}\\
If a DW solution $(\vr, \vu, S) \in C^1([0,T] \times \Td; \Bbb R^{d+2})$, $\inf_{(0,T)\times \Td} \vr > 0,$
$\vu = \vm/\vr$,
then $(\vr, \vm, S)$ is a classical solution of the Euler system, see~\cite{book}. Specifically,
\begin{eqnarray*}
&&\mathfrak{R} = 0,\  \mathfrak{E} = 0,\  \mathcal{E}(t+) = \intTd{ E(\vr, \vm, S)(t, \cdot) } \ \mbox{ for any } t \in [0,T), \\
&&\mathcal{V}_{t,x} = \delta_{[\vr(t,x), \vm(t,x), S(t,x)]}\
\mbox{ for }\
(t,x) \in (0,T) \times \Td.
\end{eqnarray*}

\item \textbf{Semigroup selection} \\
In the class of dissipative solutions one can select a solution maximizing the entropy production. Such a selection  satisfies  semigroup property with respect to time, see Breit et al.~\cite{Breit}.

\item \textbf{Vanishing viscosity limit} \\
In the class of dissipative solutions there is a solution obtained as a vanishing viscosity limit of the compressible Navier-Stokes equations, see \cite{FLSS}.
\end{itemize}

\section{Convergence via Lax equivalence principle}
\label{sec:3}

A celebrated \emph{Lax equivalence principle} states that any consistent and stable numerical method is convergent. In \cite{lax} Lax and Richtmyer  proved this property for linear numerical methods applied to  linear partial differential equations. Using the framework of DW solutions we can generalize the Lax equivalence principle to our nonlinear system of the Euler equations and prove convergence of some structure-preserving numerical methods. To this end, let us
firstly introduce \emph{consistent approximations} of the Euler equations \eqref{euler}, \eqref{entropy}, cf.~\cite{book}.

\begin{definition}[\textbf{Consistent approximation}]
\label{consistent}
A sequence $\left\{ \vr_n, \vm_n, S_n  \right\}_{n=1}^\infty$ is a \emph{consistent approximation} of the Euler system
\eqref{euler}, \eqref{entropy} in $(0,T) \times \Td$ with the initial data
$[\vr_0, \vm_0, S_0]$, $E_0 = E(\vr_0, \vm_0, S_0) \equiv \frac{1}{2} \frac{|\vm_{0}|^2}{\vr_{0}} + \vr_{0} e(\vr_{0}, S_{0}), $ if:
\begin{itemize}
\item {\bf Energy inequality:} there is a sequence $\{ \vr_{0,n}, \vm_{0,n}, S_{0,n} \}_{n=1}^\infty$,
\begin{equation} \label{MV21}
\begin{split}
\vr_{0,n} &\to \vr_0 \ \mbox{weakly in} \ L^1(\Td), \ \vm_{0,n} \to \vm_0 \ \mbox{weakly in}\ L^1(\Td; \R^d),\
\\ S_{0,n} &\to S_0 \ \mbox{weakly in}\ L^1(\Td),
\end{split}
\end{equation}
and
\begin{equation} \label{MV22}
\intO{ \left[ \frac{1}{2} \frac{|\vm_{0,n}|^2}{\vr_{0,n}} + \vr_{0,n} e(\vr_{0,n}, S_{0,n}) \right]  }
\to \intO{ \left[ \frac{1}{2} \frac{|\vm_{0}|^2}{\vr_{0}} + \vr_{0} e(\vr_{0}, S_{0}) \right]  } < \infty
\end{equation}
satisfying
\begin{equation} \label{MV20}
\begin{split}
&\intO{ \left[ \frac{1}{2} \frac{|\vm_n|^2}{\vr_n} + \vr_n e(\vr_n, S_n) \right] (\tau, \cdot) }
\\ &\leq \intO{ \left[ \frac{1}{2} \frac{|\vm_{0,n}|^2}{\vr_{0,n}} + \vr_{0,n} e(\vr_{0,n}, S_{0,n}) \right]  }
+ e^1_n
\end{split}
\end{equation}
for a.a. $0 \leq \tau \leq T$ and
$
e^1_n \to 0\ \mbox{ as } \ n \to \infty;
$
\item {\bf Minimum entropy principle:} there exists $\underline{s} \in \R$ such that
\begin{equation} \label{MV23}
S_n \geq \vr_n \underline{s} \ \mbox{ a.a. in }\ (0,T) \times \Td;
\end{equation}
\item {\bf Equation of continuity:}
\begin{equation} \label{MV24}
\int_0^T \intO{ \left[ \vr_n \partial_t \varphi + \vm_n \cdot \Grad  \varphi \right] } \dt =
- \intO{ \vr_{0,n} \varphi (0, \cdot) } + e^2_n[ \varphi ]
\end{equation}
holds for any $\varphi \in C^1_c([0,T) \times \Td)$ and
$
e^2_n[ \varphi ] \to 0 \ \mbox{as}\ n \to \infty \ \mbox{for any}\ \varphi \in C^2_c([0, T) \times \Td);
$
\item {\bf Momentum equation:}
\begin{align} \label{MV25}
\int_0^T &\intO{ \left[ \vm_n \partial_t \bfvarphi + 1_{\vr_n > 0}\frac{\vm_n \otimes \vm_n}{\vr_n}: \Grad  \bfvarphi
+ 1_{\vr_n > 0} p(\vr_n, S_n) \Div \bfvarphi \right] } \dt \nonumber\\ &=
- \intO{ \vm_{0,n} \bfvarphi (0, \cdot) } + e^3_n[ \bfvarphi ]
\end{align}
holds for any $\bfvarphi \in C^1_c([0,T) \times \Td; \R^d)$ and
$
e^3_n[ \bfvarphi ] \to 0 \ \mbox{as}\ n \to \infty \ \mbox{for any}\ \bfvarphi \in C^2_c([0, T) \times \Td; \R^d);
$
\item {\bf Entropy inequality:}
\begin{equation} \label{MV26}
\int_0^T \intO{ \left[ S_n \partial_t \varphi + 1_{\vr_n > 0} \left( S_n \frac{\vm_n}{\vr_n} \right) \cdot \Grad \varphi \right] } \dt
\leq - \intO{ S_{0,n} \varphi (0, \cdot)} + e^4_n[ \varphi ]
\end{equation}
holds for any $\varphi \in C^1_c([0,T) \times \Td)$, $\varphi \geq 0$, and
$
e^4_n[ \varphi ] \to 0 \ \mbox{as}\ n \to \infty \ \mbox{for any}\ \varphi \in C^2_c([0, T) \times \Td).
$
\end{itemize}
\end{definition}

Roughly speaking, consistent approximations satisfy weak formulation of the Euler equations \eqref{euler}, \eqref{entropy} modulo consistency errors $e^i_n[\varphi],$ $i=1,\dots,4,$ that vanish as $n\to \infty$ for sufficiently smooth test functions $\varphi.$

There are several ways to obtain consistent approximations of the Euler equations. In particular, suitable structure-preserving numerical methods can be proved to provide consistent approximations. First stability estimates that follow directly from the discrete energy balance and discrete minimum entropy estimates are
\begin{eqnarray}
\label{first}
&&{\rm ess} \sup_{\tau \in (0,T)} \| \vr_n (\tau, \cdot) \|_{L^\gamma(\Td)} \leq c({E}_0, \underline{s}), \nonumber \\
&&{\rm ess} \sup_{\tau \in (0,T)} \| S_n (\tau, \cdot) \|_{L^\gamma(\Td)} \leq c({E}_0, \underline{s}), \nonumber \\
&&{\rm ess} \sup_{\tau \in (0,T)} \| \vm_n (\tau, \cdot) \|_{L^{\frac{2 \gamma}{\gamma + 1}}(\Td; \R^d)} \leq c({E}_0, \underline{s})
\end{eqnarray}
uniformly for $n \to \infty$.

However, in order to prove consistency of a numerical method,  we typically need additional stability information. A convenient and physically reasonable way is to assume uniform boundedness of the approximate density from below and of the approximate energy from above: \\
\emph{There exist constants $\underline{\vr}, \overline{E}$ such that
\begin{equation}
\label{hypo}
\vr_n(t) \geq \underline{\vr} > 0, \ \ E_n(t) \leq \overline{E} \quad \mbox{ for all } t \in [0,T]
\end{equation}
uniformly for $n \to \infty$. }

As shown in \cite{LY} hypothesis \eqref{hypo} implies uniform boundedness of $(\vr_n, \vm_n, E_n)_{n=1}^\infty$ in $L^\infty( (0,T) \times \Td; \R^{d+2}).$ In addition, discrete  entropy stability of a numerical method is crucial not only to  obtain consistent approximation of the entropy inequality \eqref{MV26}, but also to obtain  \emph{weak BV estimates} that control discrete gradients arising in the numerical diffusion terms that are a part of the consistency errors~\cite{book, FLM_LF}.

The following theorem characterizes a weak limit of consistent approximations, cf.~\cite[Chapter~5.1.3]{book}. As we will see, weak limits bring us \emph{beyond weak entropy solutions}.

\begin{theorem}[\textbf{Existence of a DW solution}]
\label{DW_thm}
Let the initial data $(\vr_{0,n}, \vm_{0,n}, E_{0,n})_{n=1}^\infty$ satisfy
\[
\vr_{0,n} \geq \underline{\vr} > 0,\ \ E_{0,n} - \frac{1}{2} \frac{ |\vc{m}_{0,n}|^2 }{\vr_{0,n}} > 0, \quad n=1,2,\dots
\]
Let
$( \vr_n, \vc{m}_n, S_n )_{n=1}^\infty $ be a consistent approximation of the Euler equations in the sense of Definition~\ref{consistent}. Further, let stability hypothesis~\eqref{hypo} hold, i.e.
$$
\vr_n(t) \geq \underline{\vr} > 0, \ \ E_n(t) \leq \overline{E} \qquad \mbox{ for all } t \in [0,T], \mbox{ uniformly for }\ \  n\to \infty.
$$
Then up to a subsequence, as the case may be, stable consistent approximation $(\vr_n, \vm_n, S_n)_{n=1}^\infty$ generates a DW solution $(\vr, \vm, S)$ in the sense of Definition~\ref{DE3}
\begin{equation} \label{MV38}
\begin{split}
\vr_n &\to \vr \ \mbox{weakly-(*) in}\ L^\infty((0,T) \times \Td),\\
S_n &\to S \ \mbox{weakly-(*) in}\ L^\infty((0,T) \times \Td),\\
\vm_n &\to \vm \ \mbox{weakly-(*) in}\ L^\infty((0,T) \times \Td; \R^d)\qquad \mbox{ for } n \to \infty.
\end{split}
\end{equation}
Moreover,
$$
E(\vr_n, \vm_n, S_n) \to \left< \mathcal{V}_{t,x}; E(\tvr, \tvm, \tvS) \right>
\mbox{weakly-(*) in}\ L^\infty((0,T) \times \Td),\  n \to \infty.
$$
\end{theorem}

\begin{remark}
We note that in general $\left< \mathcal{V}_{t,x}; E(\tvr, \tvm, \tvS) \right>  \geq E(\vr, \vm, S)$ since the energy is a (convex) nonlinear function of its weakly convergent arguments $(\vr, \vm, S).$ By the same token, even if the consistent approximation satisfies the discrete energy balance, i.e. \eqref{MV20} holds as equality and $e_n^1 = 0$ for all $n=1,2,\dots$, we still have in the limit only the energy inequality, cf.~\eqref{E21}
\begin{equation}
\intO{E(\vr,\vm,S)(t, \cdot)} + \int_{\Td} \D \mathfrak{E} (t) \leq  \intO{E(\vr_0, \vm_0, S_0)} \quad \mbox{ for } t \in [0,T],
\end{equation}
where $\mathfrak{E} \in L^\infty(0,T; \mathcal{M}^+(\Td))$ is an energy defect defined by
\begin{equation}
\label{en_defect}
\mathfrak{E} = \left<  \mathcal{V}; E(\tvr, \tvm, \tvS) \right> -    E(\vr,\vm,S).
\end{equation}
\end{remark}
Theorem~\ref{DW_thm} provides a  general strategy to analyse convergence of any structure-preserving numerical scheme. We refer to our recent works, where consistency, stability and convergence of the following well-known numerical methods have been studied: the first order finite volume Godunov method in \cite{LY}, the first order Lax-Friedrichs method in~\cite{FLM_LF}, the viscosity finite volume method in~\cite{FLM18_brenner}, higher order discontinuous Galerkin methods in \cite{po}, higher order residual distribution methods in~\cite{abgrall} and the second order finite element flux-corrected method in~\cite{kuzmin}.

Numerical methods mentioned above were analysed as semi-discrete approximation schemes, keeping time continuous. Convergence results directly generalize to fully discrete time implicit methods. In the case of explicit time discretization, that is typical for hyperbolic conservation laws, we need to take into account that the discrete entropy inequality may fail. To cure this problem, a special technique, the so-called time relaxation, can be used. It  controls  a time step such that discrete entropy inequality of time explicit structure-preserving numerical methods is preserved, cf.~\cite{Ranocha}.

\section{Weak versus strong convergence}
\label{sec:4}

As demonstrated in Theorem~\ref{DW_thm} the approximate solutions obtained by structure-preserving numerical methods may not, in general,  converge strongly. This is typical in turbulent flows, such as the Kelvin-Helmholtz or  the Richtmyer-Meshkov problem. However, even in such cases it would be desirable to recover strong convergence at least for some ``generalized objects''.
As shown in our recent works~\cite{FLM_K, FLSW} this is indeed possible. In order  to convert weakly converging sequences to strongly converging, we apply an averaging procedure mimicking the Strong Law of Large
Numbers in probability. The resulting concept of convergence was named $\mathcal{K}$-convergence,
according to the Koml\'os result~\cite{Kom}. Its functional analytic background
goes back to the classical results of Banach and Saks~\cite{BS}. Application to the Young measure $(\mathcal{V}_{t,x})_{(t,x) \in (0,T) \times \Td}$
relies on the Prokhorov theorem for random processes and allows to obtain compactness of the empirical measures.

Thus, we do not consider single numerical realisations, but observable quantities, such as the mean or deviation that are obtained by averaging over different mesh resolutions. The latter are also referred in the literature as the Ces\`aro averages over  different mesh resolutions.  To be consistent, we keep writing $(\vr_n, \vm_n, S_n)_{n=1}^\infty$ instead of more precise $(\vr_{h_n}, \vm_{h_n}, S_{h_n})_{h_n \searrow 0},$ where $h_n \in (0,1)$ is a mesh step, $h_n \to 0$ as $n\to \infty.$

\begin{theorem}[\textbf{$\mathcal{K}$-convergence}]\cite{FLSW}
\label{K}

Let the assumptions of Theorem~\ref{DW_thm} hold. Let $(\vr_n, \vm_n, S_n)_{n=1}^\infty$ be a stable, consistent approximation of the Euler equations \eqref{euler}, \eqref{entropy}.

Then up to a subsequence $(\vr_n, \vm_n, S_n)_{n=1}^\infty$ converges strongly to a DW solution $(\vr, \vm, S)$ in the following sense.
\begin{itemize}
\item {\bf{Strong convergences of Ces\` aro averages}}
\begin{equation} \label{WAS11}
\begin{split}
&\frac{1}{N} \sum_{k=1}^N \vr_{n_k} \to \vr \ \mbox{ as }\ N \to \infty \ \mbox{in} \ L^q((0,T) \times \Td)
\ \mbox{for any}\ 1 \leq q < \infty, \\
&\frac{1}{N} \sum_{k=1}^N \vm_{n_k} \to \vm \ \mbox{ as }\ N \to \infty \ \mbox{in} \ L^q((0,T) \times \Td; \R^d) \mbox{ for any }\ 1 \leq q < \infty, \\
&\frac{1}{N} \sum_{k=1}^N S_{n_k} \to S \ \mbox{ as }\ N \to \infty \ \mbox{in} \ L^q( (0,T) \times \Td)
\ \mbox{for any}\ 1 \leq q < \infty, \\
&\frac{1}{N} \sum_{k=1}^N E(\vr_{n_k}, \vm_{n_k}, S_{n_k}) \to \left< \mathcal{V}_{t,x}, E(\tvr, \tvm, \tvS) \right> \\
&\phantom{mmmmmmmm} \mbox{ as }\ N \to \infty \ \mbox{in} \ L^q( (0,T) \times \Td)
\ \mbox{for any}\ 1 \leq q < \infty;
\end{split}
\end{equation}

\item {\bf{Strong convergence  to Young measure in the Wasserstein metric}}
\begin{equation} \label{WAS12}
{W_r} \left[ \frac{1}{N} \sum_{k=1}^N \delta_{[\vr_{n_k},
\vm_{n_k}, S_{n_k} ]} ; \mathcal{V}_{t,x} \right] \to 0
\ \mbox{as}\ N \to \infty \ \ \mbox{in}\ L^{s}((0,T) \times \Td)
\end{equation}
for any $1 \leq s < r < \infty$.
Here $W_r$ denotes the Wasserstein metric of order $r.$
\end{itemize}

\end{theorem}

We close this section by mentioning that the weak convergence of stable consistent approximations (and of numerical solutions generated by structure-preserving methods) directly turns to the strong convergence in the following situations.

\begin{itemize}
\item Let the hypothesis of Theorem~\ref{DW_thm} hold. Suppose that the Euler system \eqref{euler}, \eqref{entropy}
admits a strong (Lipshitz-continuous) solution $(\vr, \vm, S)$ defined on $[0,T]$.

Then
\begin{eqnarray*}
&&\vr_n \to \vr, \ \vm_n \to \vm, \ S_n \to S \ \mbox{ and } E(\vr_n, \vm_n, S_n) \to  E(\vr, \vm, S)
\\
&& \phantom{mmmmmmm}  \mbox{ strongly in }
L^q((0,T) \times \Td), \ 1 \leq q < \infty, \mbox{ for } n \to \infty.
\end{eqnarray*}
In this case the convergence rate of a numerical approximation $(\vr_n, \vm_n, S_n)$ can be computed by means of the relative energy
$E(\vr_n, \vm_n, S_n\Big | \vr, \vm, S),$ see \cite{Yuhuan_error}.

\item Let the hypotheses of Theorem~\ref{DW_thm} hold. Let the weak limit $(\vr, \vm, S)$ in \eqref{MV38}
be a weak solution of the Euler system \eqref{euler}, \eqref{entropy}.

Then
$$
\mathcal{V}_{t,x} = \delta_{[\vr(t,x), \vm(t,x), S(t,x)]} \ \mbox{ for a.a. } (t,x) \in (0,T) \times \Td
$$
and up to a subsequence
\begin{eqnarray*}
&&\vr_n \to \vr, \ \vm_n \to \vm, \ S_n \to S \ \mbox{ and } E(\vr_n, \vm_n, S_n) \to E(\vr, \vm, S)
\\
&&\phantom{mmmmmm} \mbox{ strongly in }  L^q((0,T) \times \Td), \ 1 \leq q < \infty, \ \mbox{ for } n\to \infty.
\end{eqnarray*}
\end{itemize}

\section{Numerical simulations}
\label{sec:5}
The aim of this section is to illustrate theoretical results presented in the previous sections by means of two finite volume methods: the viscous finite volume (VFV) method, studied in \cite{FLM18_brenner}, and the second order generalized Riemann problem (GRP) method, proposed in \cite{ben2006,ben2007}.
To this end, let us consider the well-known Kelvin-Helmholtz problem \cite{Helmhotz, Kelvin, FLSW} on $\mathbb{T}^2$ with the initial data
\[
(\vr, u_1, u_2,p)(x)=\left\{
\begin{array}{ll}
(2, -0.5, 0, 2.5), & \text{if}\ I_1<x_2<I_2,\\
(1, 0.5, 0, 2.5), & \text{otherwise,}
\end{array}
\right.
\]
where $I_1$ and $I_2$ are the two perturbed interface profiles given by
\[
I_j = I_j(x,\omega):=J_j+\epsilon Y_j(x,\omega), \quad j=1,2,
\]
 where  $J_1=0.25$ and $J_2=0.75$ and
\[
Y_j(x,\omega)=\sum_{m=1}^{M}a_j^m(\omega)\cos(b_j^m(\omega)+2m\pi x_1),\quad j=1,2.
\]
Here $a_j^m=a_j^m(\omega)\in[0,1]$ and $b_j^m=b_j^m(\omega)\in[-\pi,\pi]$, $j=1,2$, $m=1,\ldots,M$ are (fixed) random numbers as in our previous paper \cite{FLSW}. The coefficients $a_j^m$ have been normalized such that $\sum_{m=1}^{M}a_j^m=1$ to guarantee that $|I_j(x,\omega)-J_j|\leq \epsilon$ for $j=1,2$. We have set $M=10$ and $\epsilon=0.01.$


We shall work with the numerical solutions $\{\vr_n, \vm_n, S_n\}_{n \nearrow \infty}$ computed on a mesh with $n_x \times n_x$ regular cells,  $n_x=2^{5+n}$,  $n=1,\dots,6$, and their Ces\`aro averages $\{\widetilde{\vr}_N, \widetilde{\vm}_N, \widetilde{S}_N\}_{N \nearrow \infty}$ given by
\[\widetilde{\vr}_N = \frac1N \sum_{n=1}^N\vr_n, \; \; \; \widetilde{\vm}_N = \frac1N \sum_{n=1}^N \vm_n,\; \; \; \widetilde{S}_N = \frac1N \sum_{n=1}^N S_n.\]
Following Theorems~\ref{DW_thm} and \ref{K} we can approximate the Reynolds stress defect $\mathfrak{R}_N$ and the energy defect $\mathfrak{E}_N$ by
\begin{equation*}
\begin{aligned}
\mathfrak{R}_N &= \frac1N \sum_{n=1}^N \left( \frac{\vm_n \otimes \vm_n}{\vr_n} + p(\vr_n, S_n) \mathbb{I}  \right)
-   \left( \frac{\widetilde{\vm}_N \otimes \widetilde{\vm}_N}{\widetilde{\vr}_N} + p(\widetilde{\vr}_N, \widetilde{S}_N)\mathbb{I}  \right),\\
\mathfrak{E}_N &=  \frac1N \sum_{n=1}^N \left( \frac{|\vm_n|^2}{2\vr_n} + \vr_n e(\vr_n, S_n)  \right) -  \left( \frac{|\widetilde{\vm}_N|^2}{2\widetilde{\vr}_N} + \widetilde{\vr}_N e(\widetilde{\vr}_N, \widetilde{S}_N)  \right).
\end{aligned}
\end{equation*}

\begin{figure}[h]
\centering
\begin{subfigure}{0.245\linewidth} \centering
\includegraphics[width=1.1\textwidth]{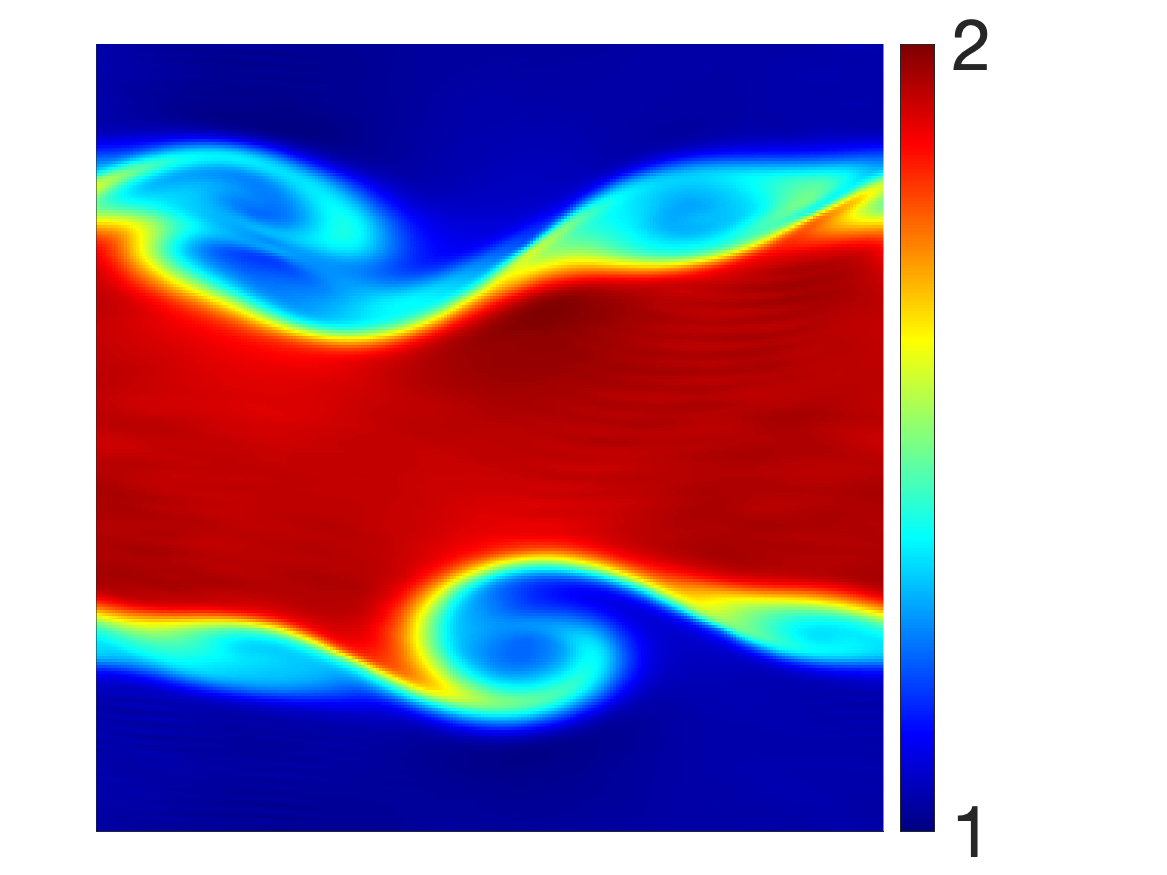}
\end{subfigure}
\begin{subfigure}{0.245\linewidth} \centering
\includegraphics[width=1.1\textwidth]{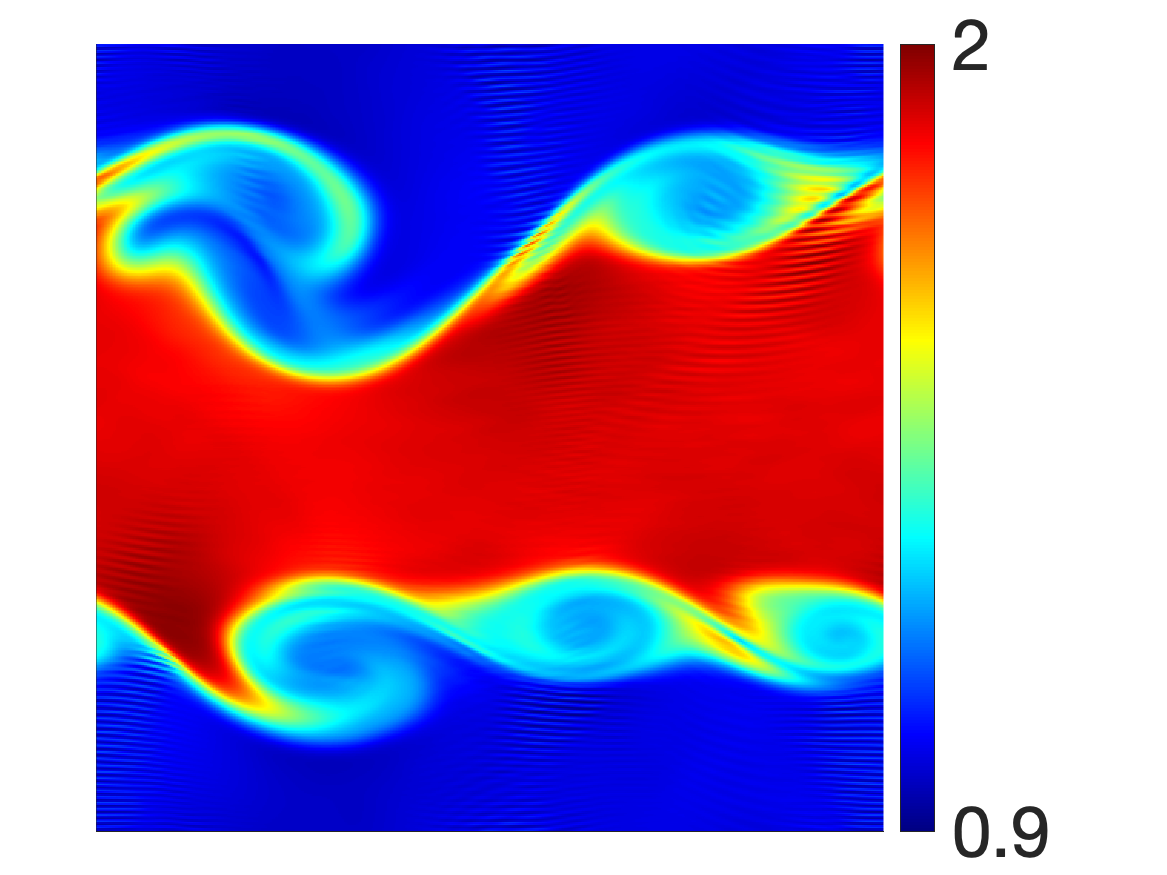}
\end{subfigure}
\begin{subfigure}{0.245\linewidth} \centering
\includegraphics[width=1.1\textwidth]{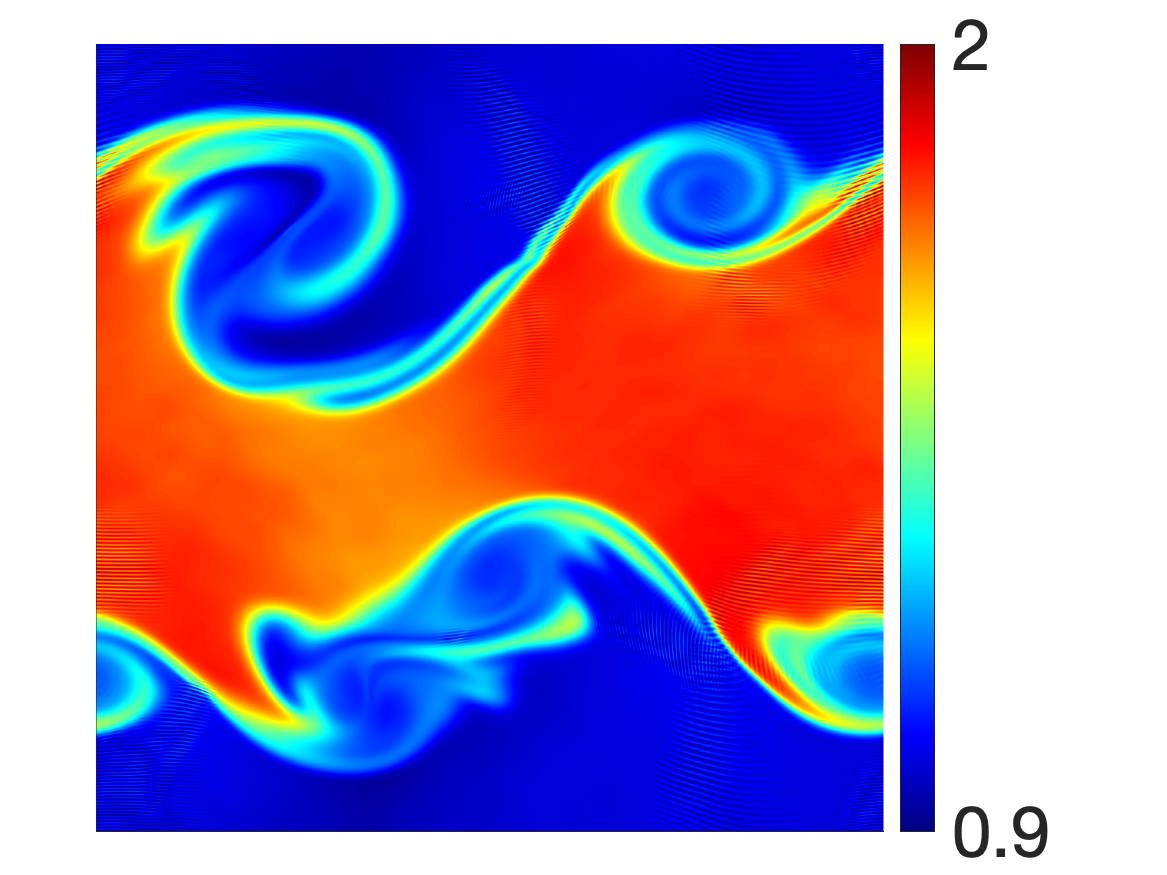}
\end{subfigure}
\begin{subfigure}{0.245\linewidth} \centering
\includegraphics[width=1.1\textwidth]{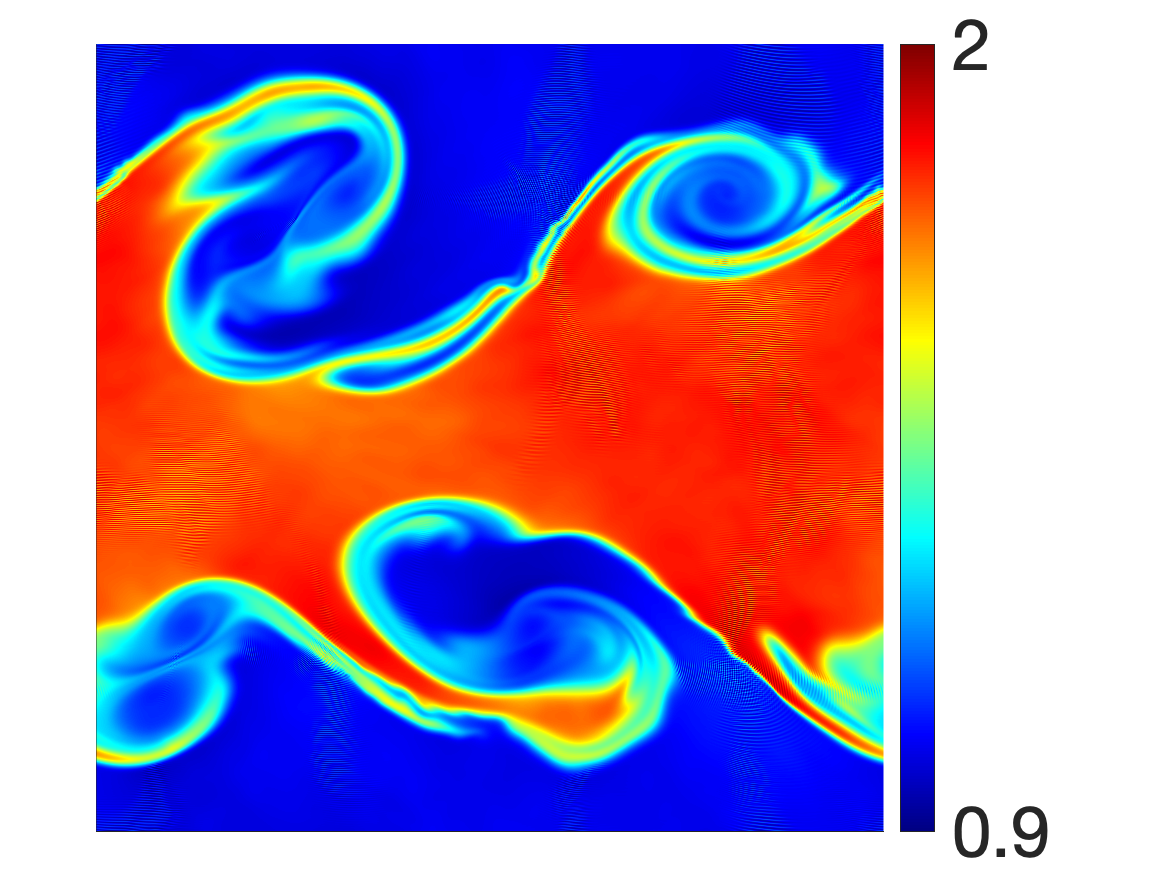}
\end{subfigure}
\\
\begin{subfigure}{0.245\linewidth} \centering
\includegraphics[width=1.1\textwidth]{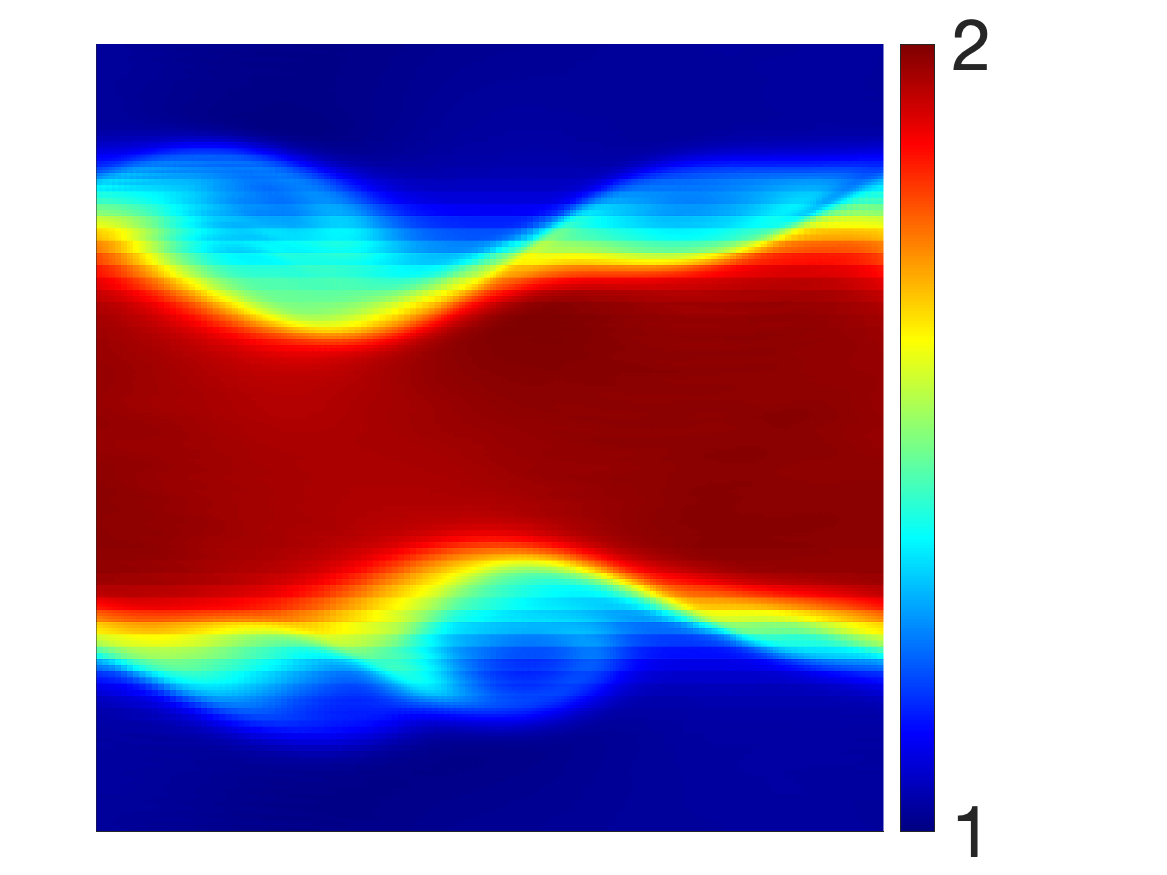}
\end{subfigure}
\begin{subfigure}{0.245\linewidth} \centering
\includegraphics[width=1.1\textwidth]{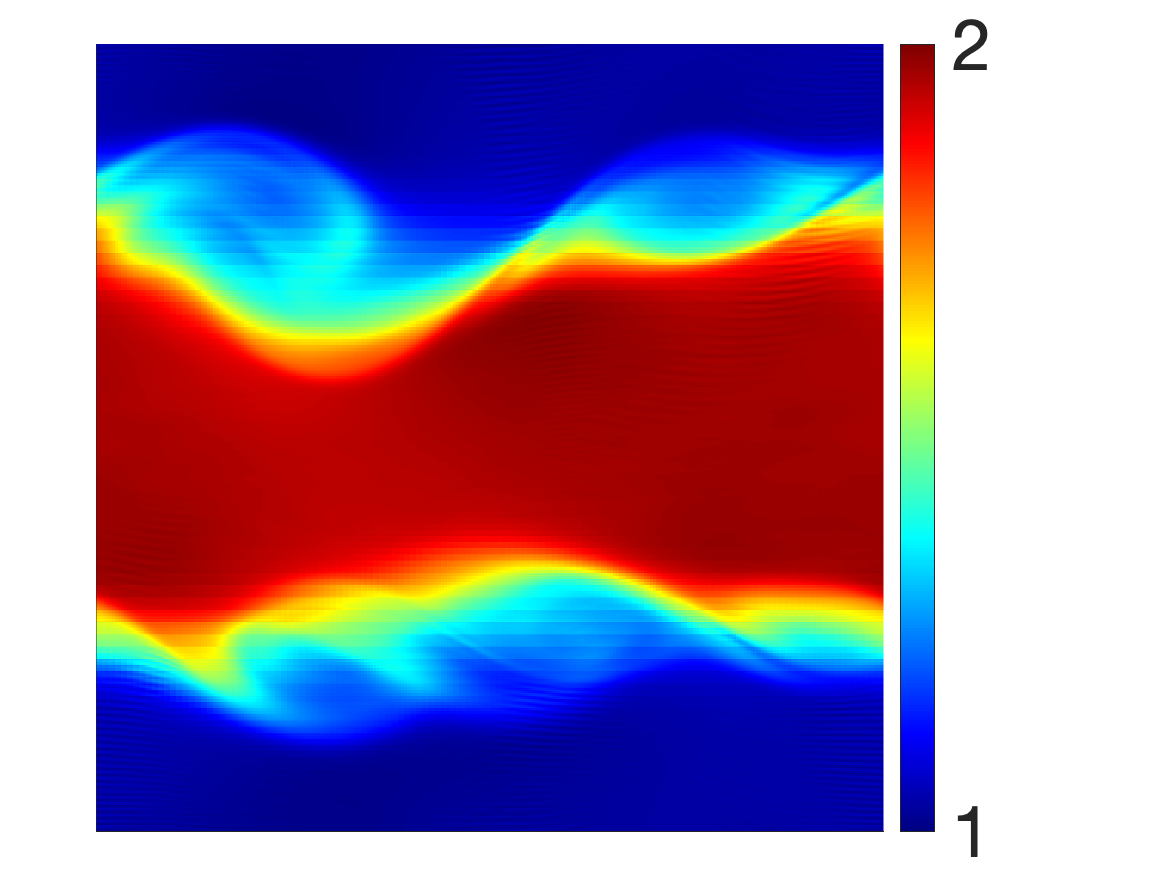}
\end{subfigure}
\begin{subfigure}{0.245\linewidth} \centering
\includegraphics[width=1.1\textwidth]{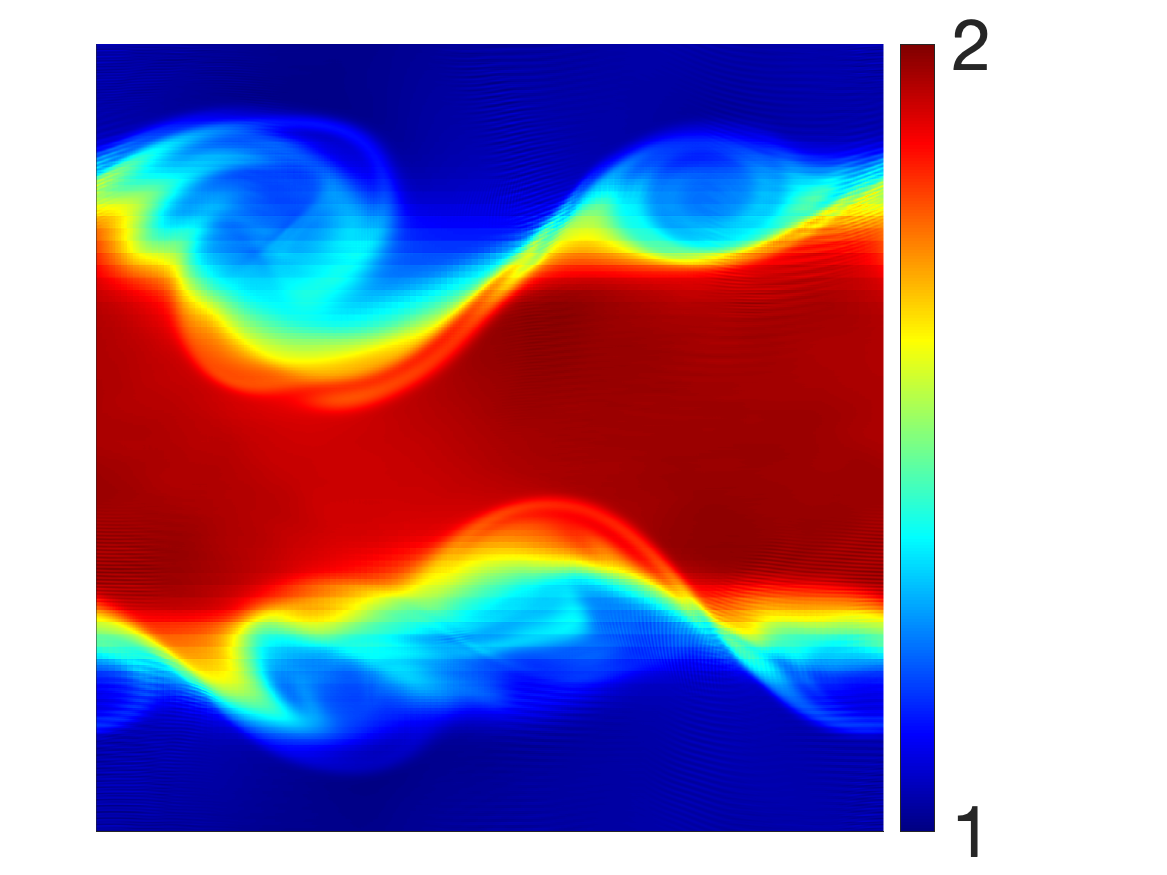}
\end{subfigure}
\begin{subfigure}{0.245\linewidth} \centering
\includegraphics[width=1.1\textwidth]{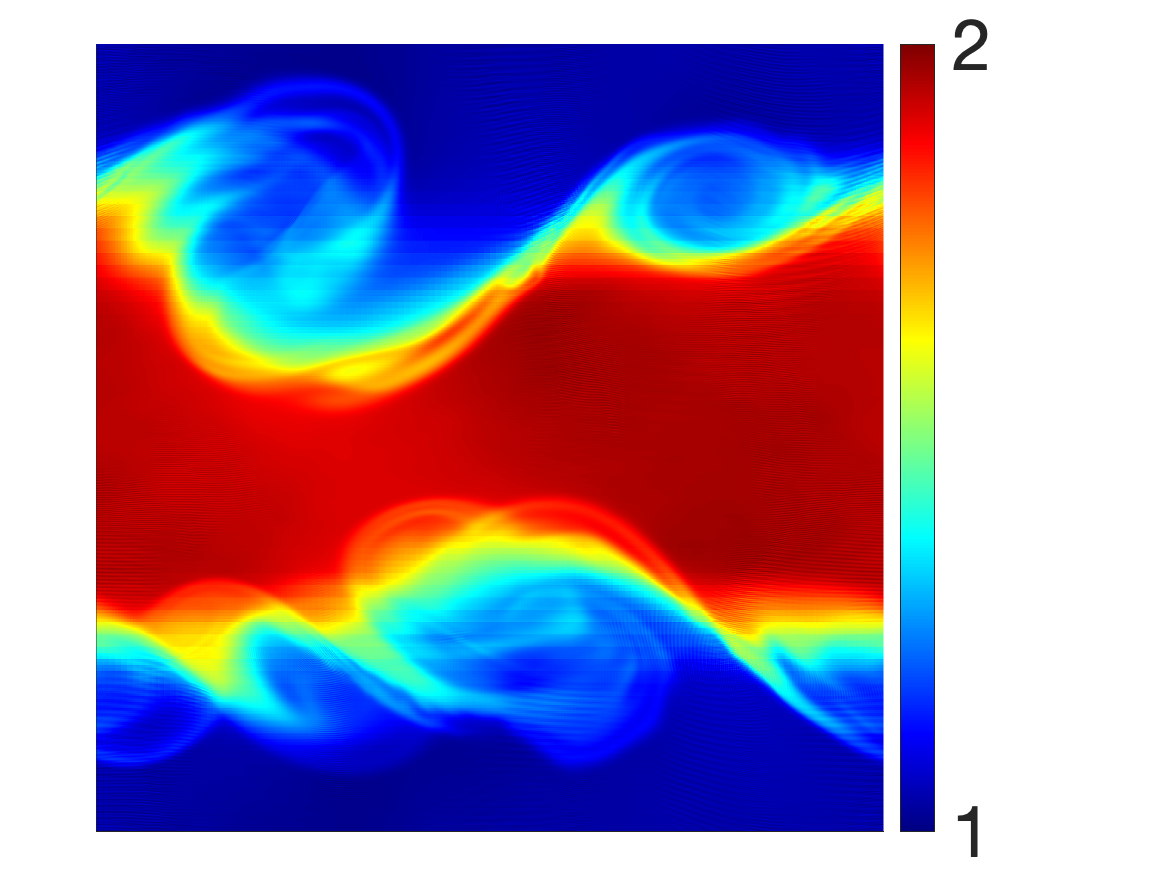}
\end{subfigure}
\caption{Density $\vr_n, n=3,\dots,6$ (top) and the Ces\`aro average of the density $\widetilde{\vr}_N,N=3,\dots,6$ (bottom) computed by the VFV scheme.}
\label{fig_sols_vfv-1}
\end{figure}

\begin{figure}[h]
\centering
\begin{subfigure}{0.245\linewidth} \centering
\includegraphics[width=1.1\textwidth]{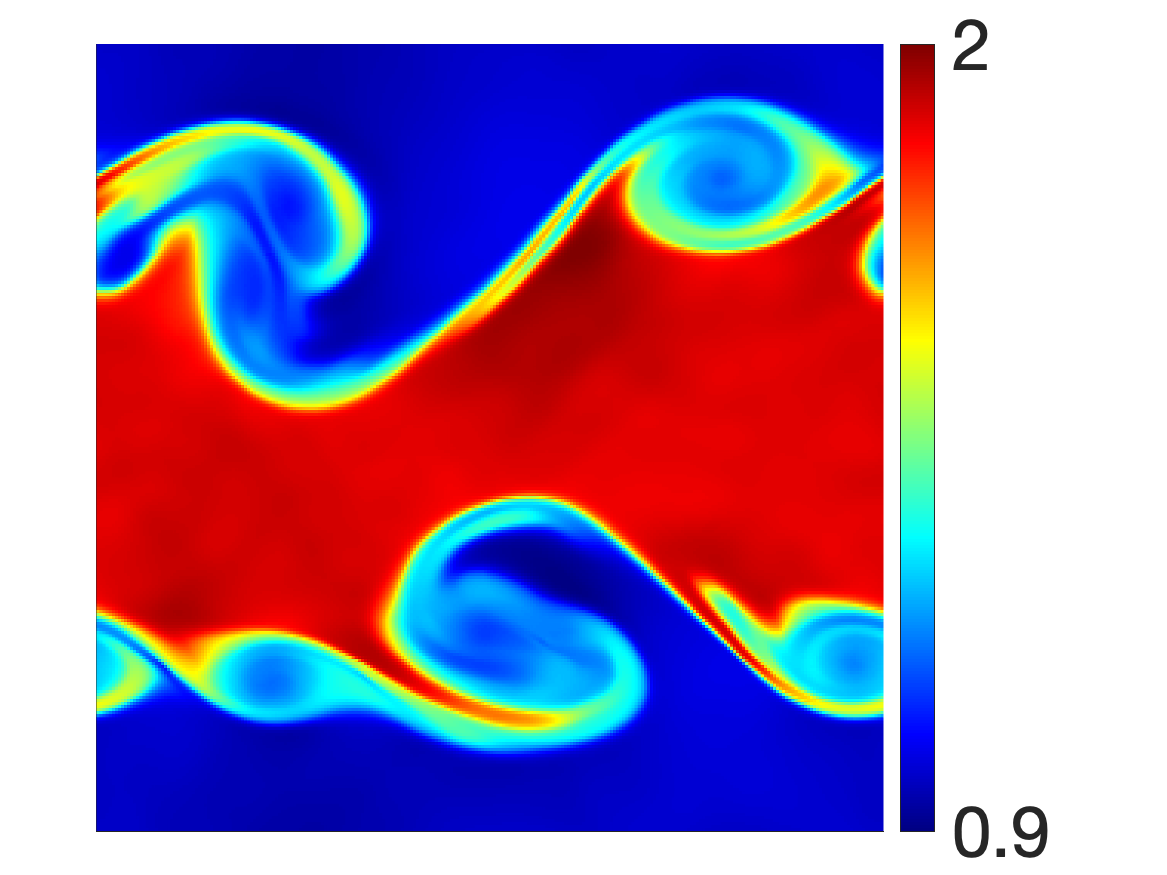}
\end{subfigure}
\begin{subfigure}{0.245\linewidth} \centering
\includegraphics[width=1.1\textwidth]{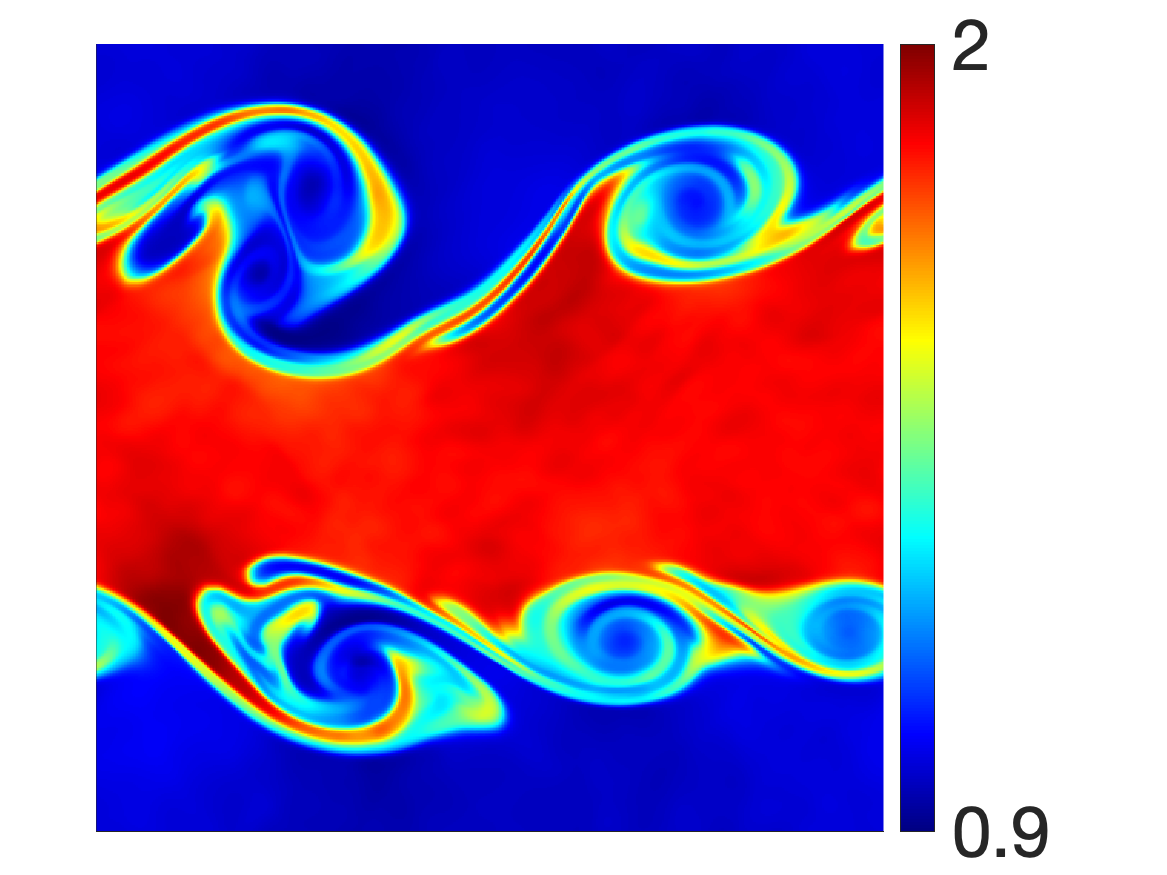}
\end{subfigure}
\begin{subfigure}{0.245\linewidth} \centering
\includegraphics[width=1.1\textwidth]{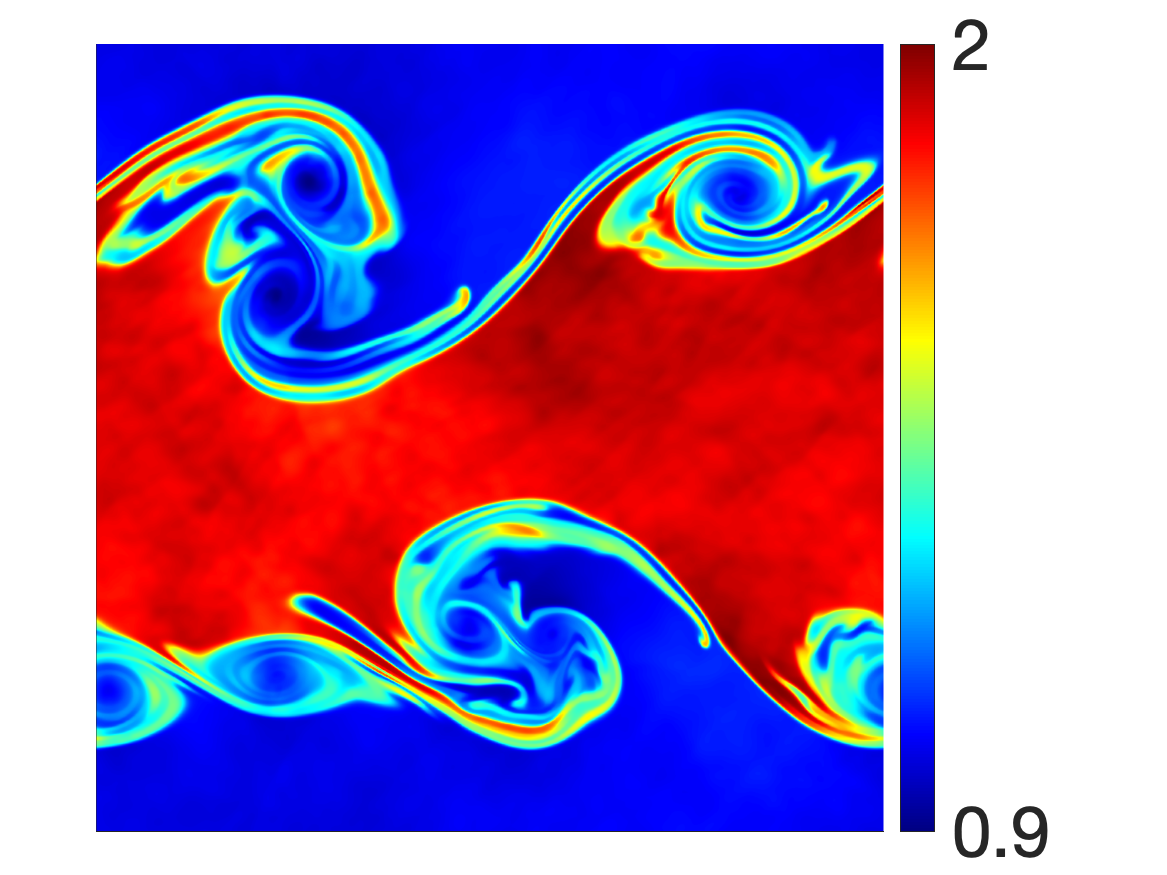}
\end{subfigure}
\begin{subfigure}{0.245\linewidth} \centering
\includegraphics[width=1.1\textwidth]{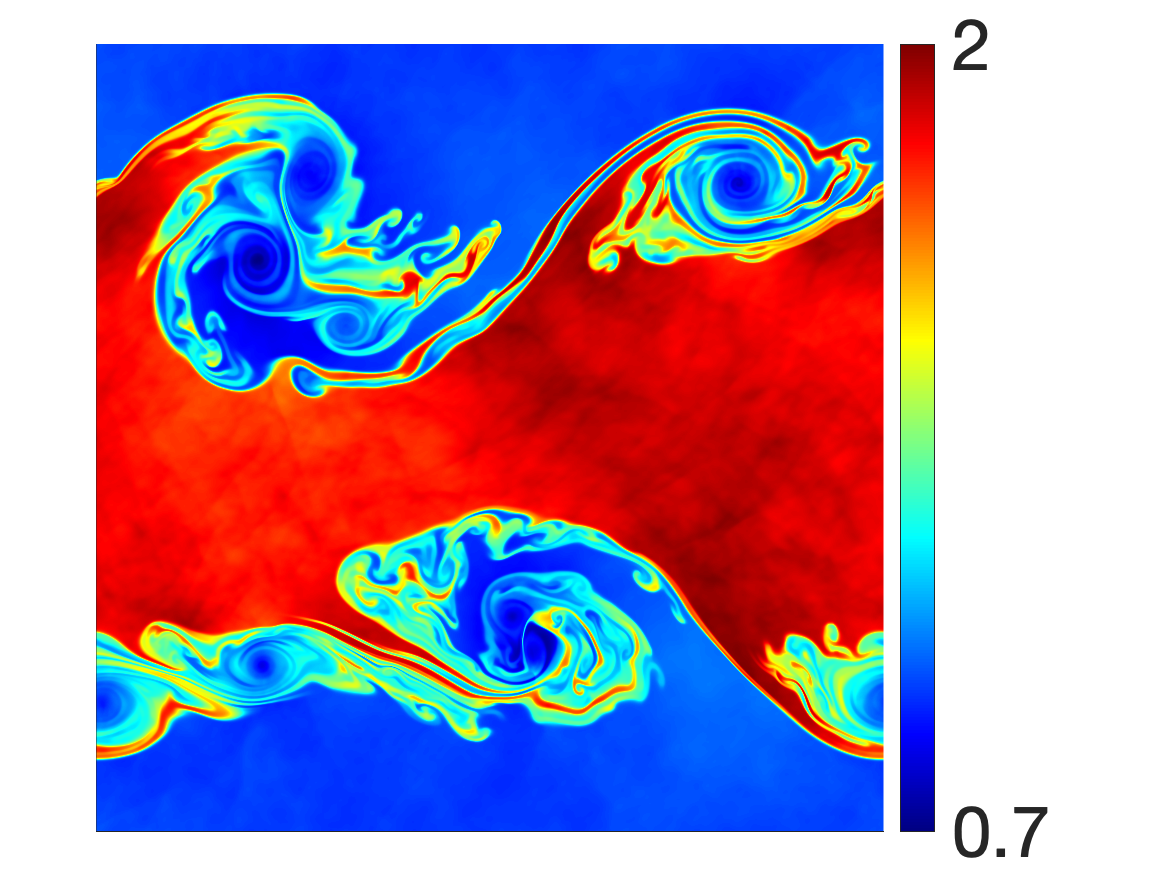}
\end{subfigure}
\\
\begin{subfigure}{0.245\linewidth} \centering
\includegraphics[width=1.1\textwidth]{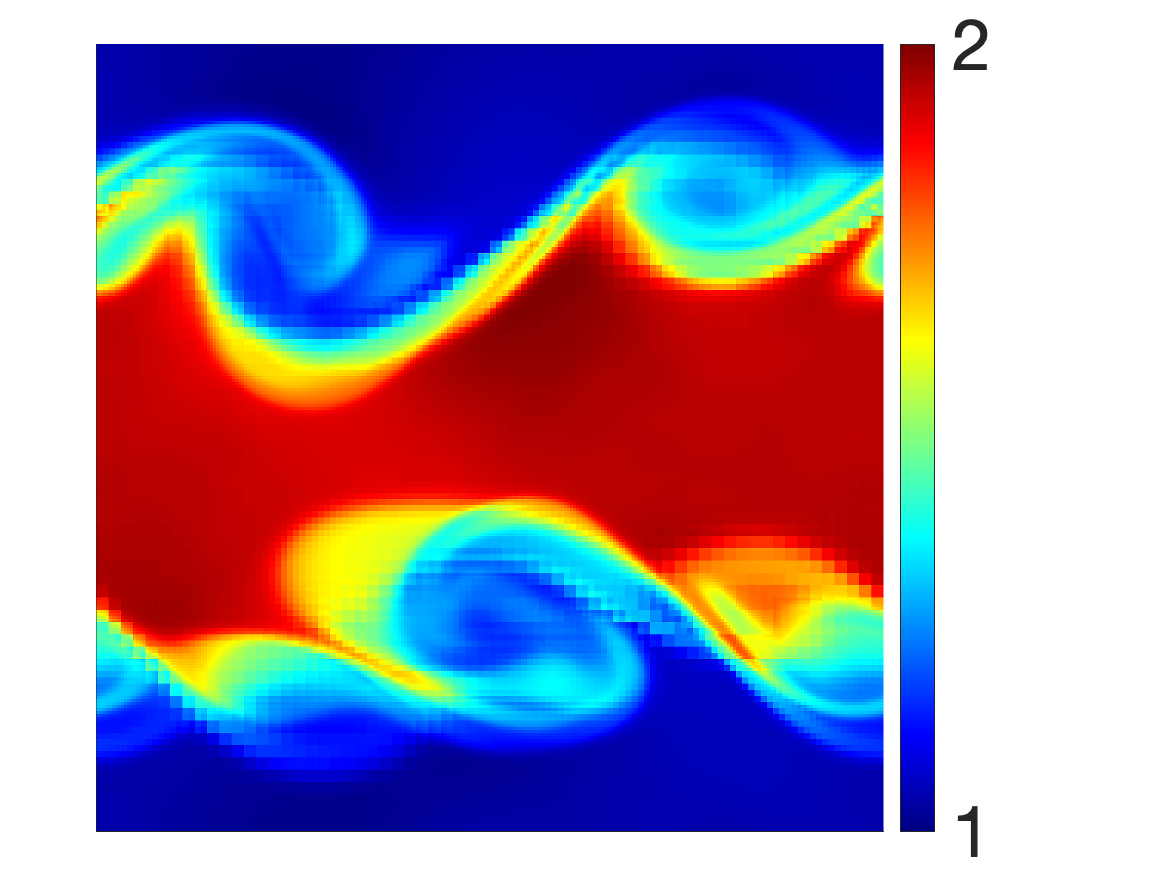}
\end{subfigure}
\begin{subfigure}{0.245\linewidth} \centering
\includegraphics[width=1.1\textwidth]{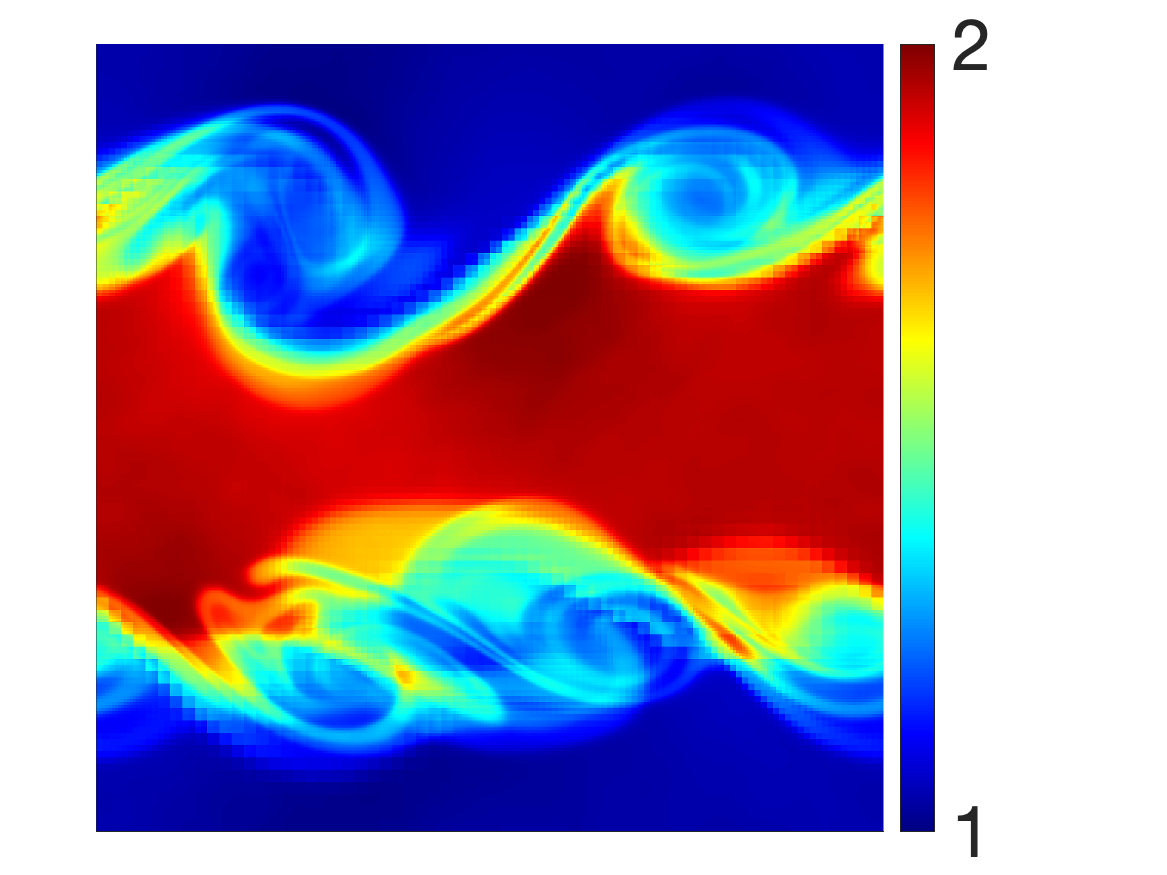}
\end{subfigure}
\begin{subfigure}{0.245\linewidth} \centering
\includegraphics[width=1.1\textwidth]{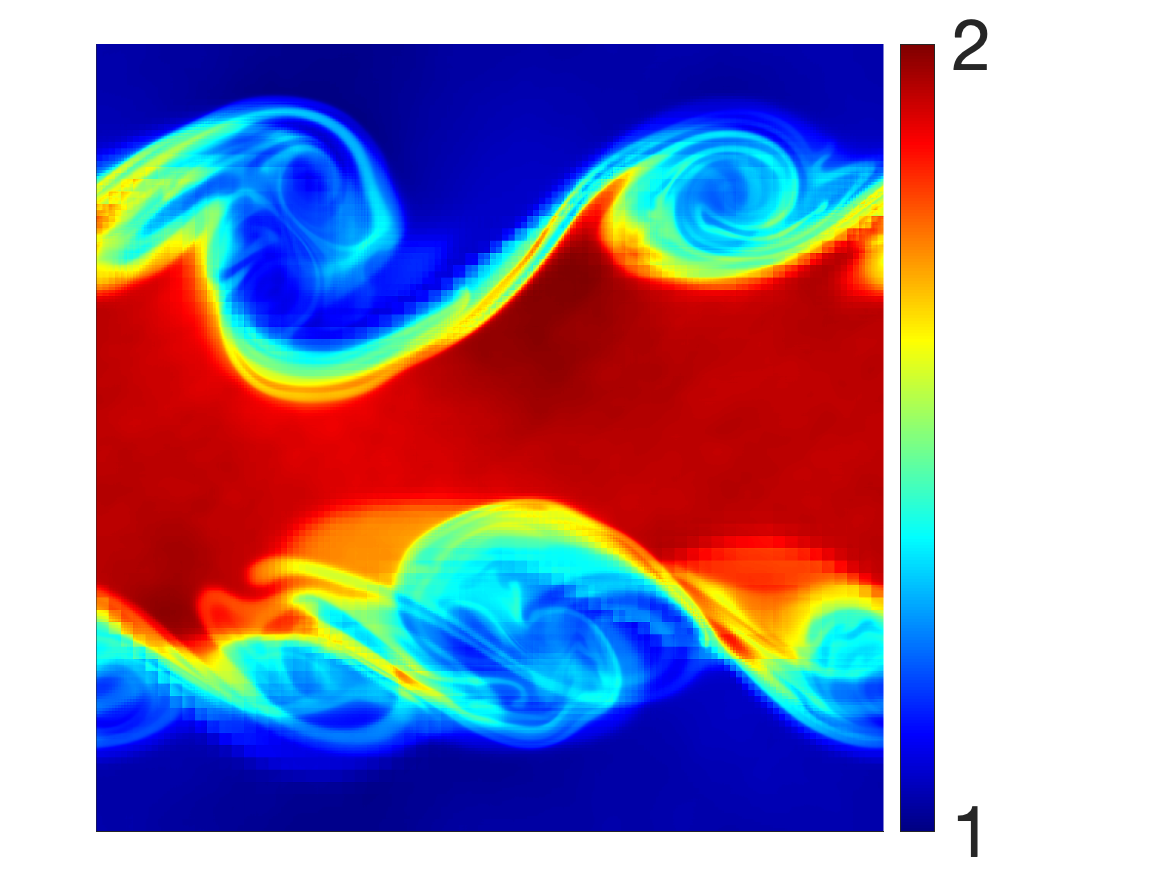}
\end{subfigure}
\begin{subfigure}{0.245\linewidth} \centering
\includegraphics[width=1.1\textwidth]{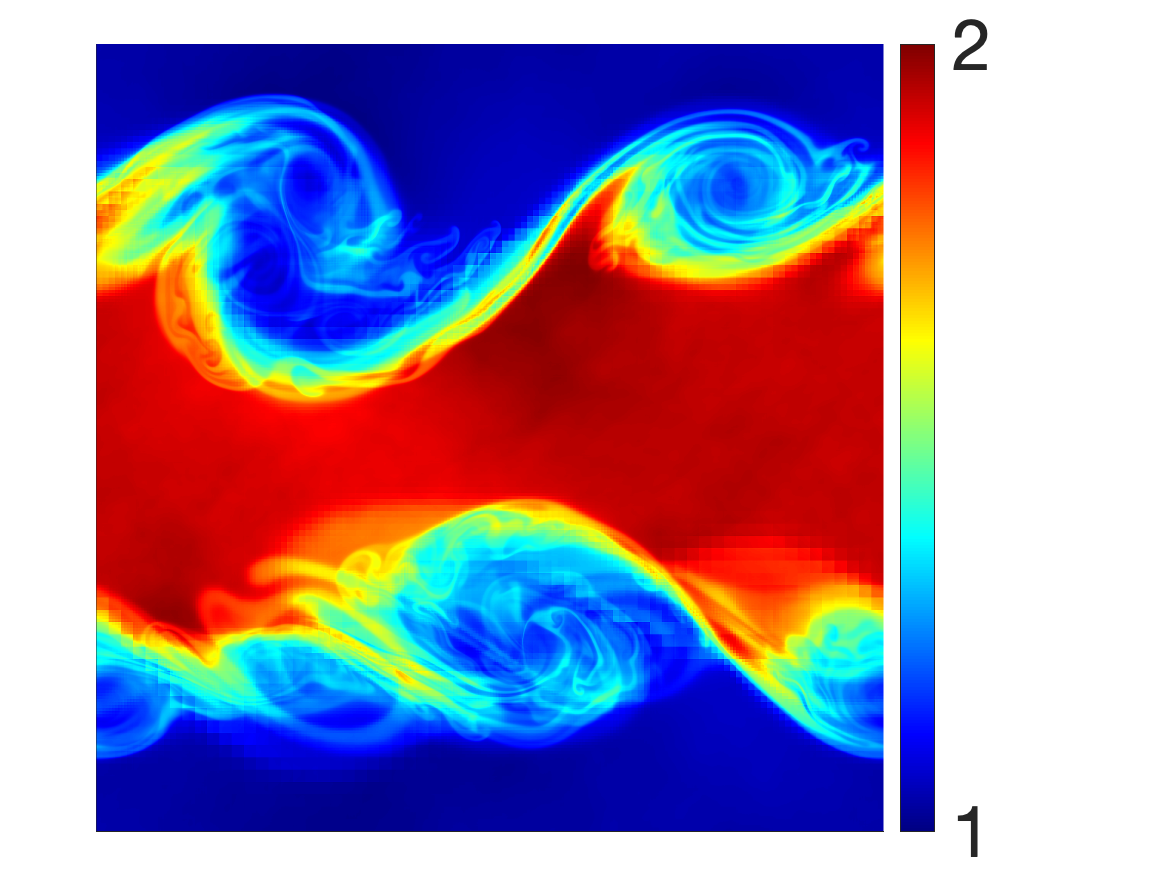}
\end{subfigure}
\caption{Density $\vr_n, n=3,\dots,6$ (top) and the Ces\`aro average of the density $\widetilde{\vr}_N, N=3,\dots,6$ (bottom) computed by the GRP scheme.}
\label{fig_sols_yue}
\end{figure}

\begin{figure}[h!]
\centering
\begin{subfigure}{0.48\textwidth}
\includegraphics[width=\textwidth]{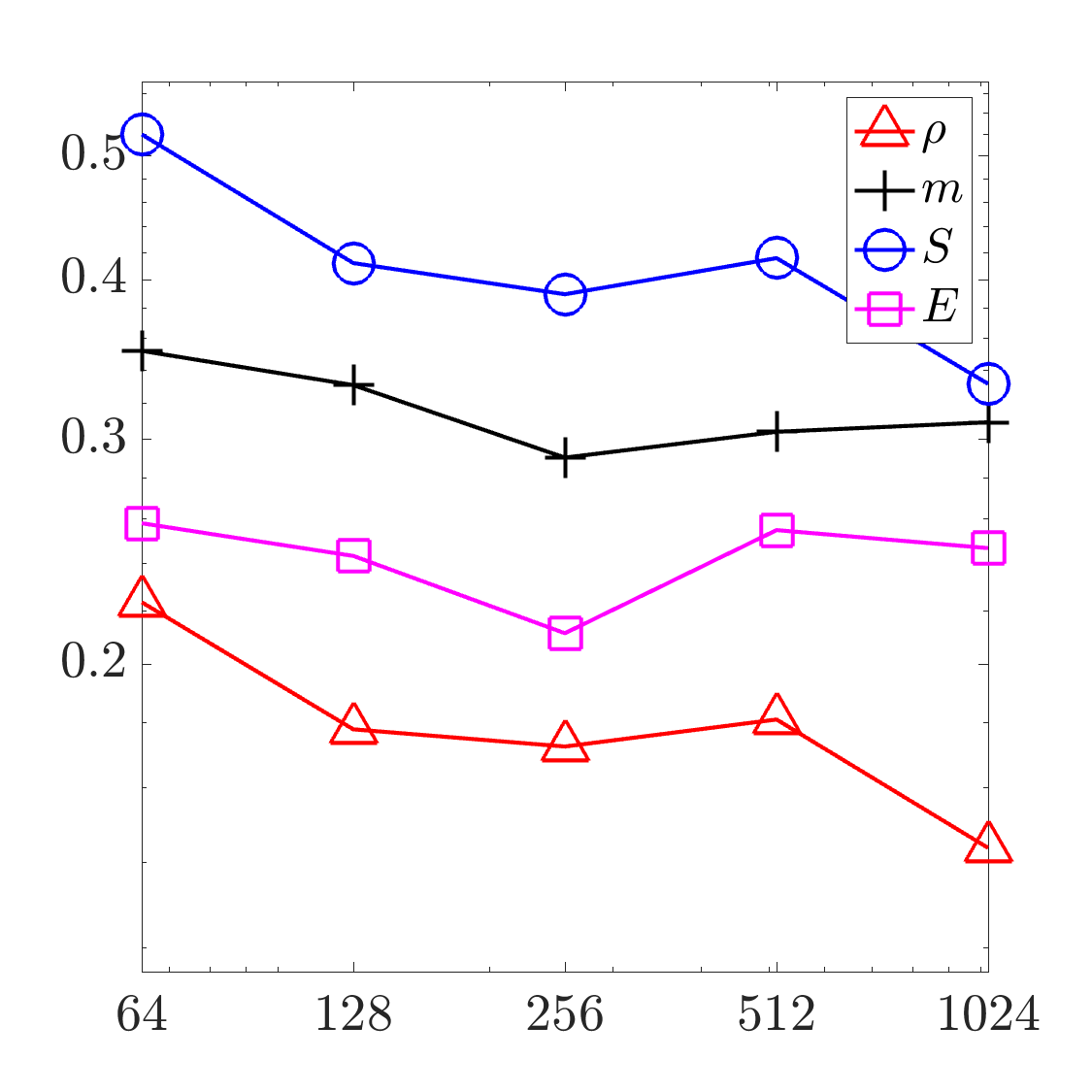}
\caption{VFV}
\end{subfigure}
\begin{subfigure}{0.48\textwidth}
\includegraphics[width=\textwidth]{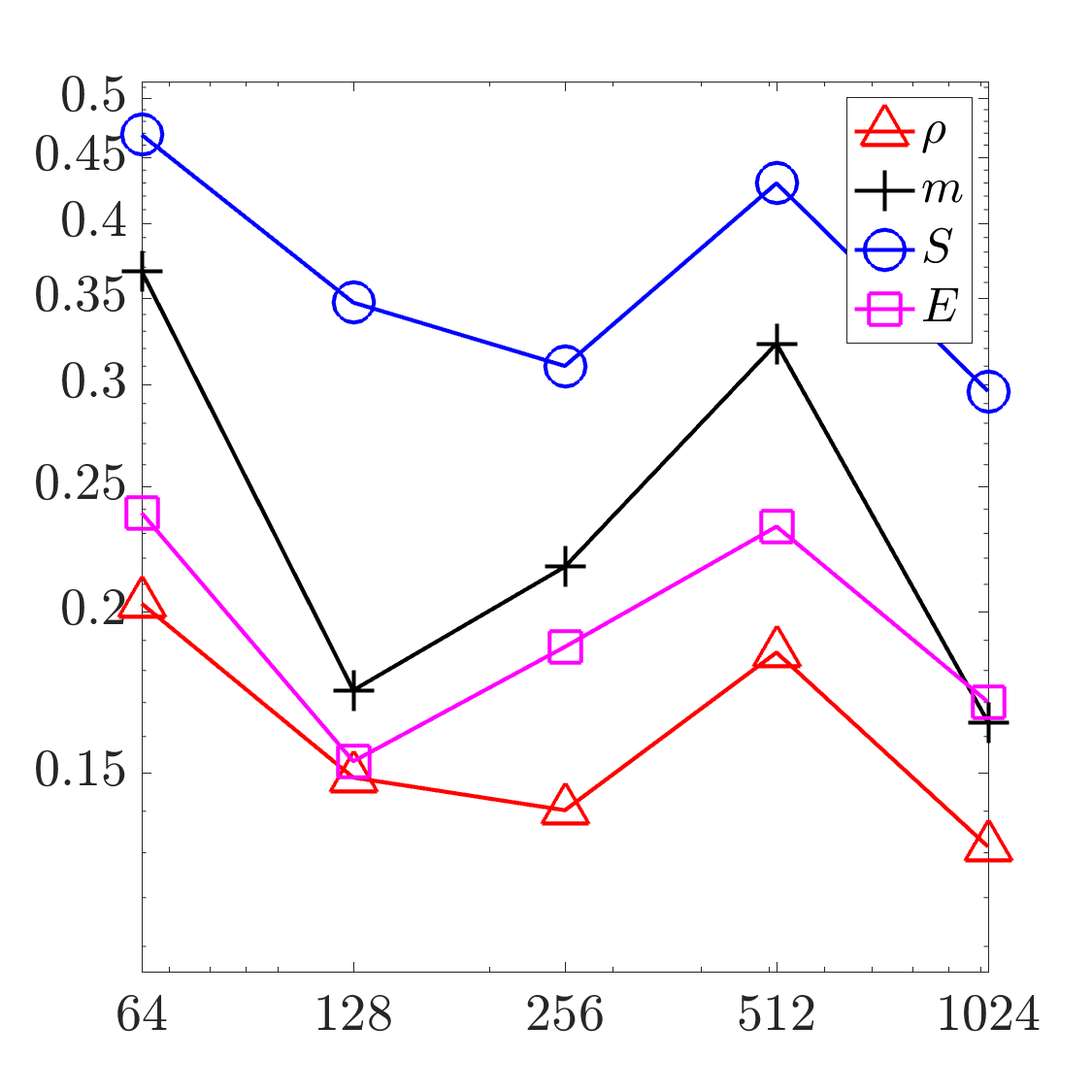}
\caption{GRP scheme}
\end{subfigure}
\caption{$L^1$-errors for ${\vr}_n$, ${\vm}_n$, ${S}_n$,  and ${E}_n$, computed for different mesh resolutions.}
\label{fig_err_ori}
\end{figure}

\begin{figure}[h!]
\centering
\begin{subfigure}{0.48\textwidth}
\includegraphics[width=\textwidth]{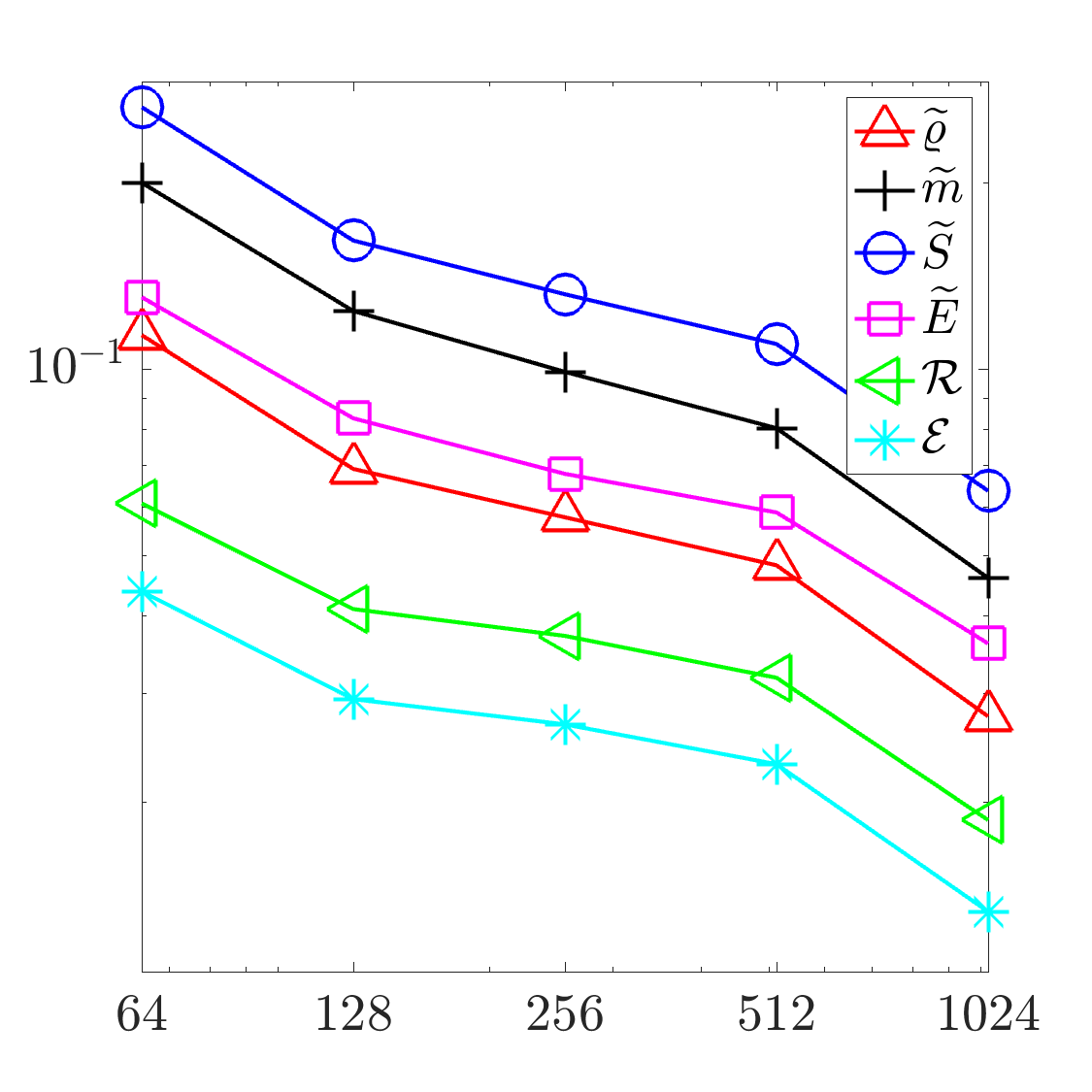}
\caption{VFV}
\end{subfigure}
\begin{subfigure}{0.48\textwidth}
\includegraphics[width=\textwidth]{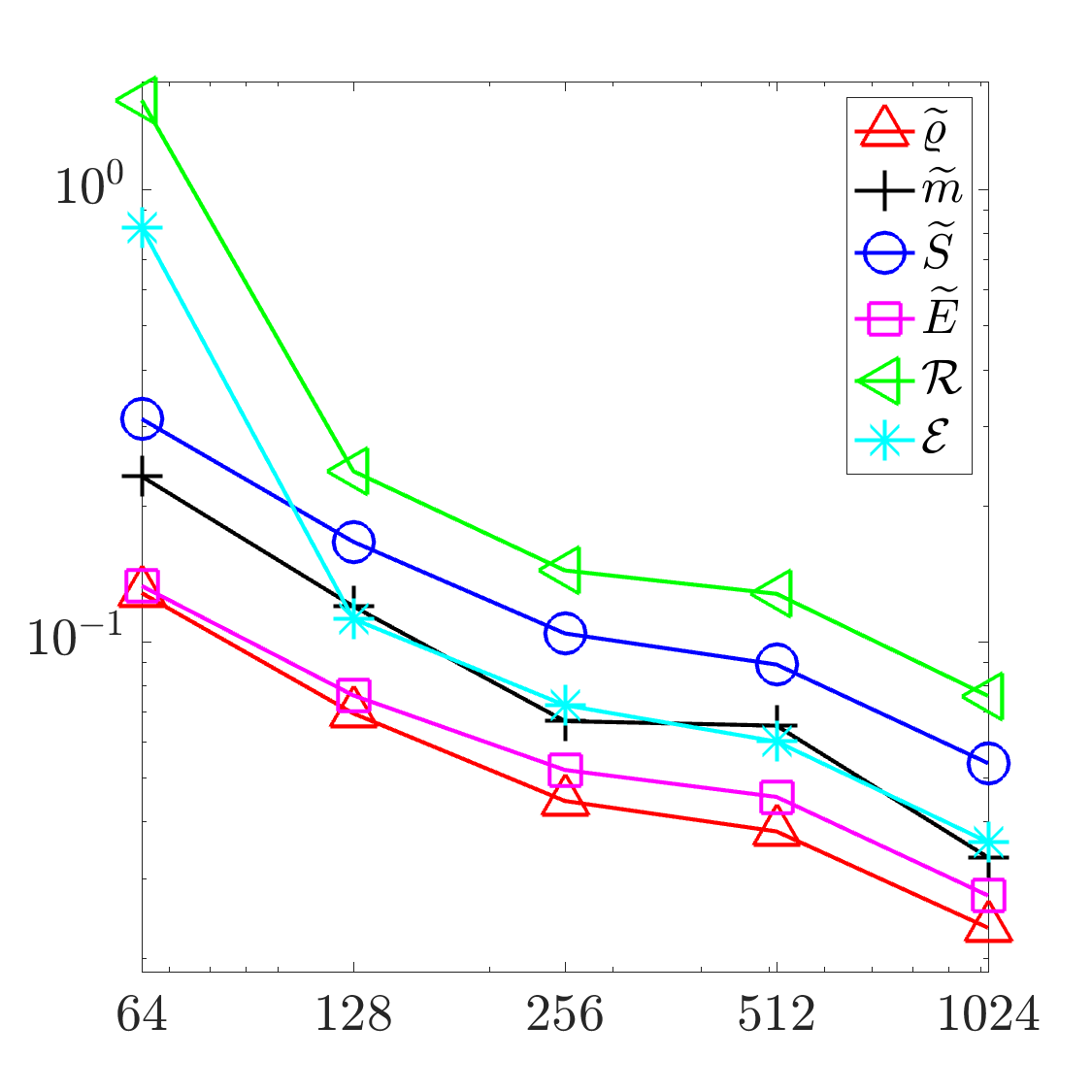}
\caption{GRP scheme}
\end{subfigure}
\caption{$L^1$-errors for $\widetilde{\vr}_N$, $\widetilde{\vm}_N$, $\widetilde{S}_N$,  $\widetilde{E}_N$,  $\mathfrak{R}_N$ (labeled as ``$\mathcal{R}$" in the plots), and  $\mathfrak{E}_{N}$ (labeled as ``$\mathcal{E}$" in the plots), computed for different mesh resolutions. }
\label{fig_err}
\end{figure}

In Figures~\ref{fig_sols_vfv-1} and \ref{fig_sols_yue} we present the numerical densities $\vr_n$, $n=3,\dots,6$, and their Ces\`aro averages $\widetilde{\vr}_N$, $N=3,\dots,6$, obtained by  the VFV method and the GRP method, respectively.
We can clearly observe fine vortex structures arising for refined meshes.
Fig.~\ref{fig_err_ori} illustrates the weak convergence of numerical solutions with respect to mesh refinement. Indeed, errors computed in the $L^1$-norm only oscillate and do not converge. In agreement with Theorem~\ref{K}, their Ces\`aro averages do converge, see Fig.~\ref{fig_err}.

\begin{figure}[h]
\centering
\begin{subfigure}{0.245\linewidth} \centering
\includegraphics[width=1\textwidth]{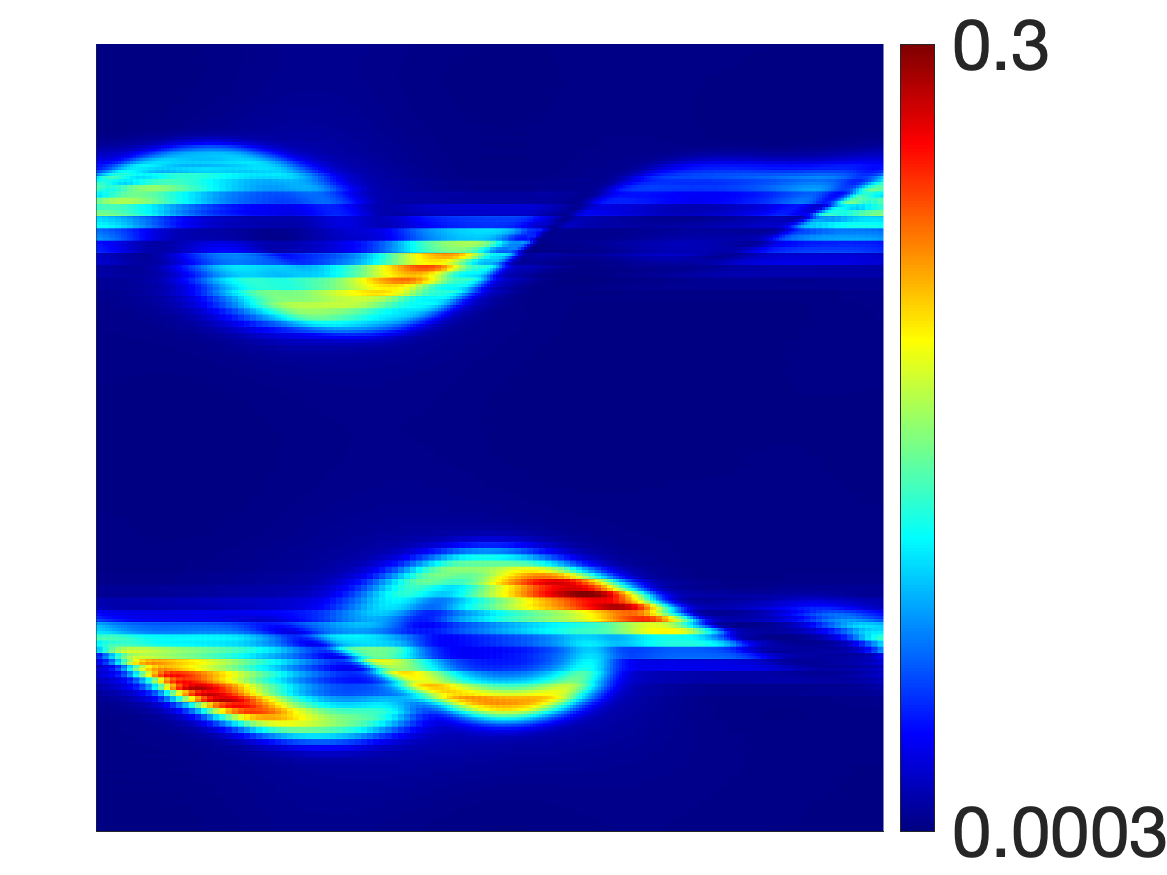}
\end{subfigure}
\begin{subfigure}{0.245\linewidth} \centering
\includegraphics[width=1\textwidth]{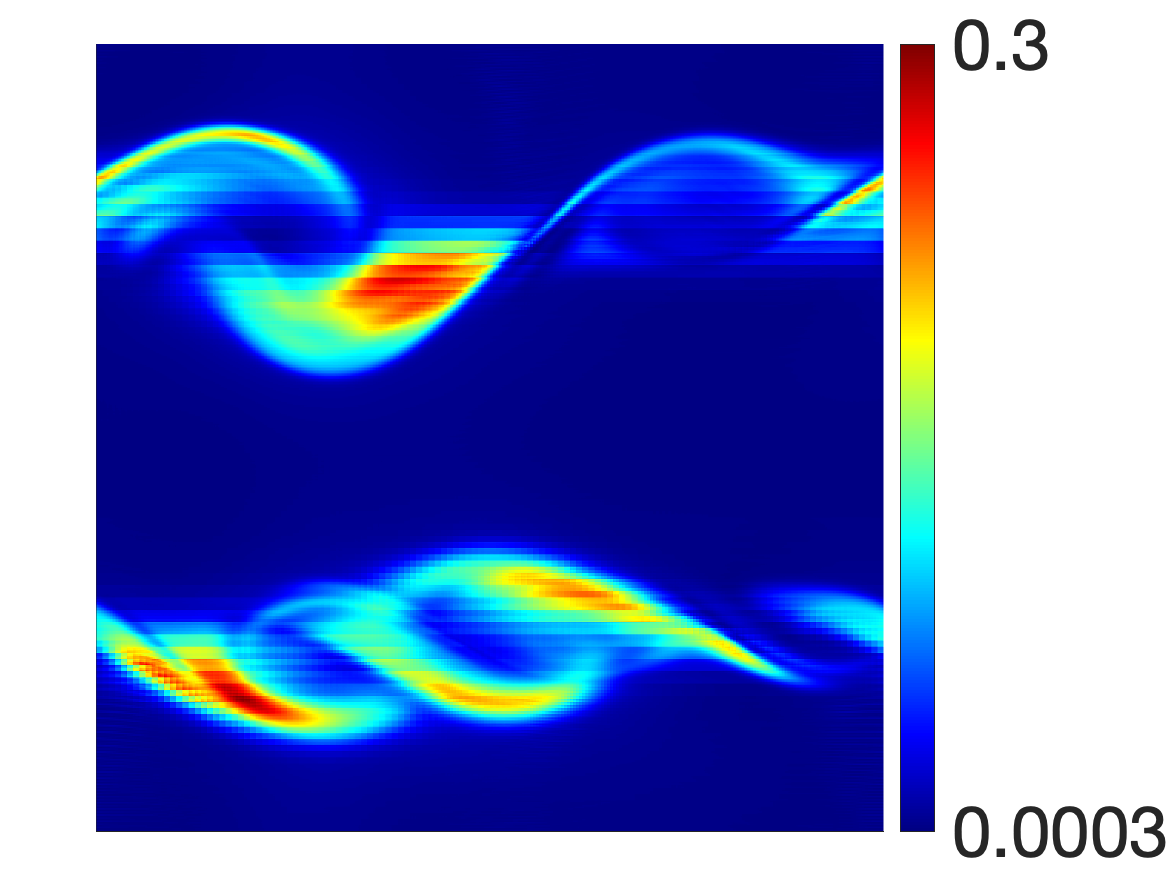}
\end{subfigure}
\begin{subfigure}{0.245\linewidth} \centering
\includegraphics[width=1\textwidth]{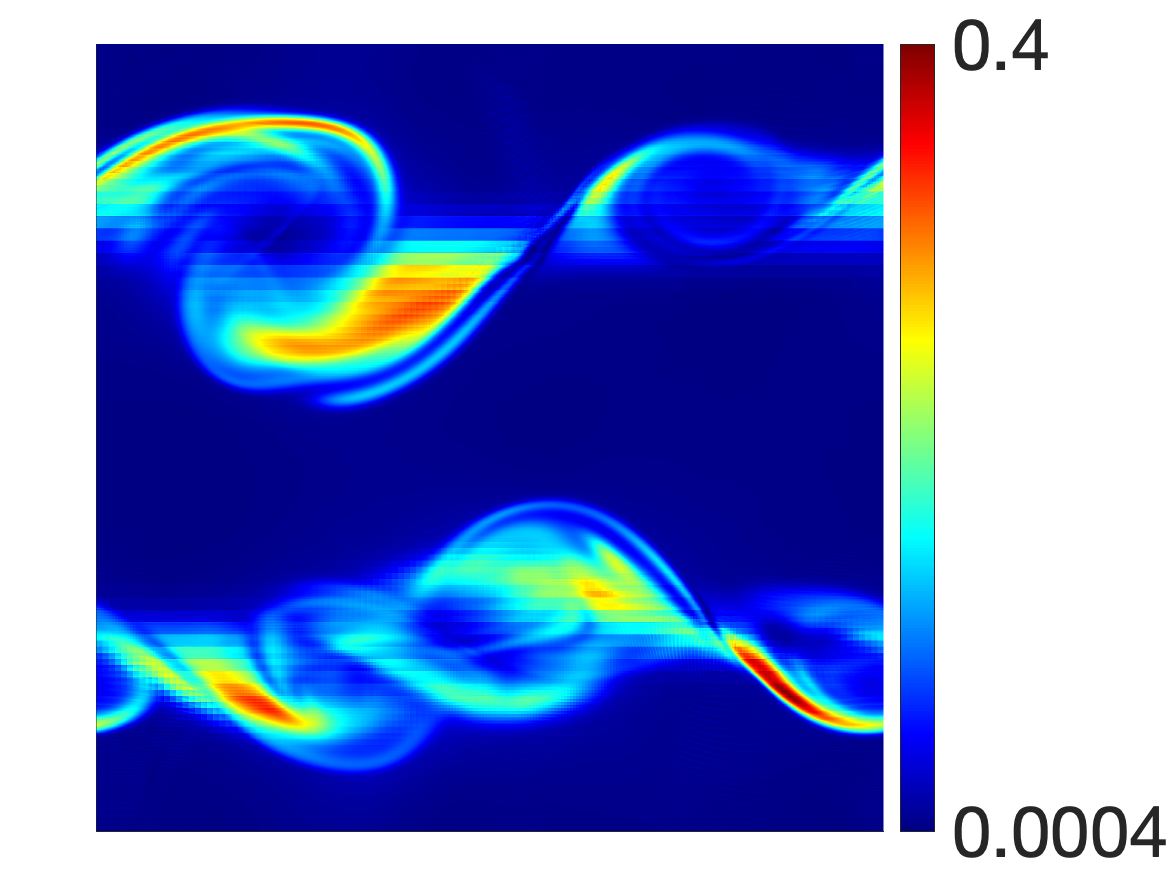}
\end{subfigure}
\begin{subfigure}{0.245\linewidth} \centering
\includegraphics[width=1\textwidth]{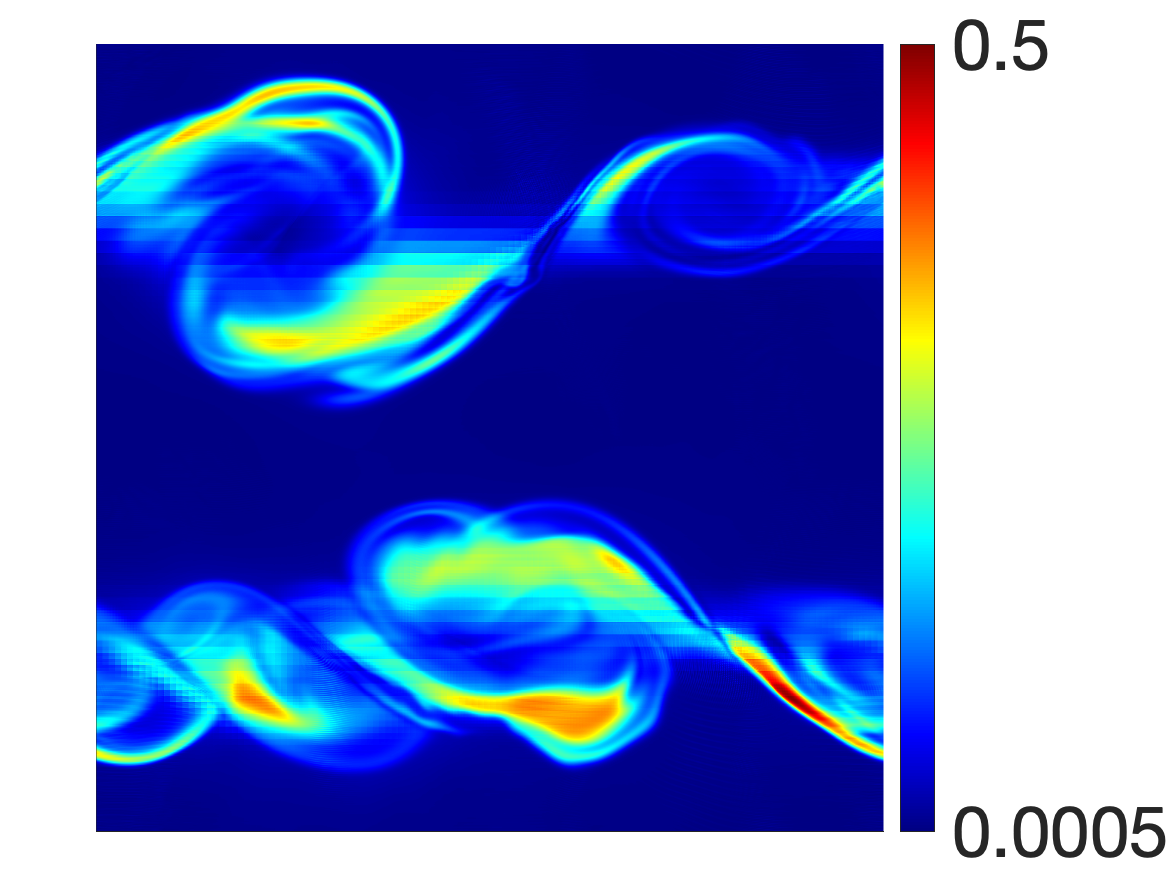}
\end{subfigure}
\\
\begin{subfigure}{0.245\linewidth} \centering
\includegraphics[width=1\textwidth]{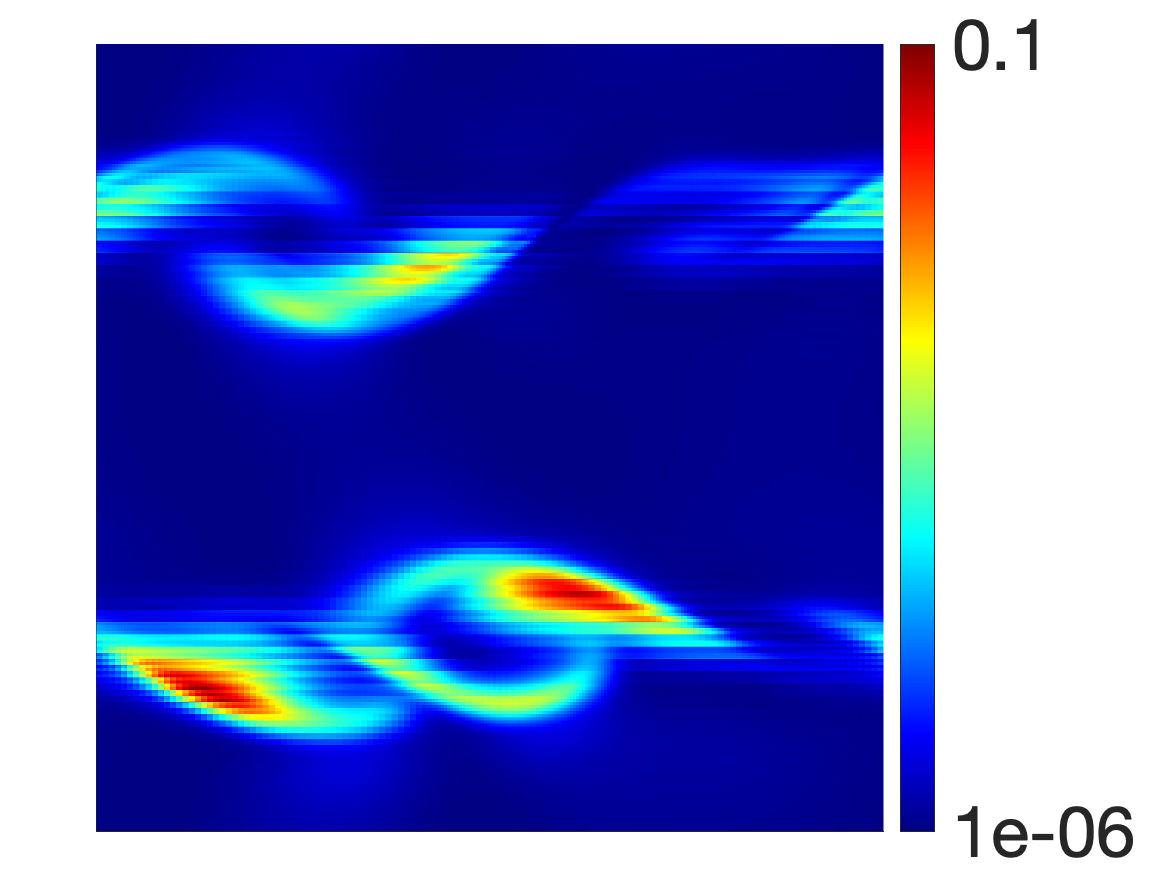}
\end{subfigure}
\begin{subfigure}{0.245\linewidth} \centering
\includegraphics[width=1\textwidth]{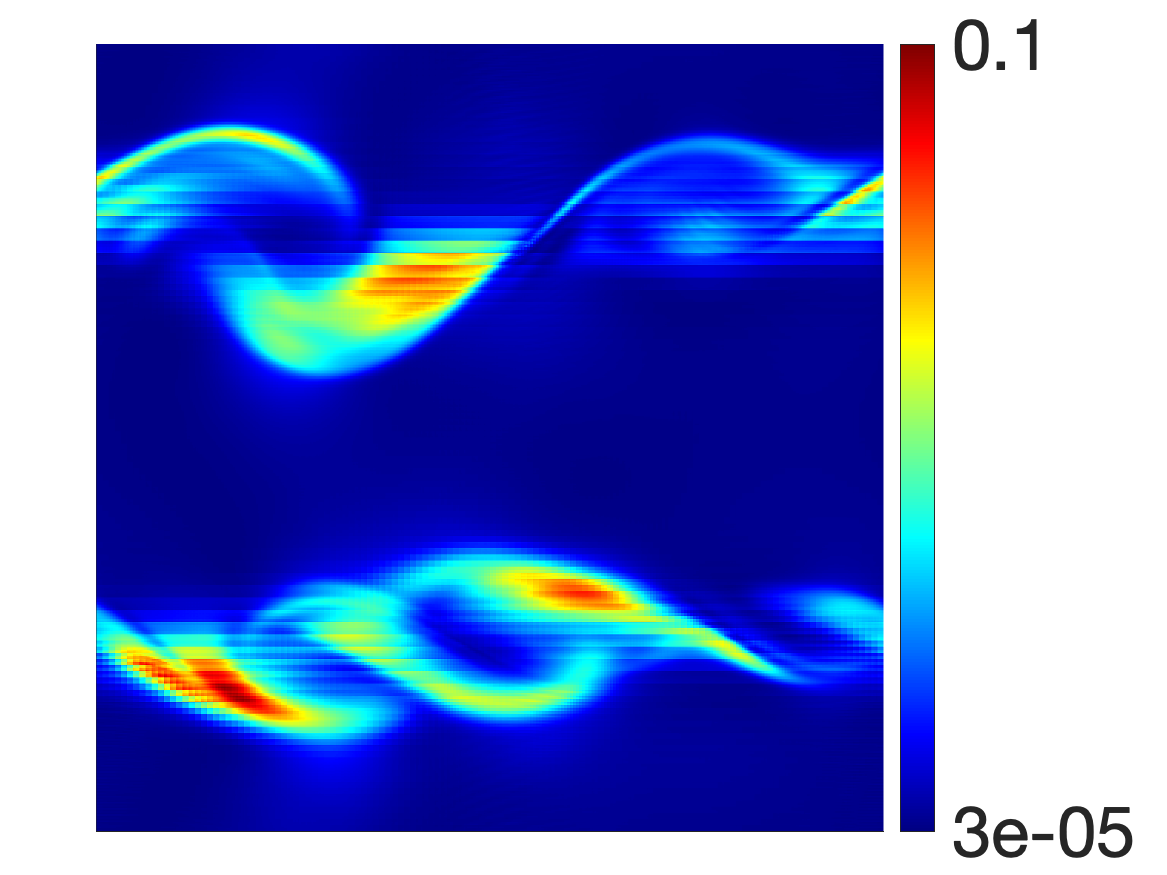}
\end{subfigure}
\begin{subfigure}{0.245\linewidth} \centering
\includegraphics[width=1\textwidth]{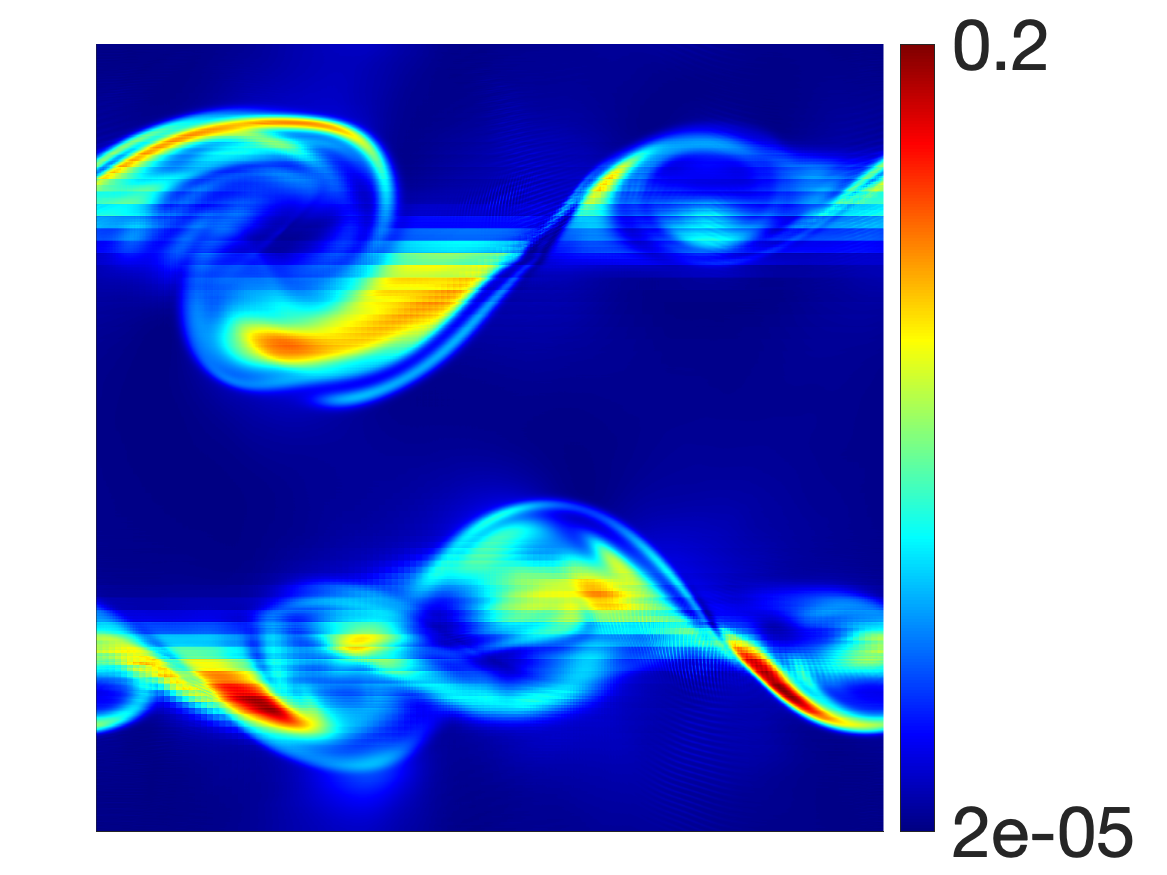}
\end{subfigure}
\begin{subfigure}{0.245\linewidth} \centering
\includegraphics[width=1\textwidth]{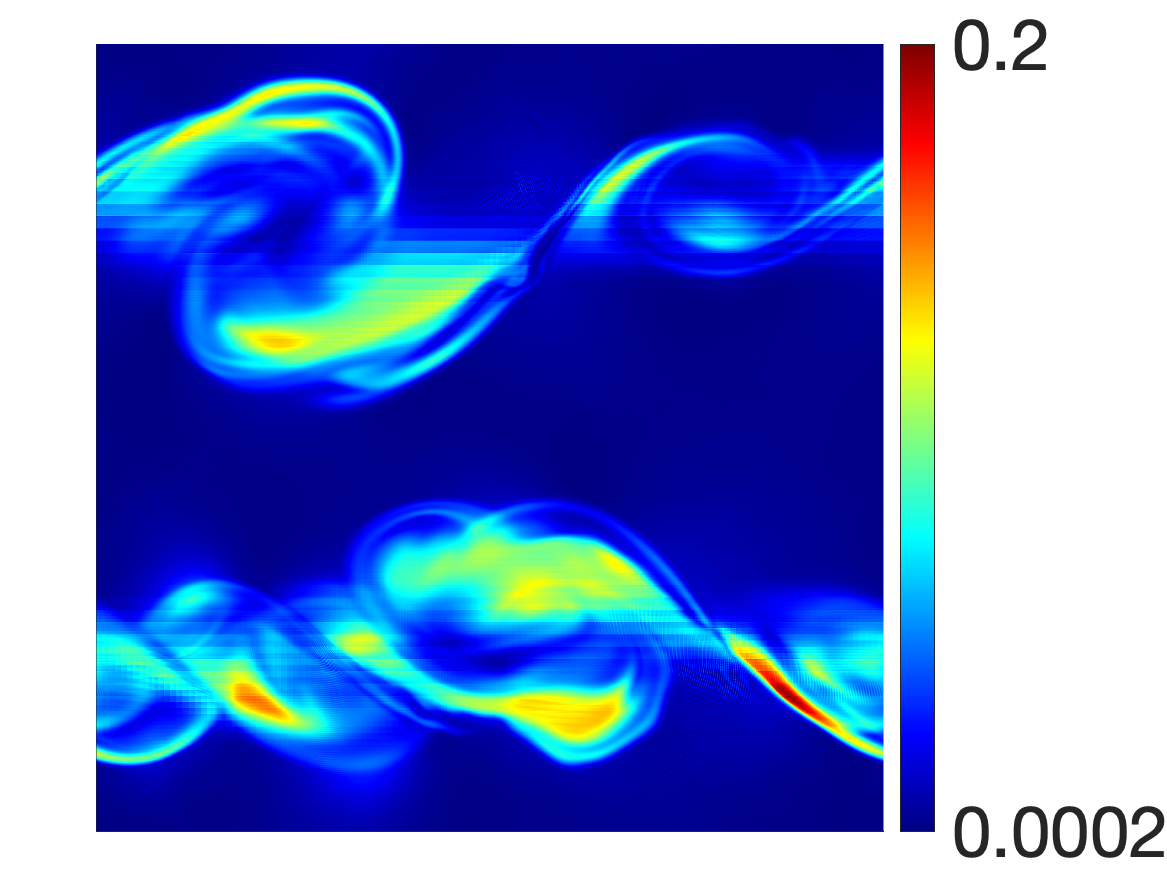}
\end{subfigure}
\\
\begin{subfigure}{0.245\linewidth} \centering
\includegraphics[width=1\textwidth]{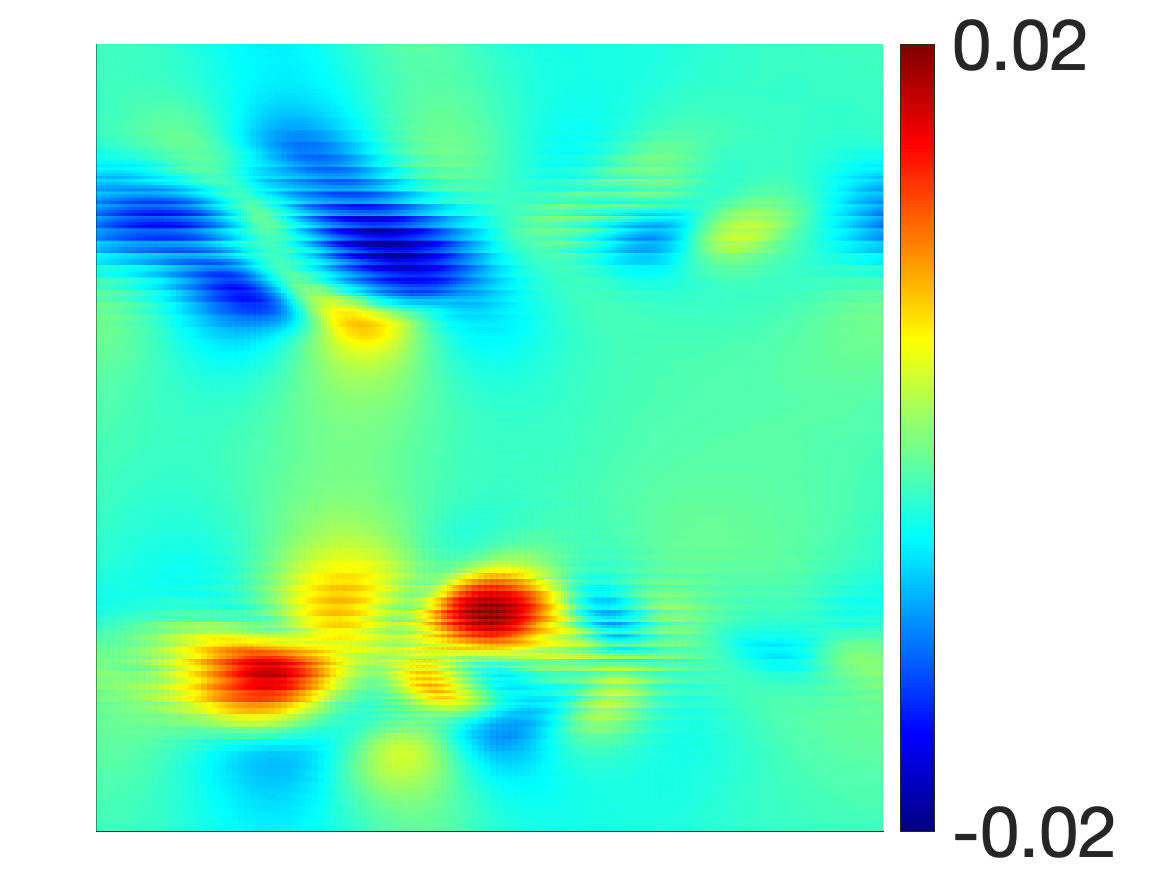}
\end{subfigure}
\begin{subfigure}{0.245\linewidth} \centering
\includegraphics[width=1\textwidth]{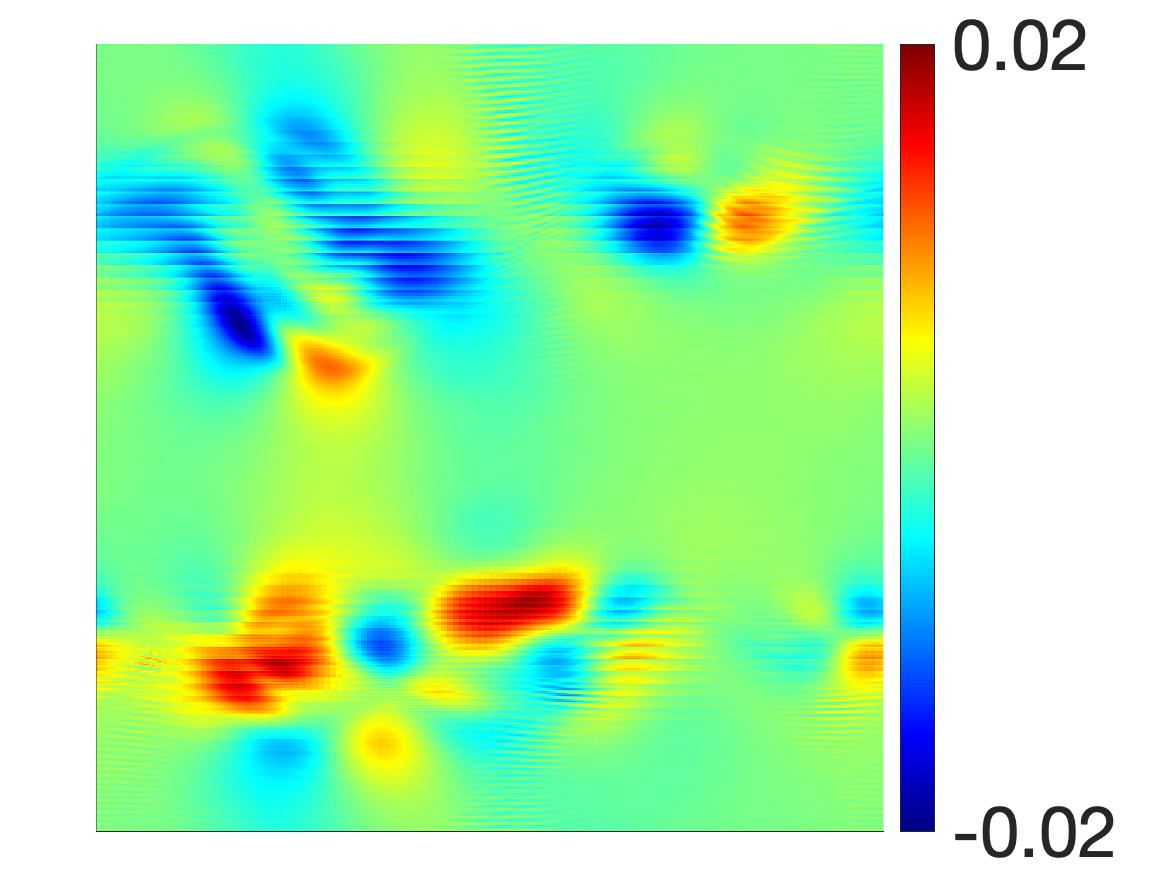}
\end{subfigure}
\begin{subfigure}{0.245\linewidth} \centering
\includegraphics[width=1\textwidth]{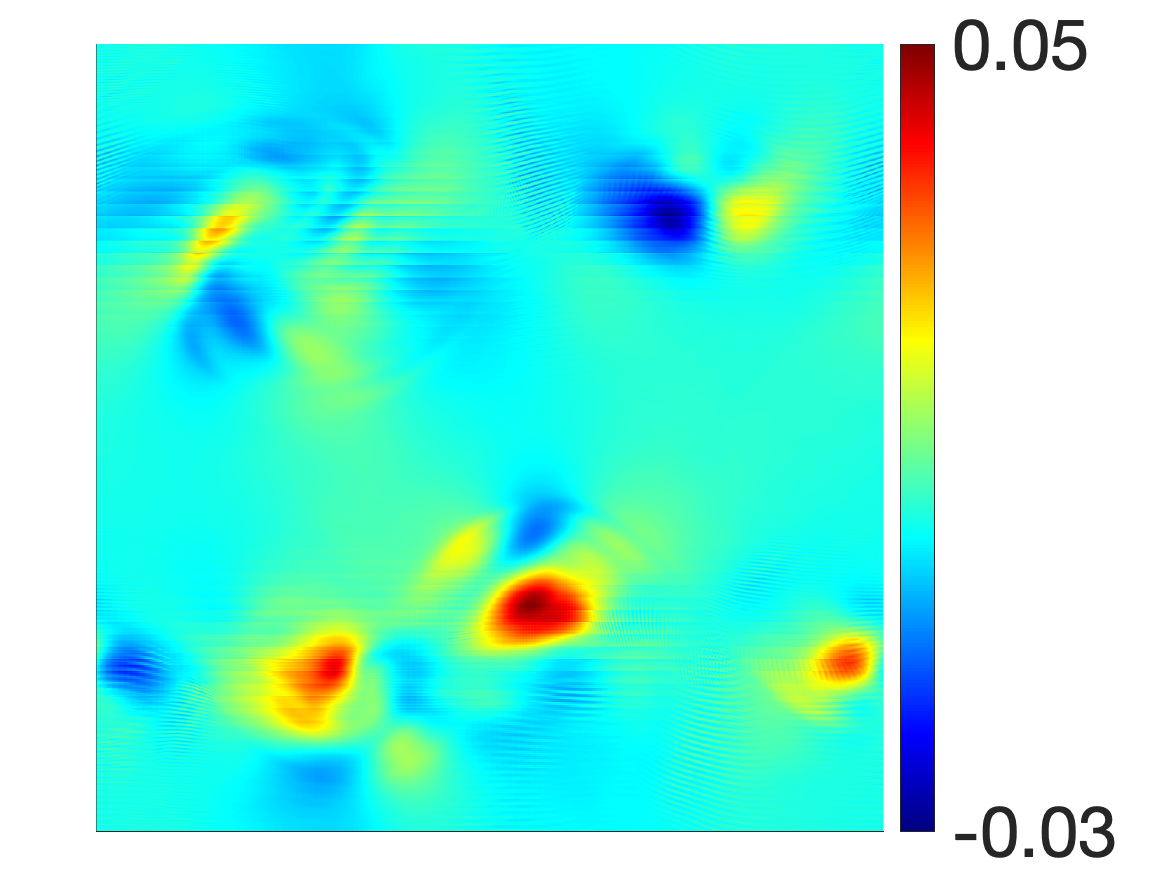}
\end{subfigure}
\begin{subfigure}{0.245\linewidth} \centering
\includegraphics[width=1\textwidth]{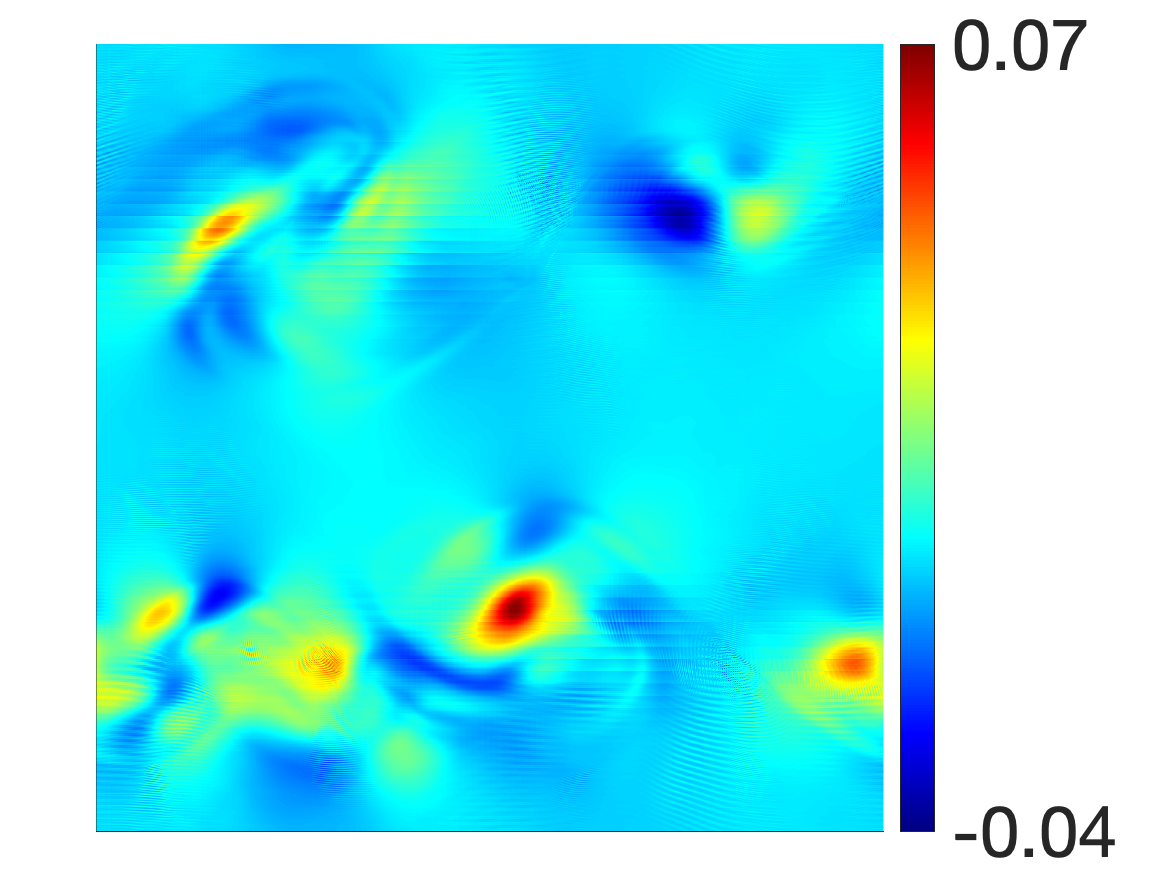}
\end{subfigure}
\\
\begin{subfigure}{0.245\linewidth} \centering
\includegraphics[width=1\textwidth]{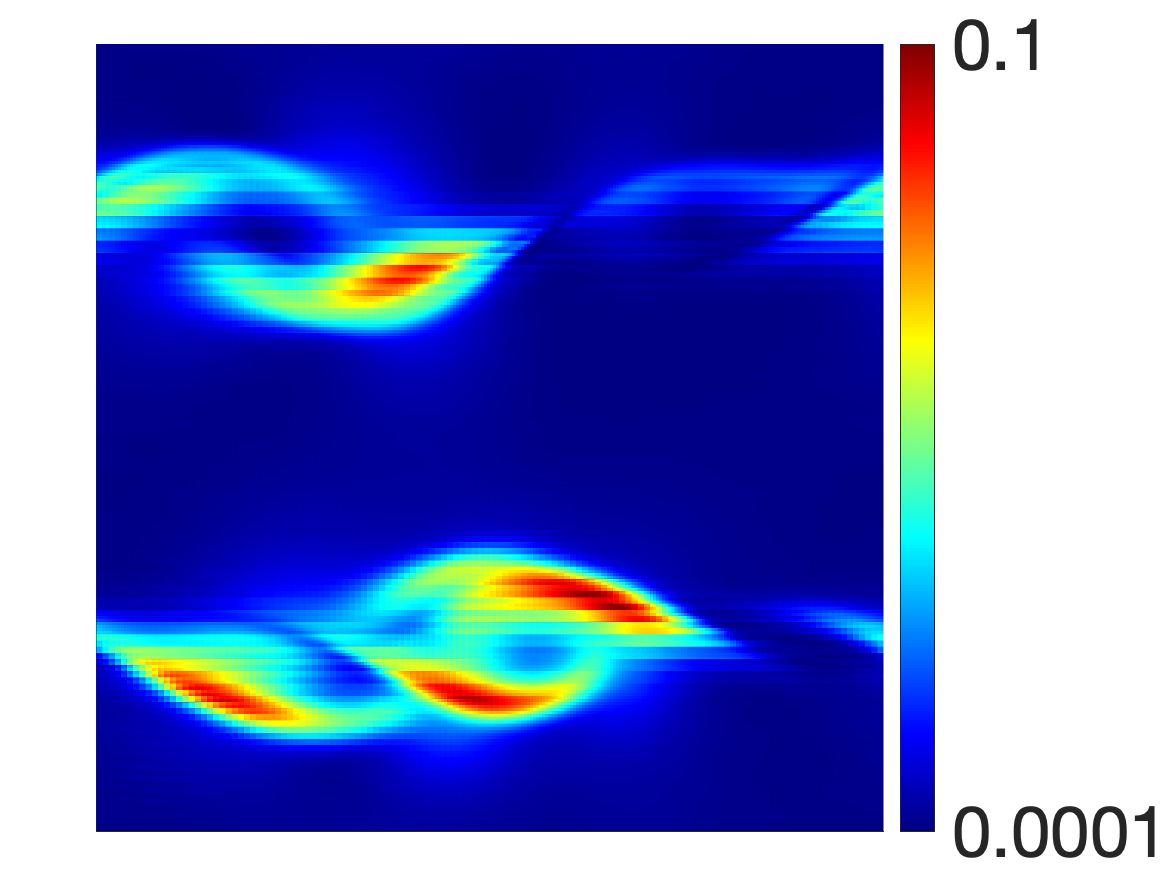}
\end{subfigure}
\begin{subfigure}{0.245\linewidth} \centering
\includegraphics[width=1\textwidth]{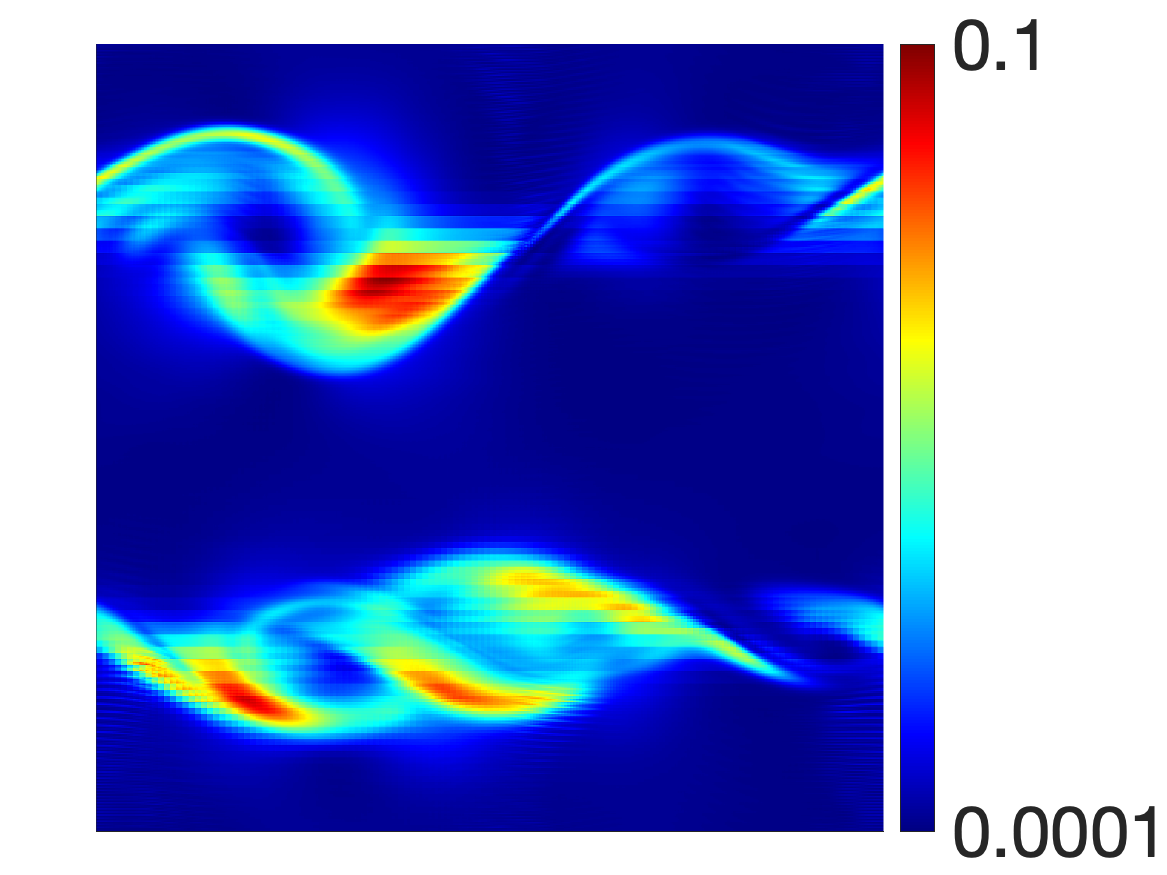}
\end{subfigure}
\begin{subfigure}{0.245\linewidth} \centering
\includegraphics[width=1\textwidth]{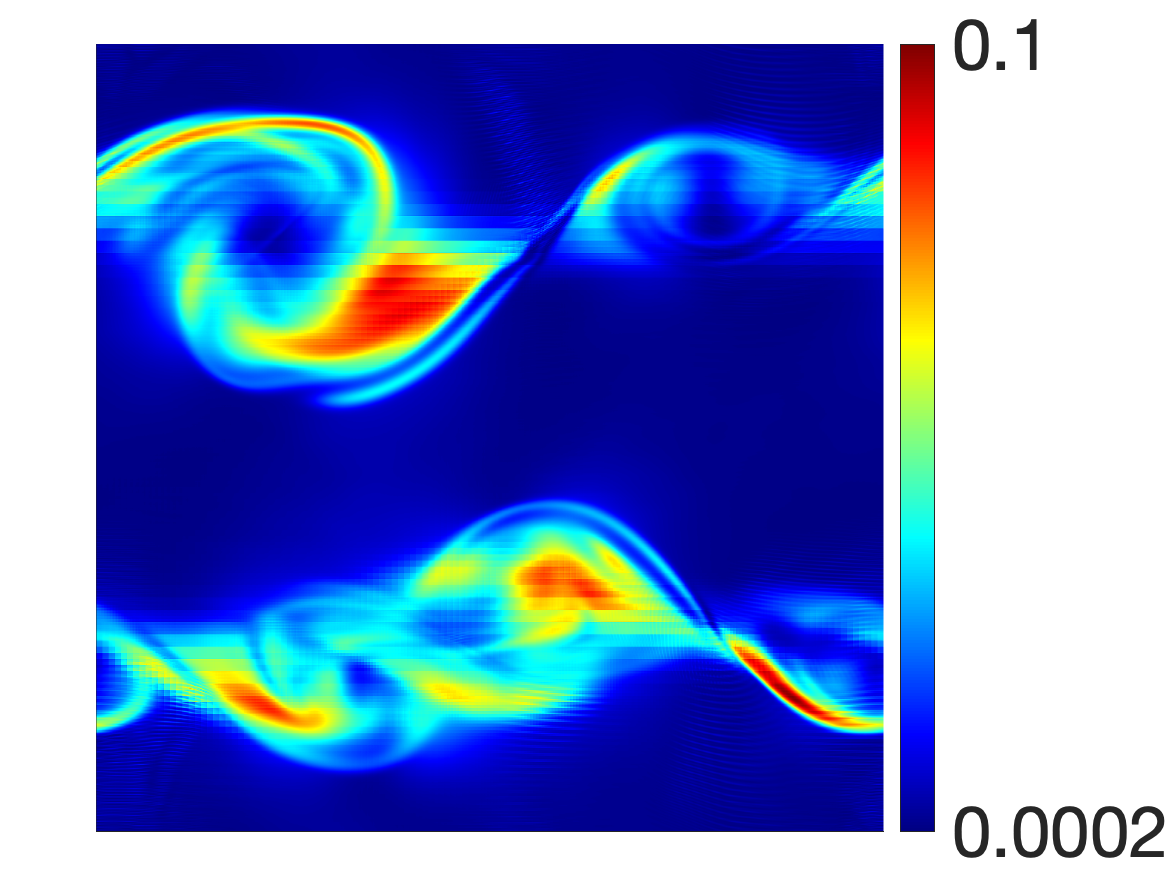}
\end{subfigure}
\begin{subfigure}{0.245\linewidth} \centering
\includegraphics[width=1\textwidth]{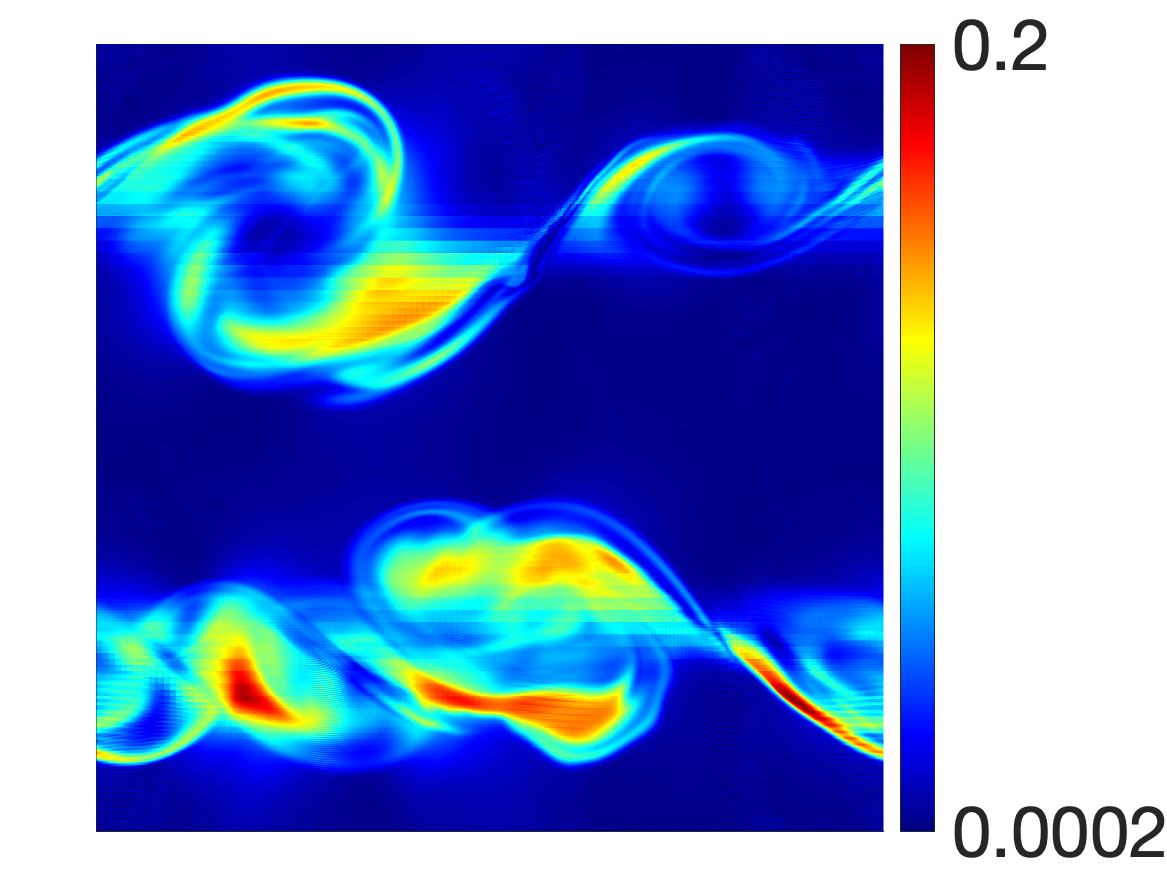}
\end{subfigure}
\\
\begin{subfigure}{0.245\linewidth} \centering
\includegraphics[width=1\textwidth]{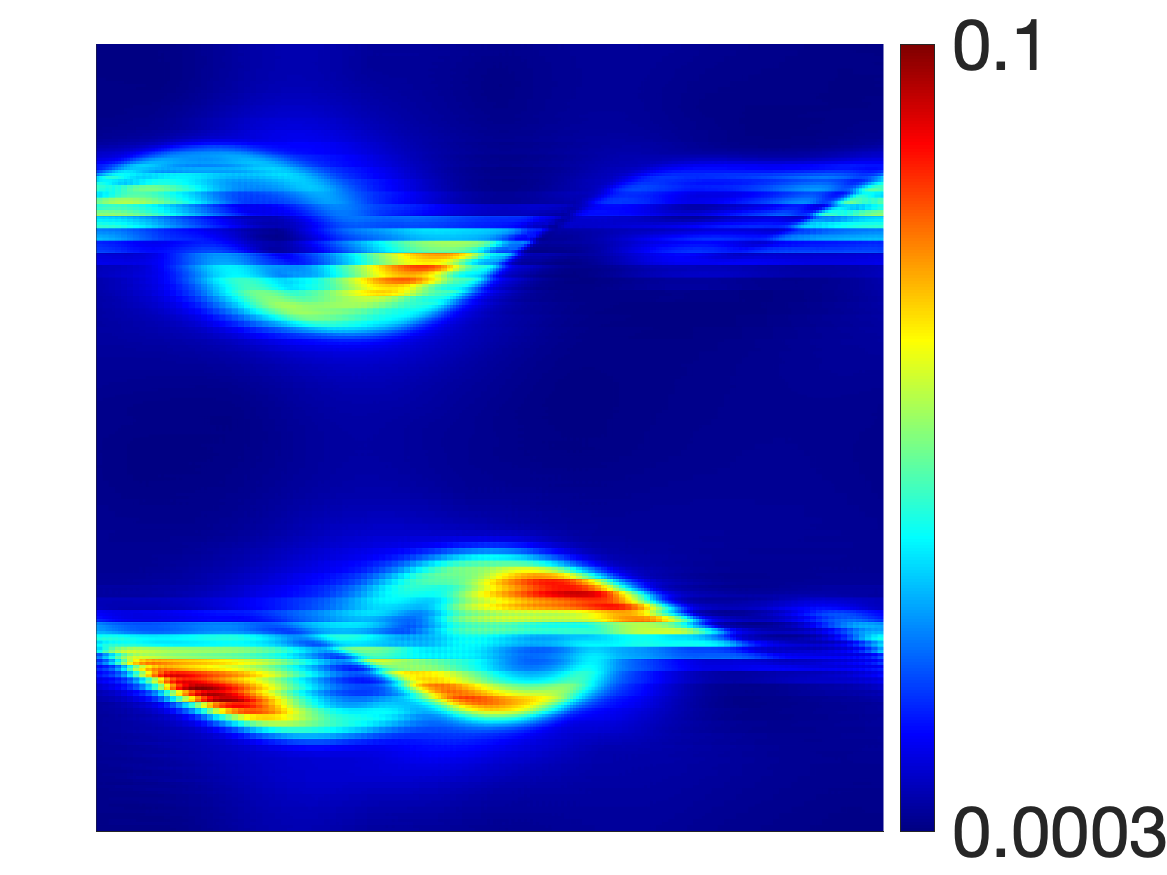}
\end{subfigure}
\begin{subfigure}{0.245\linewidth} \centering
\includegraphics[width=1\textwidth]{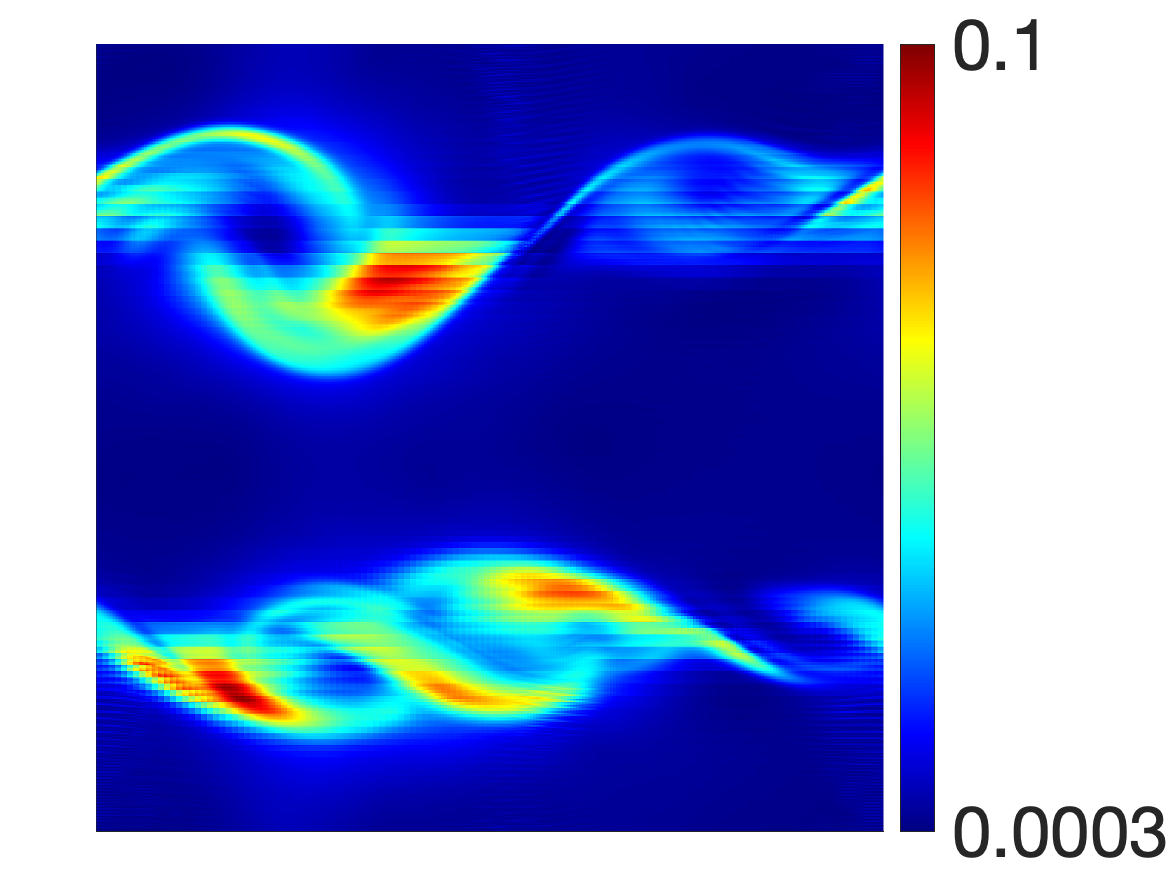}
\end{subfigure}
\begin{subfigure}{0.245\linewidth} \centering
\includegraphics[width=1\textwidth]{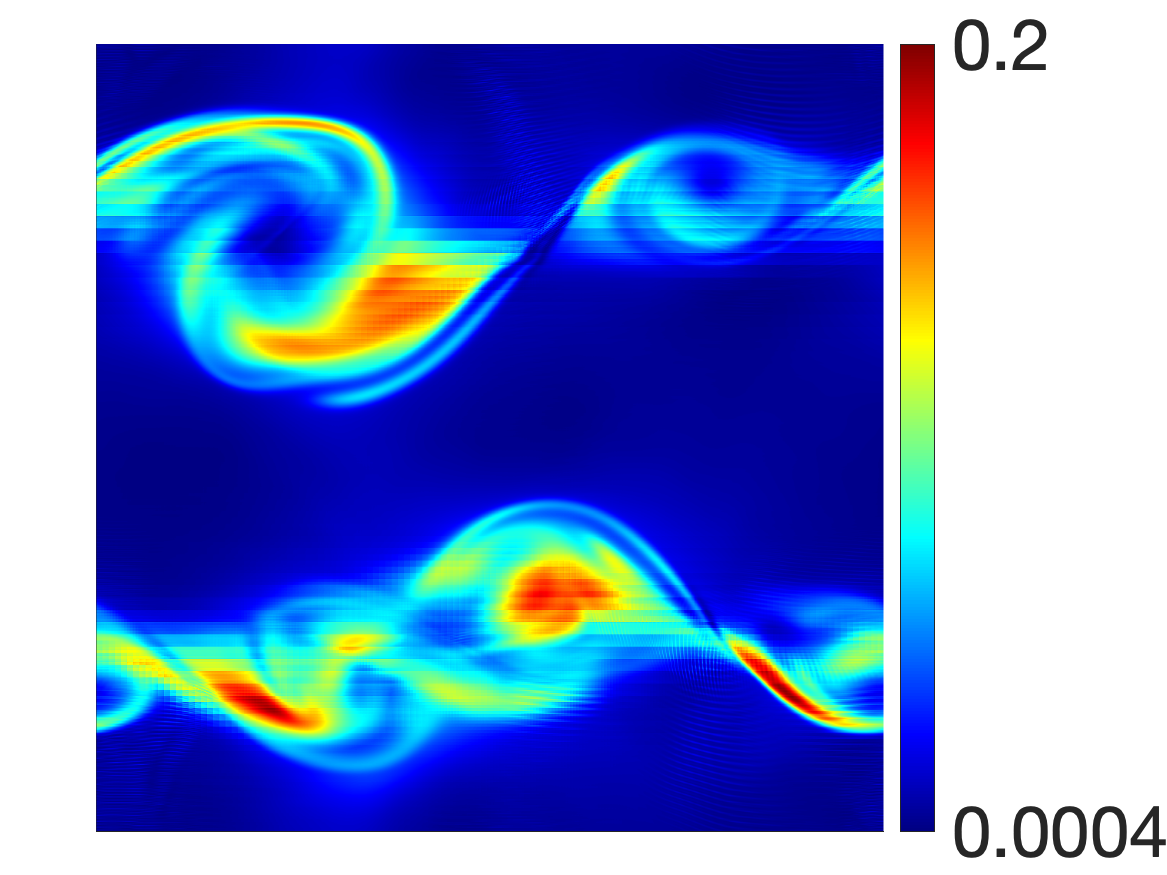}
\end{subfigure}
\begin{subfigure}{0.245\linewidth} \centering
\includegraphics[width=1\textwidth]{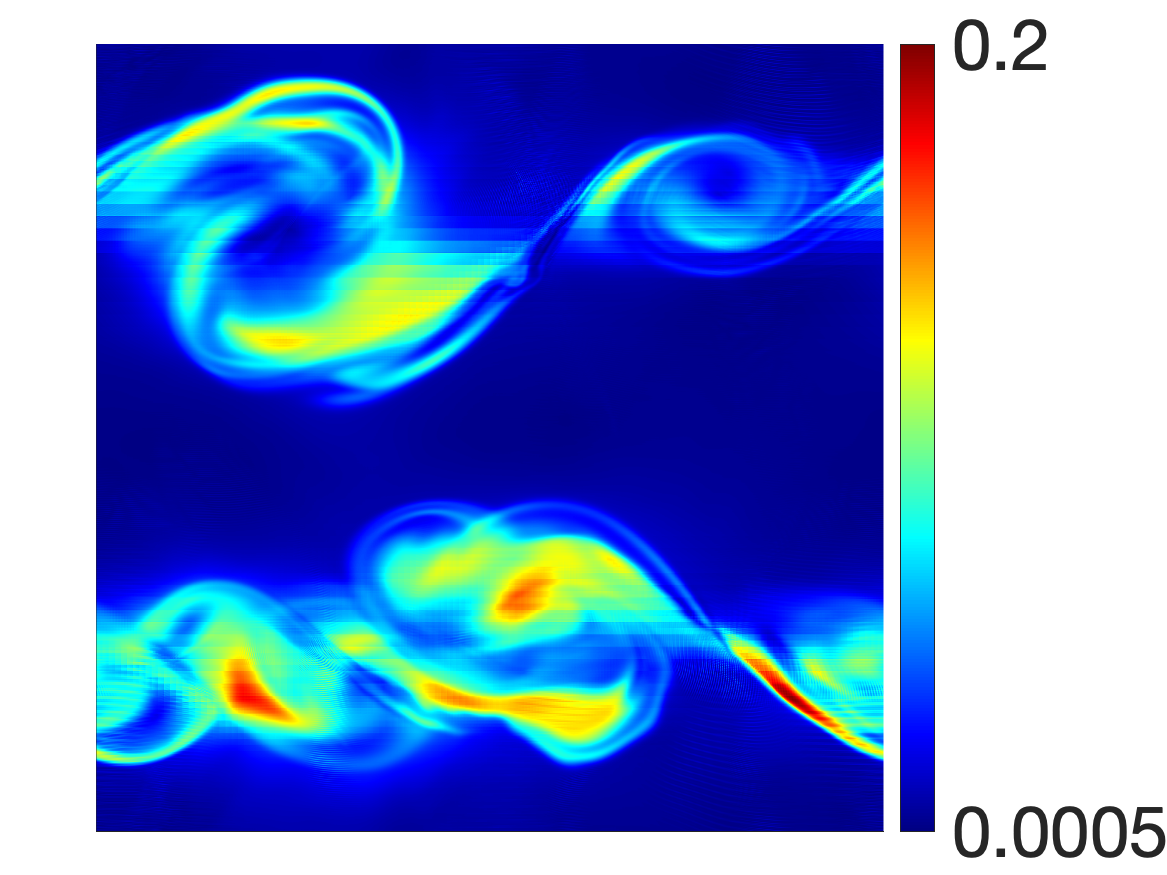}
\end{subfigure}
\\
\begin{subfigure}{0.245\linewidth} \centering
\includegraphics[width=1\textwidth]{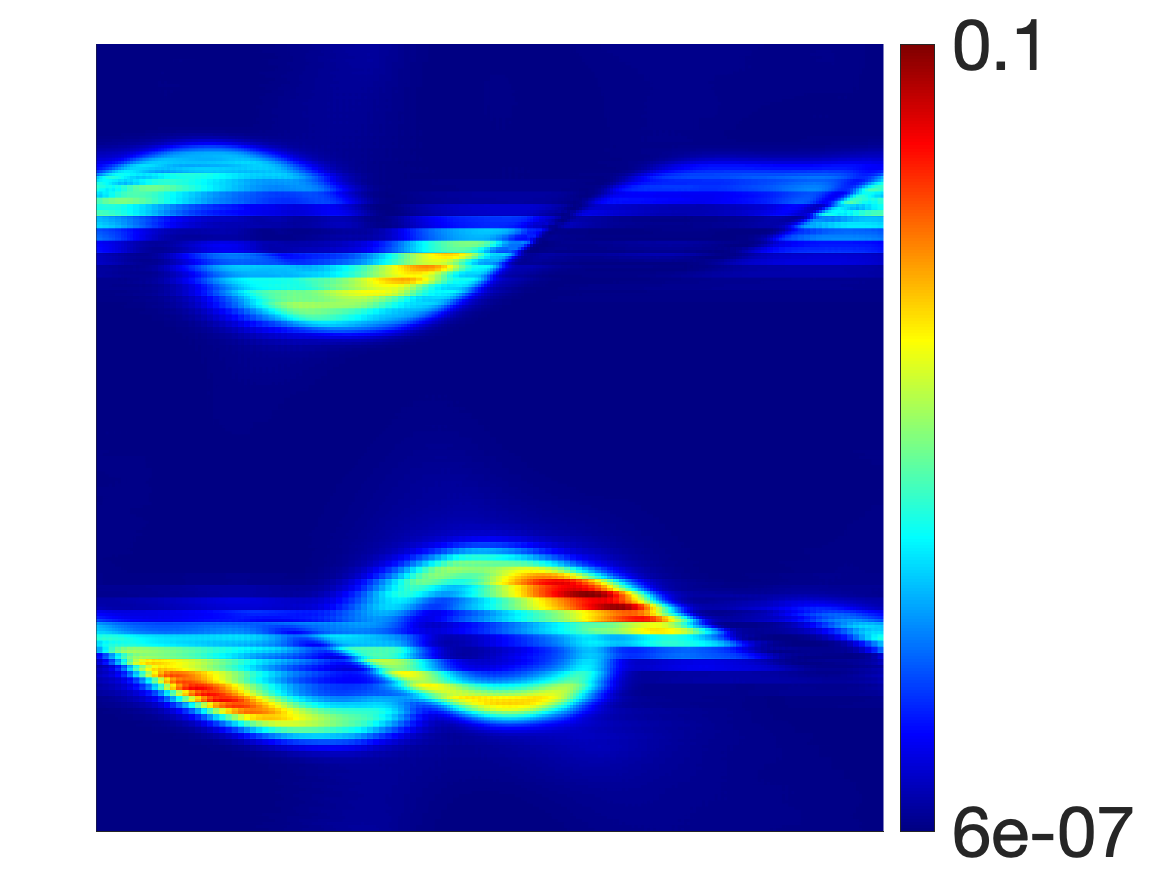}
\caption{$N=3$}
\end{subfigure}
\begin{subfigure}{0.245\linewidth} \centering
\includegraphics[width=1\textwidth]{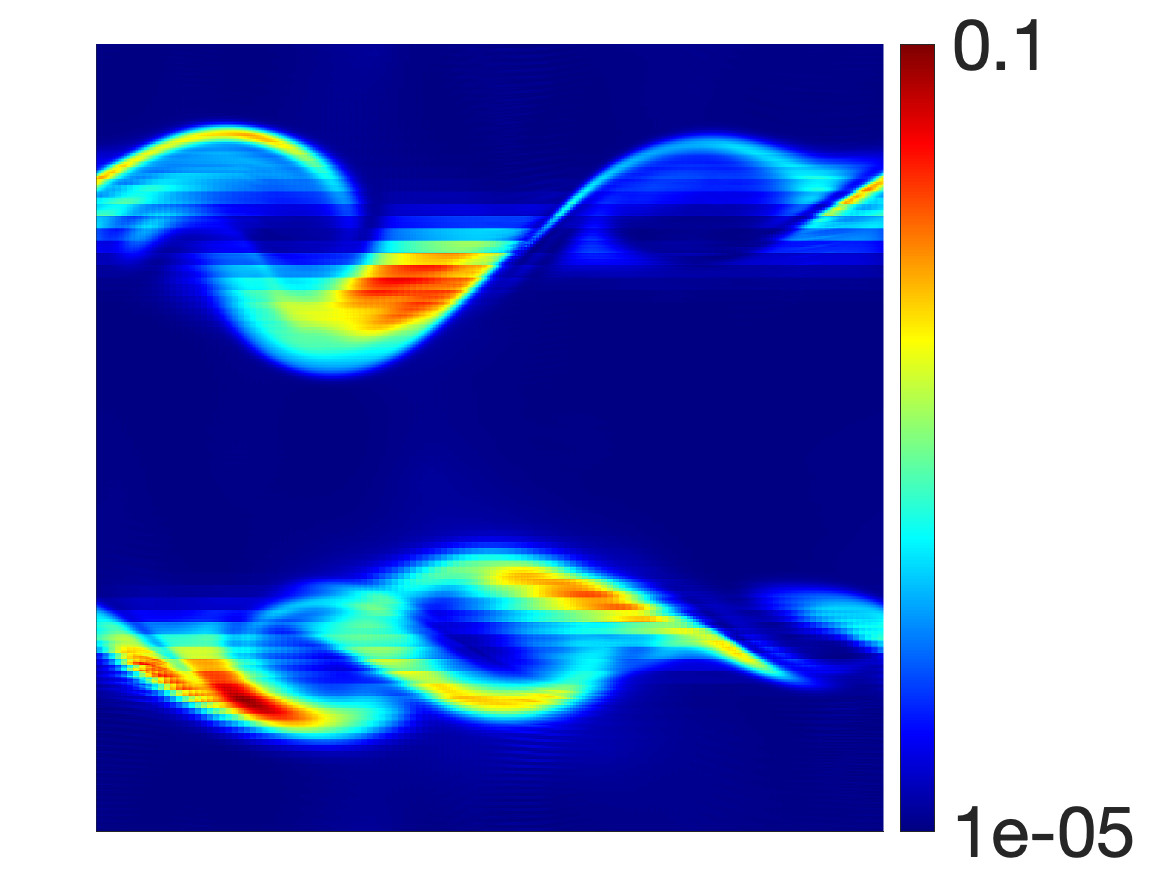}
\caption{$N=4$}
\end{subfigure}
\begin{subfigure}{0.245\linewidth} \centering
\includegraphics[width=1\textwidth]{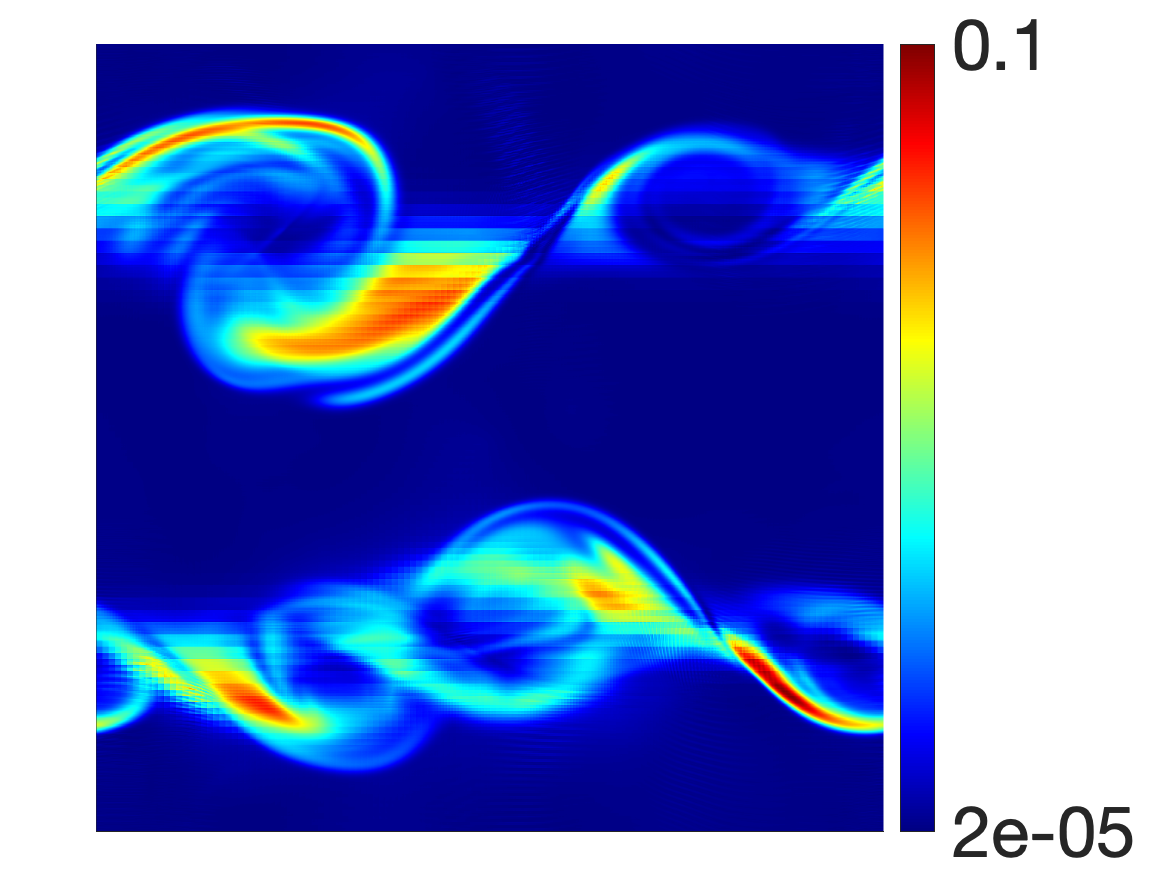}
\caption{ $N=5$}
\end{subfigure}
\begin{subfigure}{0.245\linewidth} \centering
\includegraphics[width=1\textwidth]{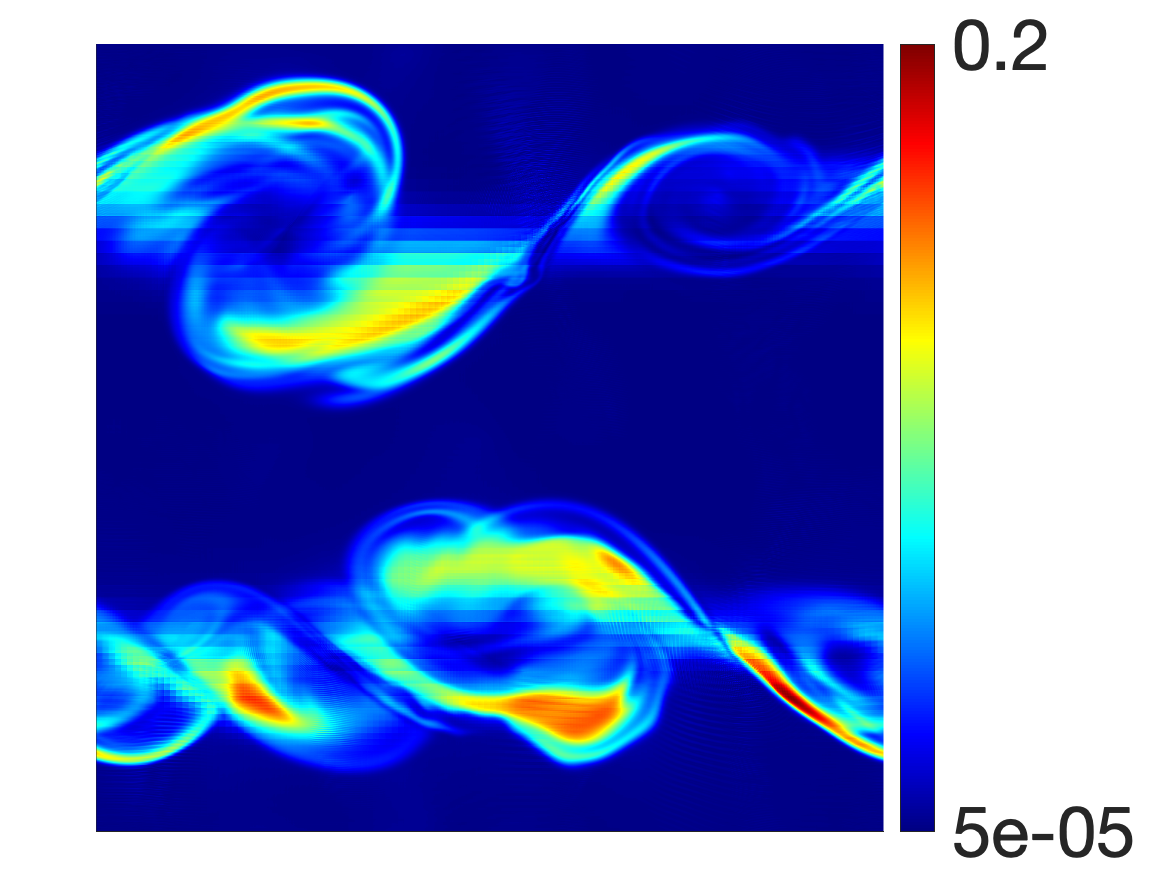}
\caption{$N=6$}
\end{subfigure}
\caption{Defects computed from the VFV solution.  From top to  bottom: $\mathfrak{E}_{N}$,  $\mathfrak{R}^{11}_N$, $\mathfrak{R}^{12}_N$, $\mathfrak{R}^{22}_N$,  $\lambda_1(\mathfrak{R}_N)$ $\lambda_2(\mathfrak{R}_N),$ $N=3, \dots, 6.$}
\label{fig-Rey-1}
\end{figure}
\begin{figure}[h!]
\centering
\begin{subfigure}{0.245\linewidth} \centering
\includegraphics[width=1\textwidth]{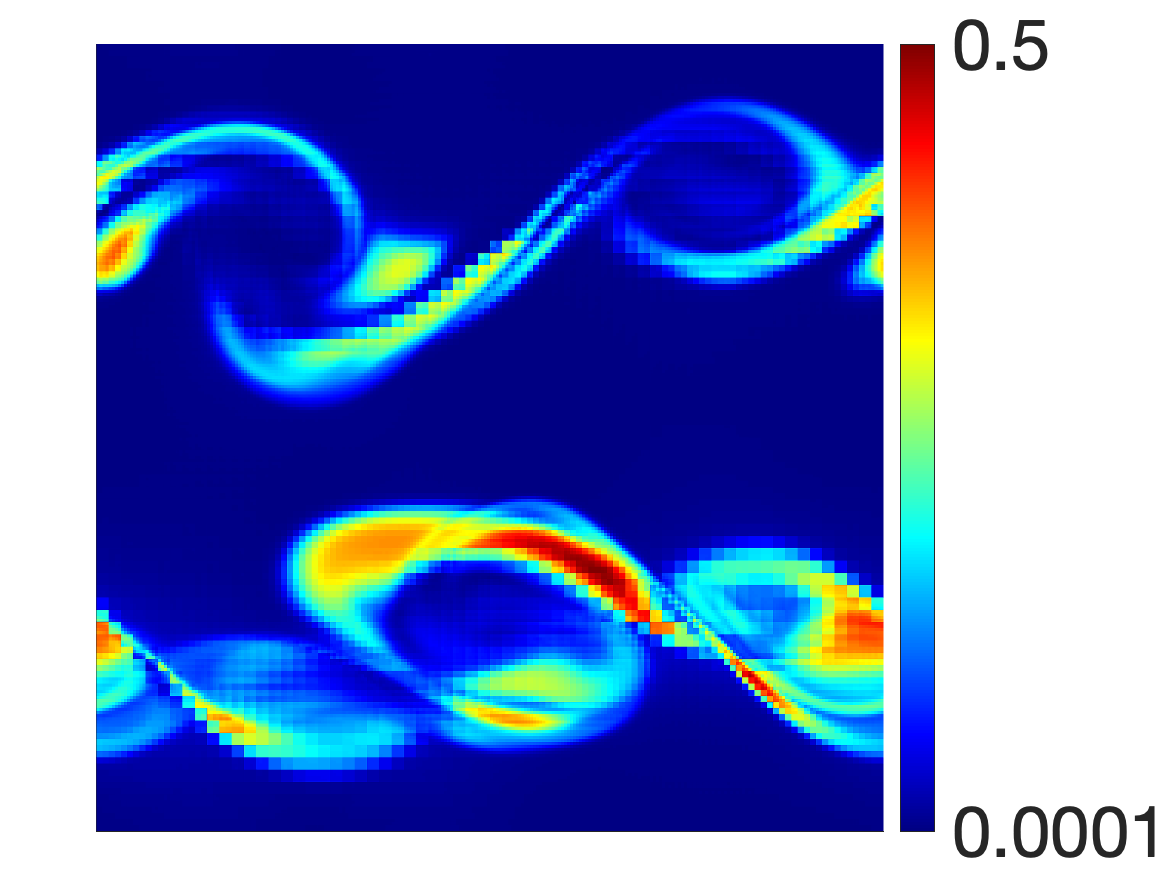}
\end{subfigure}
\begin{subfigure}{0.245\linewidth} \centering
\includegraphics[width=1\textwidth]{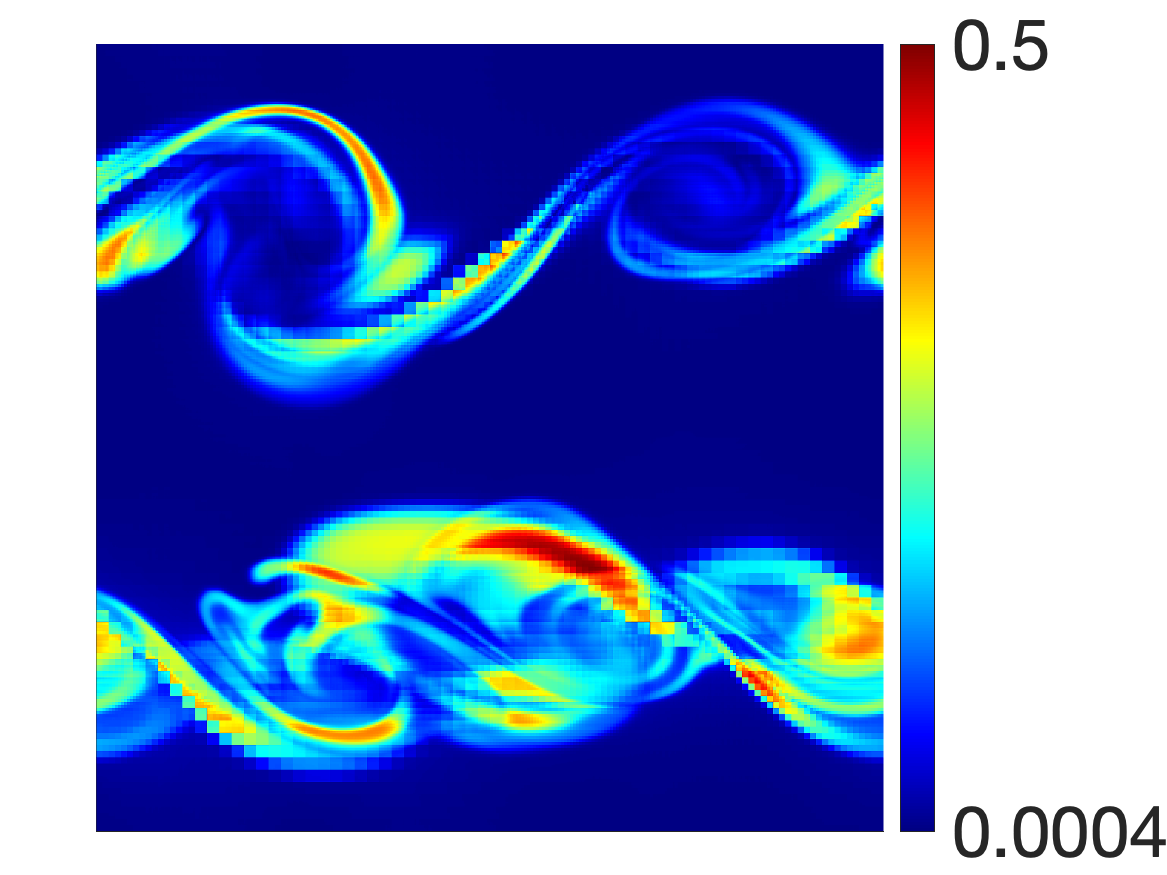}
\end{subfigure}
\begin{subfigure}{0.245\linewidth} \centering
\includegraphics[width=1\textwidth]{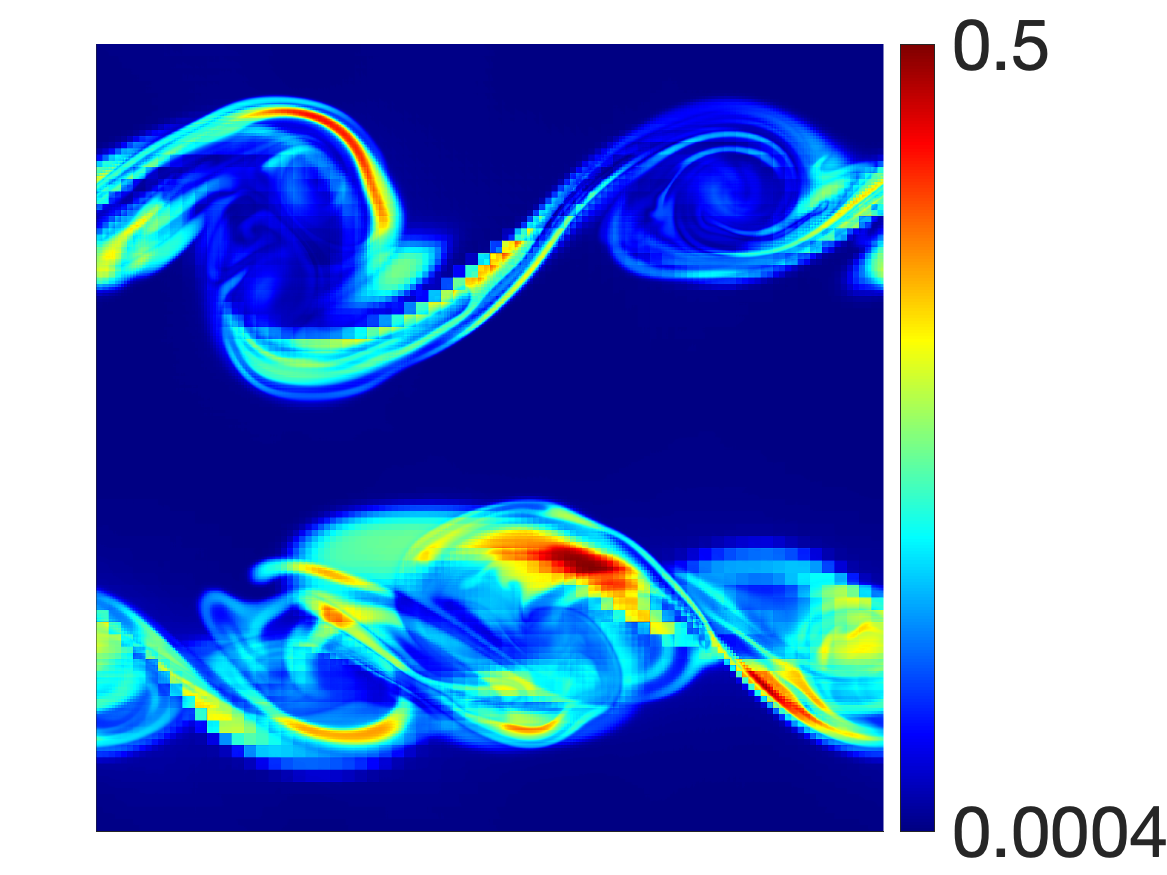}
\end{subfigure}
\begin{subfigure}{0.245\linewidth} \centering
\includegraphics[width=1\textwidth]{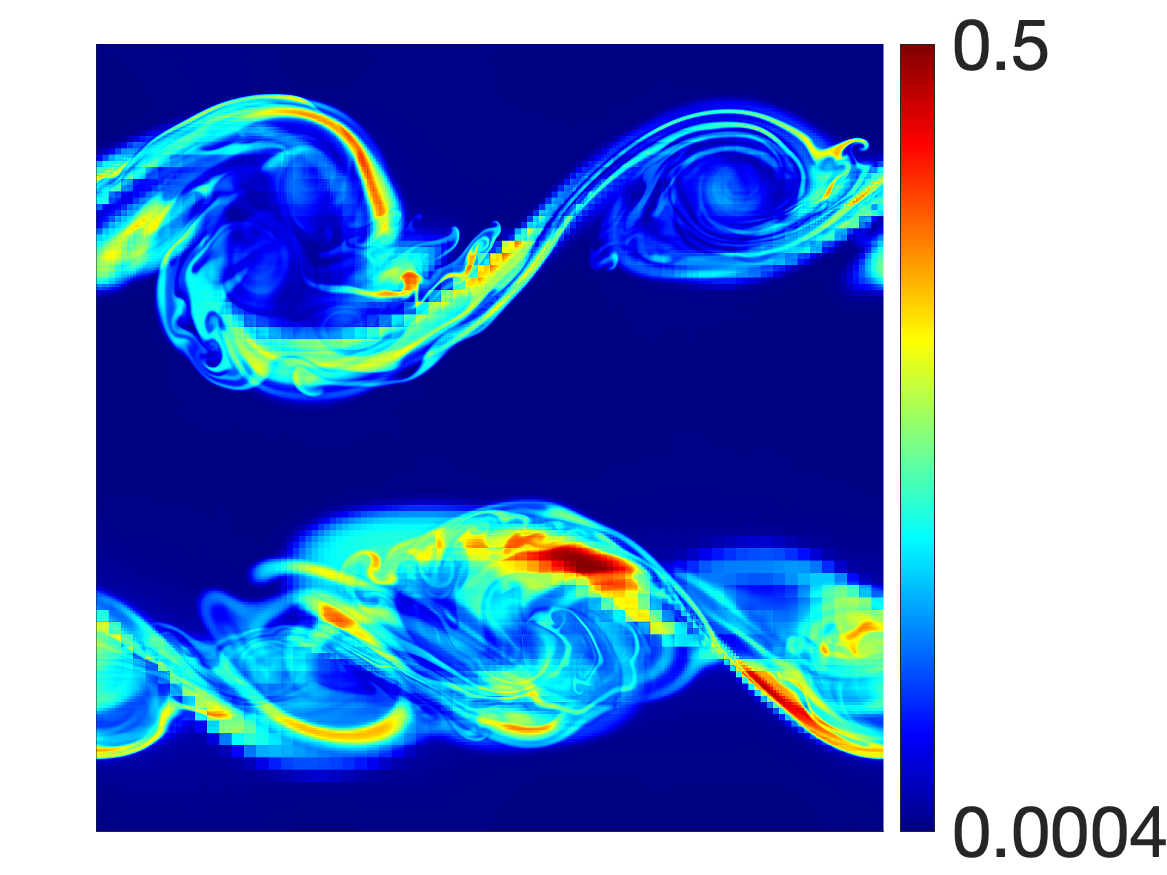}
\end{subfigure}
\\
\begin{subfigure}{0.245\linewidth} \centering
\includegraphics[width=1\textwidth]{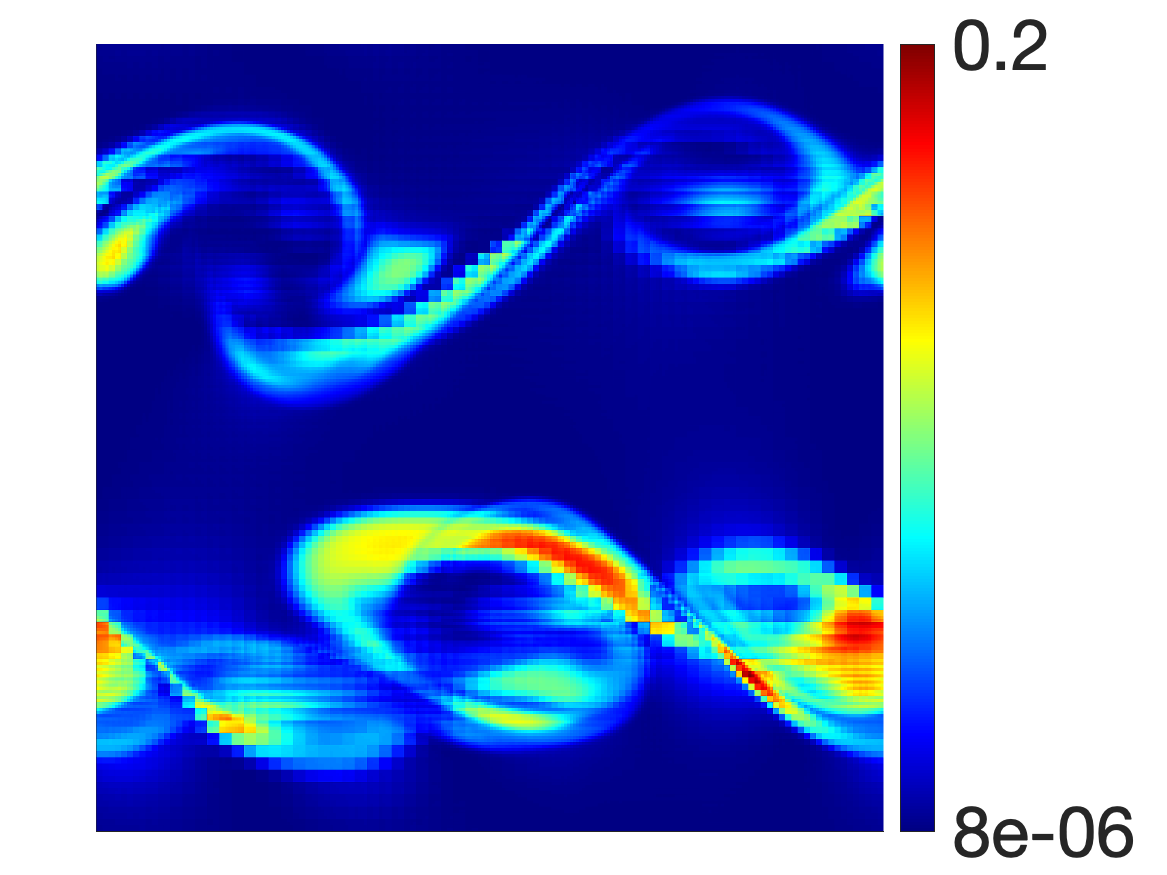}
\end{subfigure}
\begin{subfigure}{0.245\linewidth} \centering
\includegraphics[width=1\textwidth]{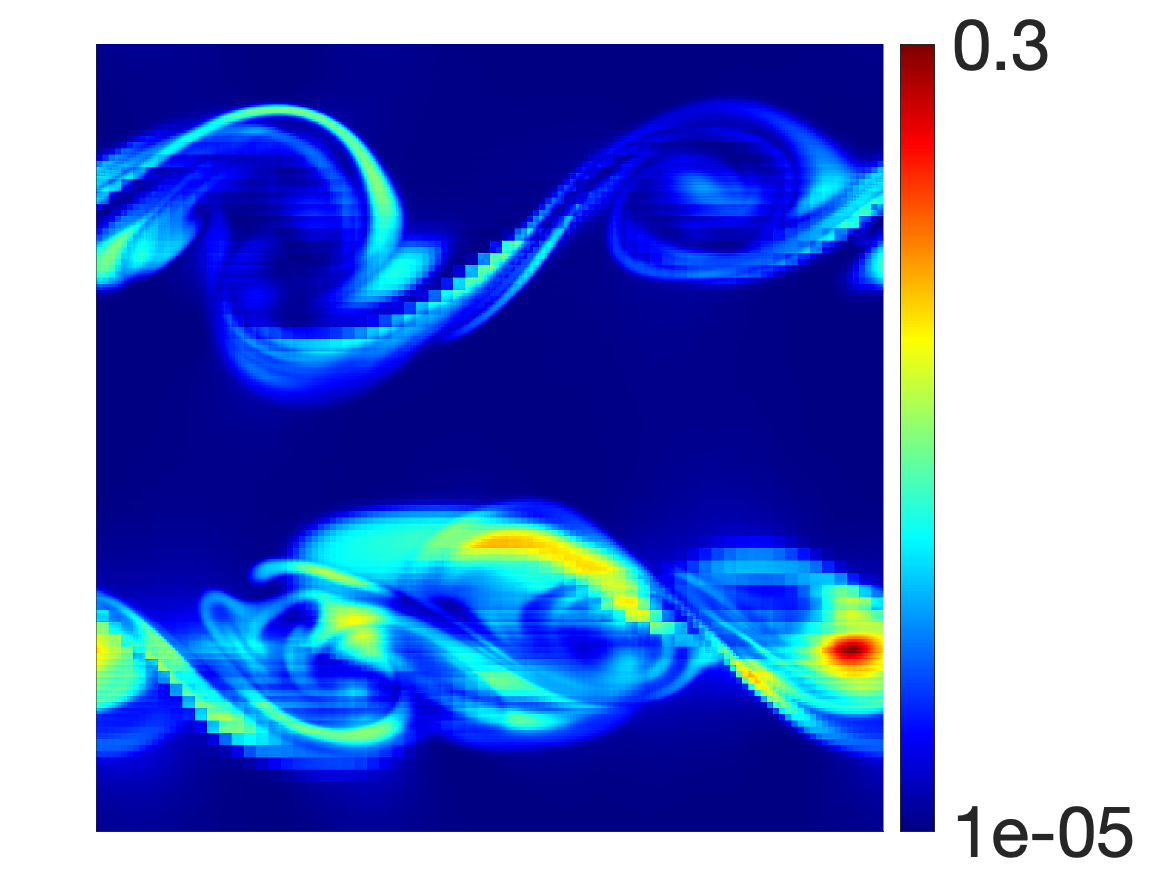}
\end{subfigure}
\begin{subfigure}{0.245\linewidth} \centering
\includegraphics[width=1\textwidth]{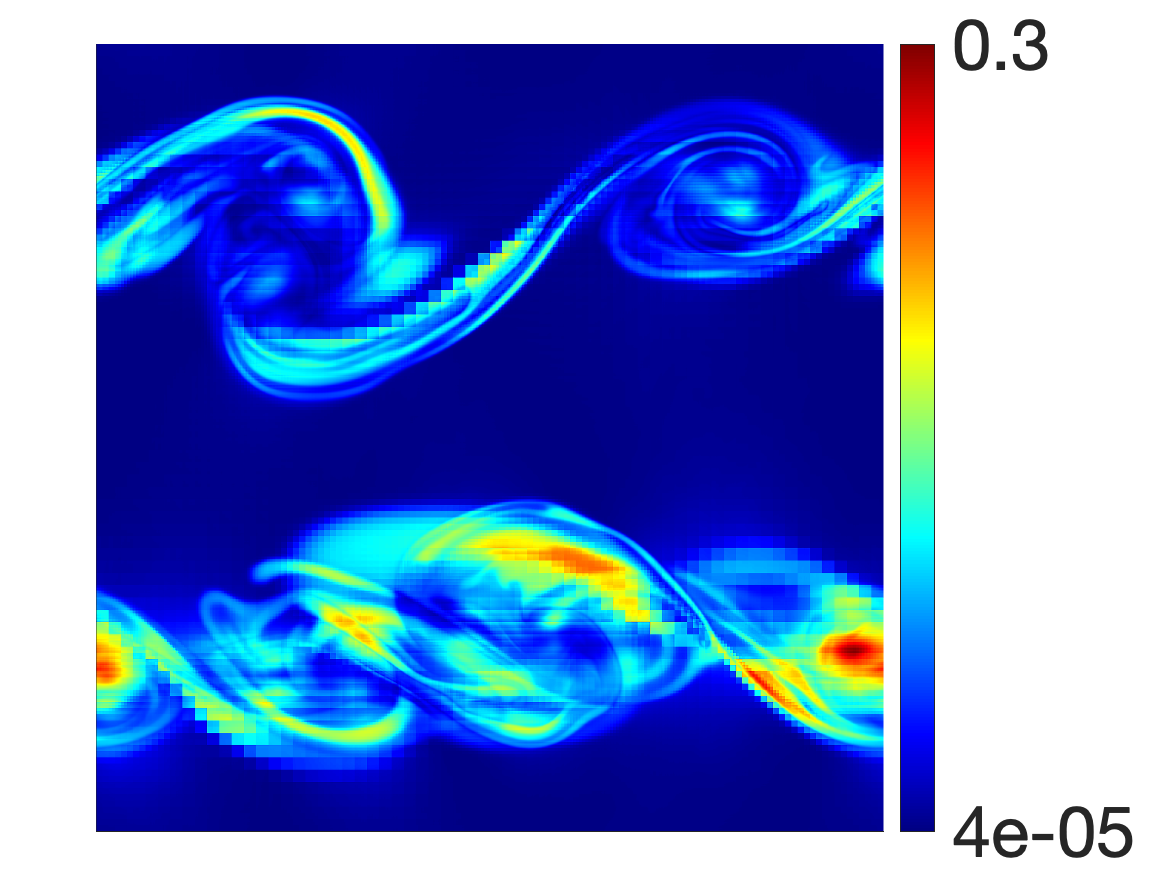}
\end{subfigure}
\begin{subfigure}{0.245\linewidth} \centering
\includegraphics[width=1\textwidth]{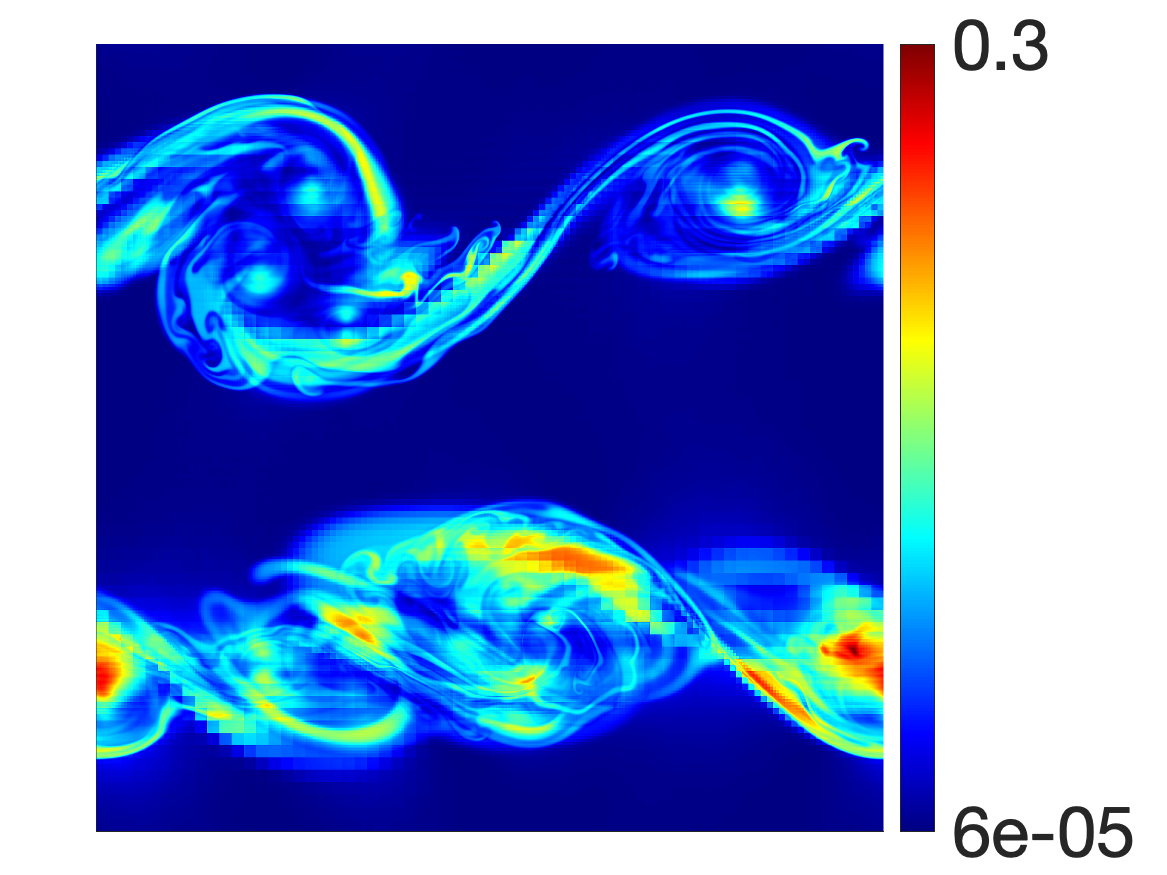}
\end{subfigure}
\\
\begin{subfigure}{0.245\linewidth} \centering
\includegraphics[width=1\textwidth]{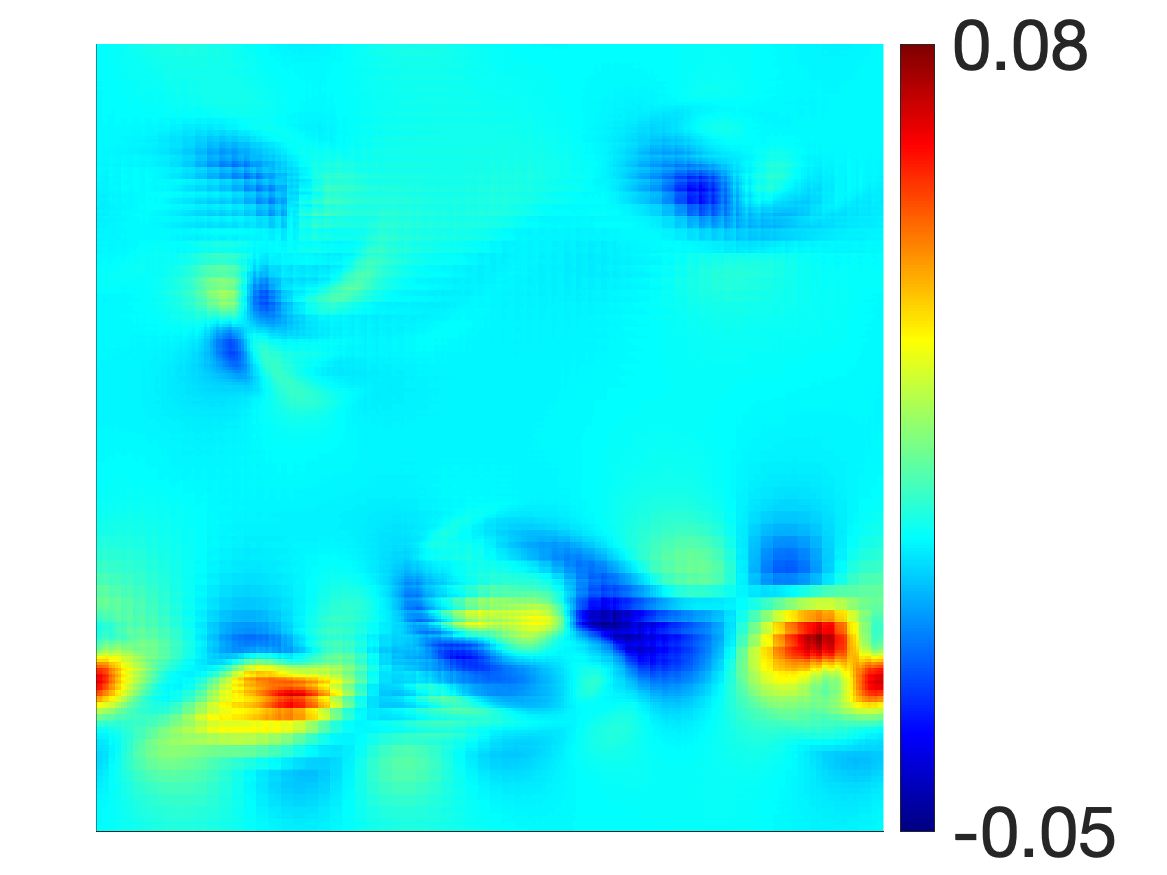}
\end{subfigure}
\begin{subfigure}{0.245\linewidth} \centering
\includegraphics[width=1\textwidth]{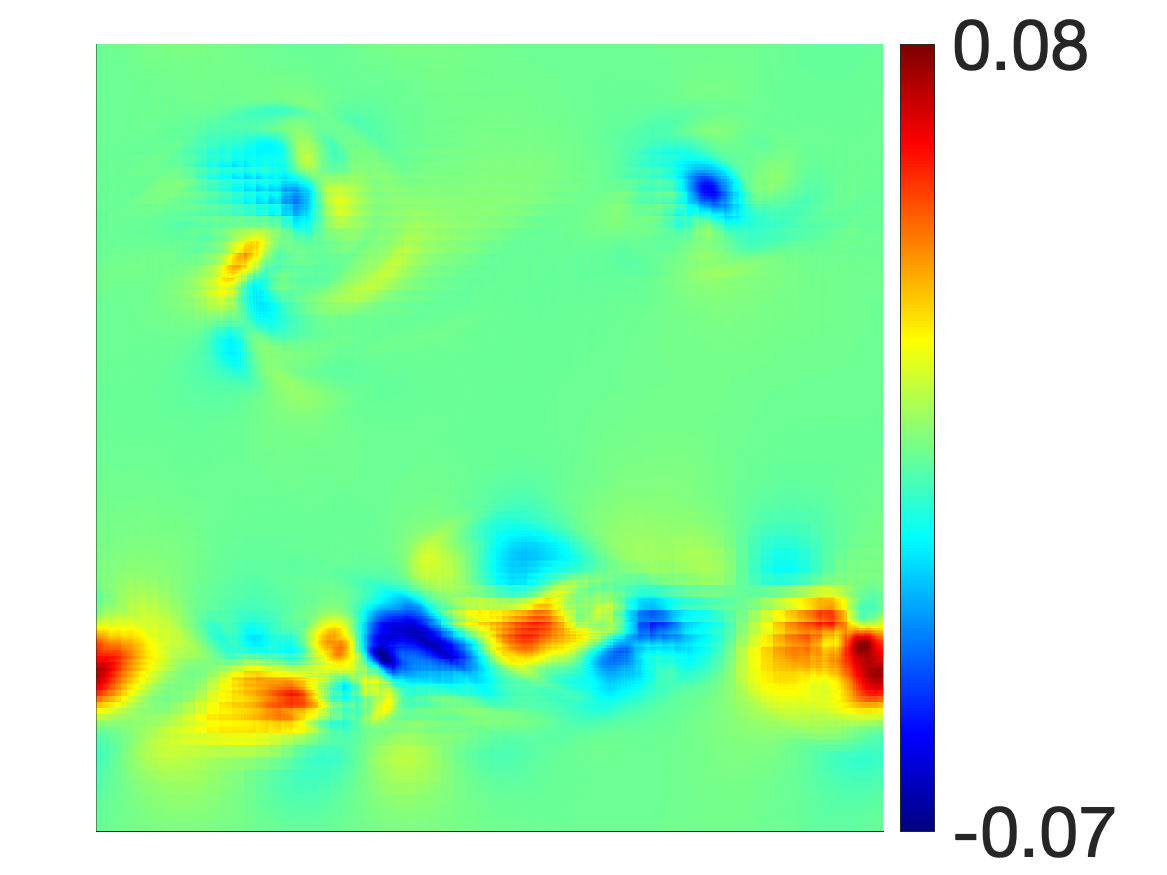}
\end{subfigure}
\begin{subfigure}{0.245\linewidth} \centering
\includegraphics[width=1\textwidth]{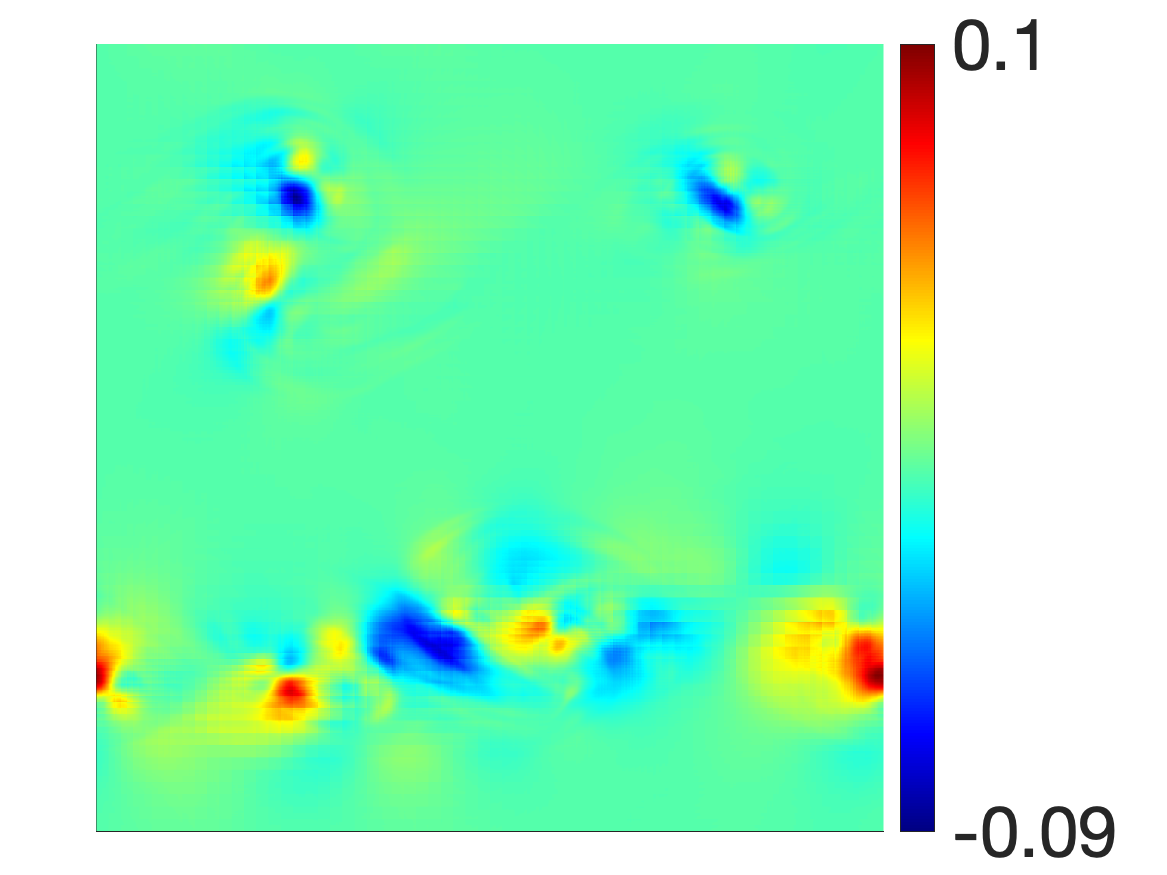}
\end{subfigure}
\begin{subfigure}{0.245\linewidth} \centering
\includegraphics[width=1\textwidth]{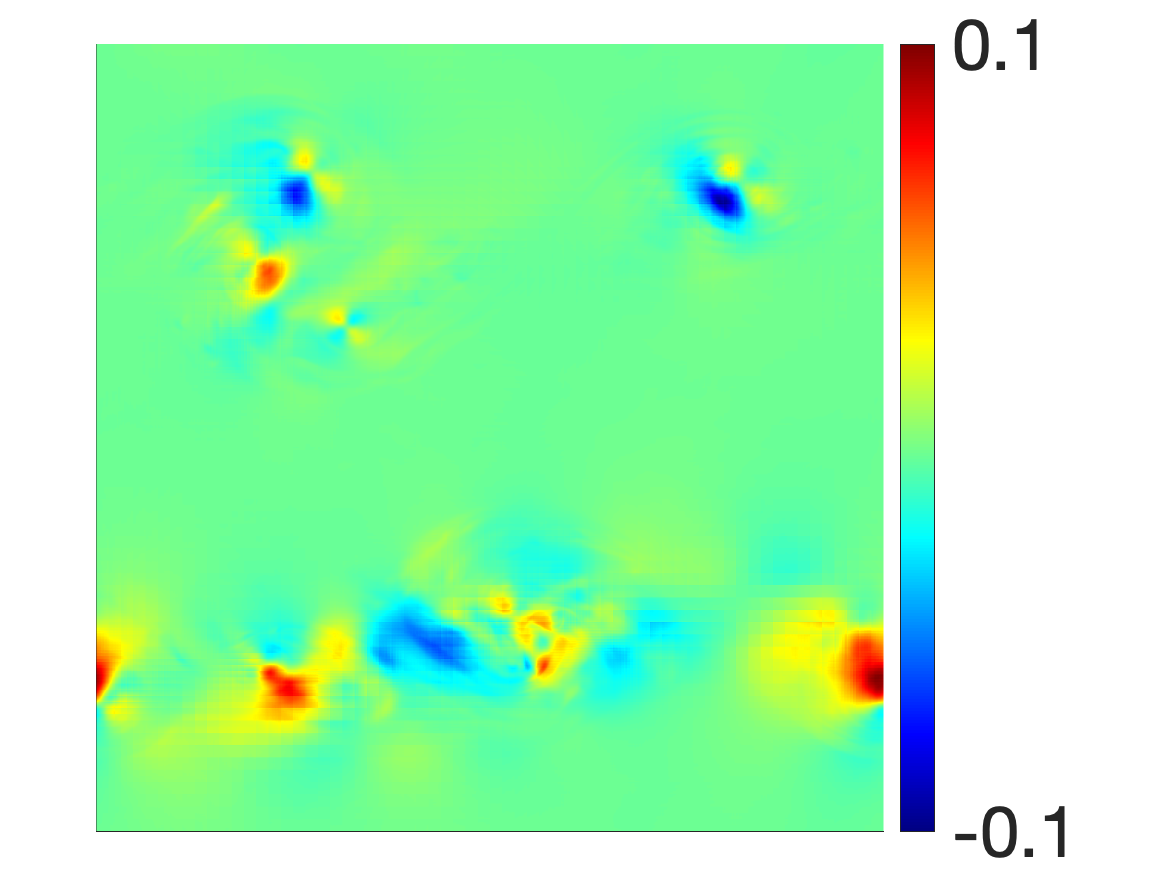}
\end{subfigure}
\\
\begin{subfigure}{0.245\linewidth} \centering
\includegraphics[width=1\textwidth]{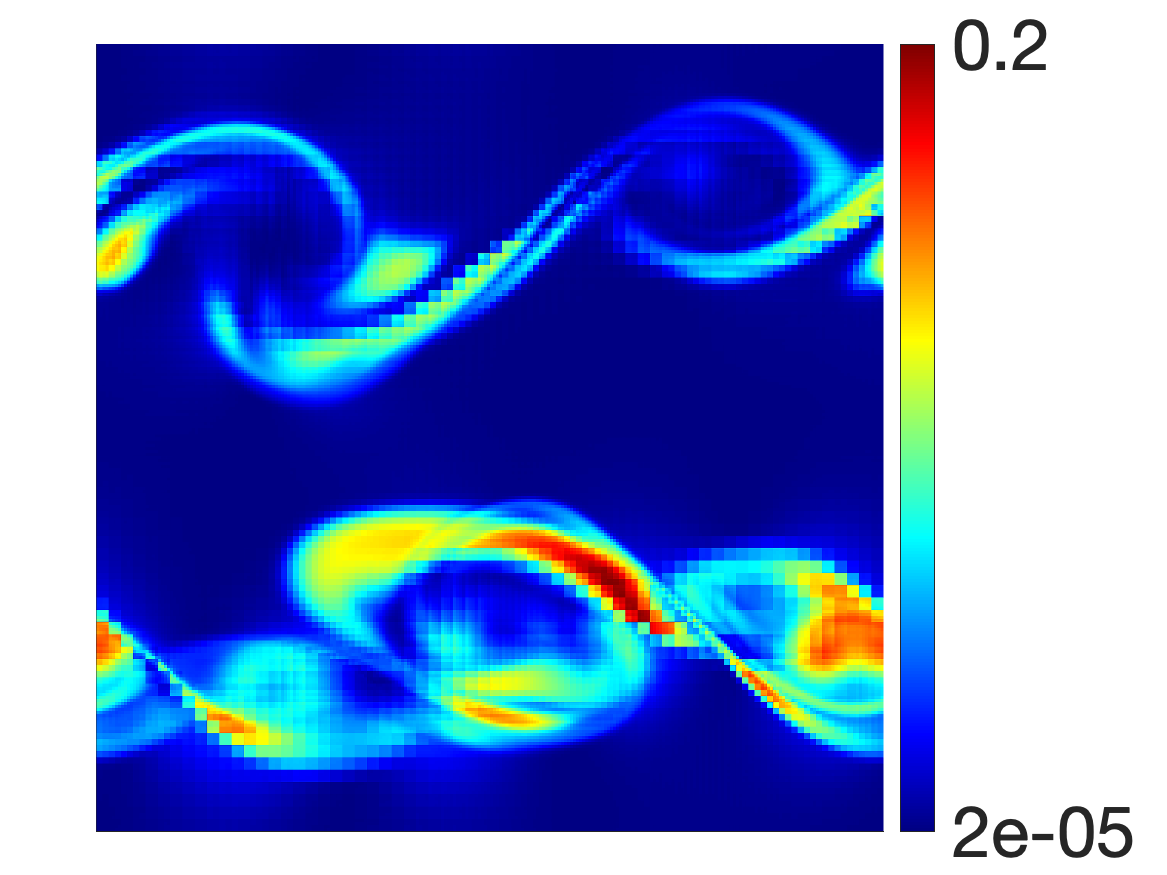}
\end{subfigure}
\begin{subfigure}{0.245\linewidth} \centering
\includegraphics[width=1\textwidth]{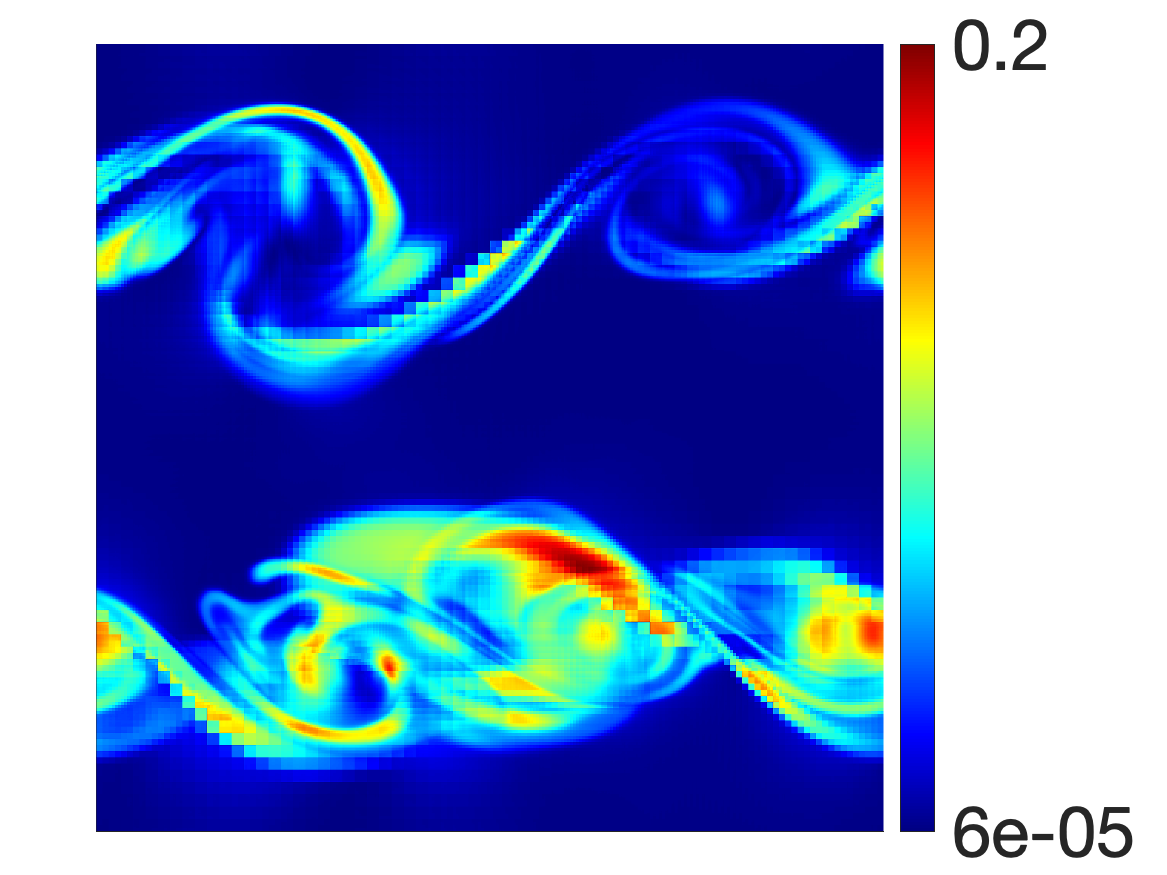}
\end{subfigure}
\begin{subfigure}{0.245\linewidth} \centering
\includegraphics[width=1\textwidth]{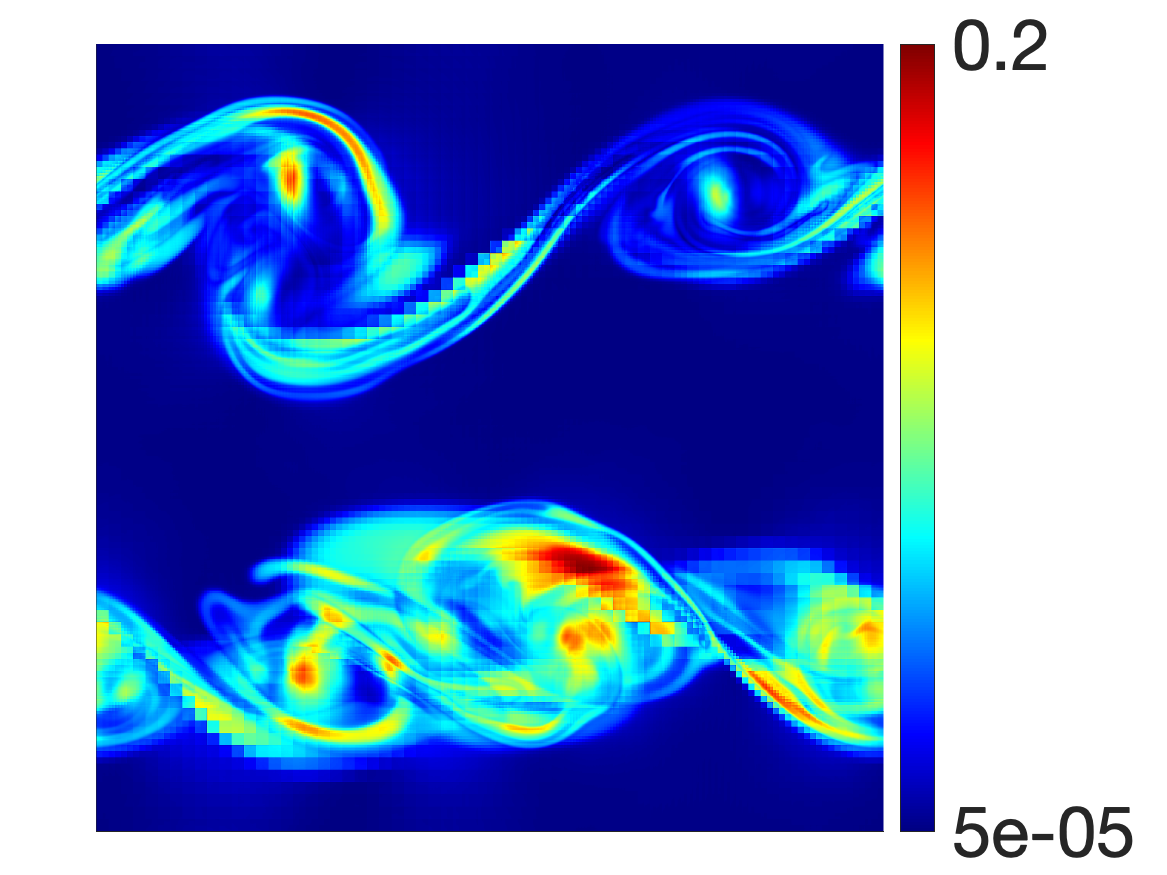}
\end{subfigure}
\begin{subfigure}{0.245\linewidth} \centering
\includegraphics[width=1\textwidth]{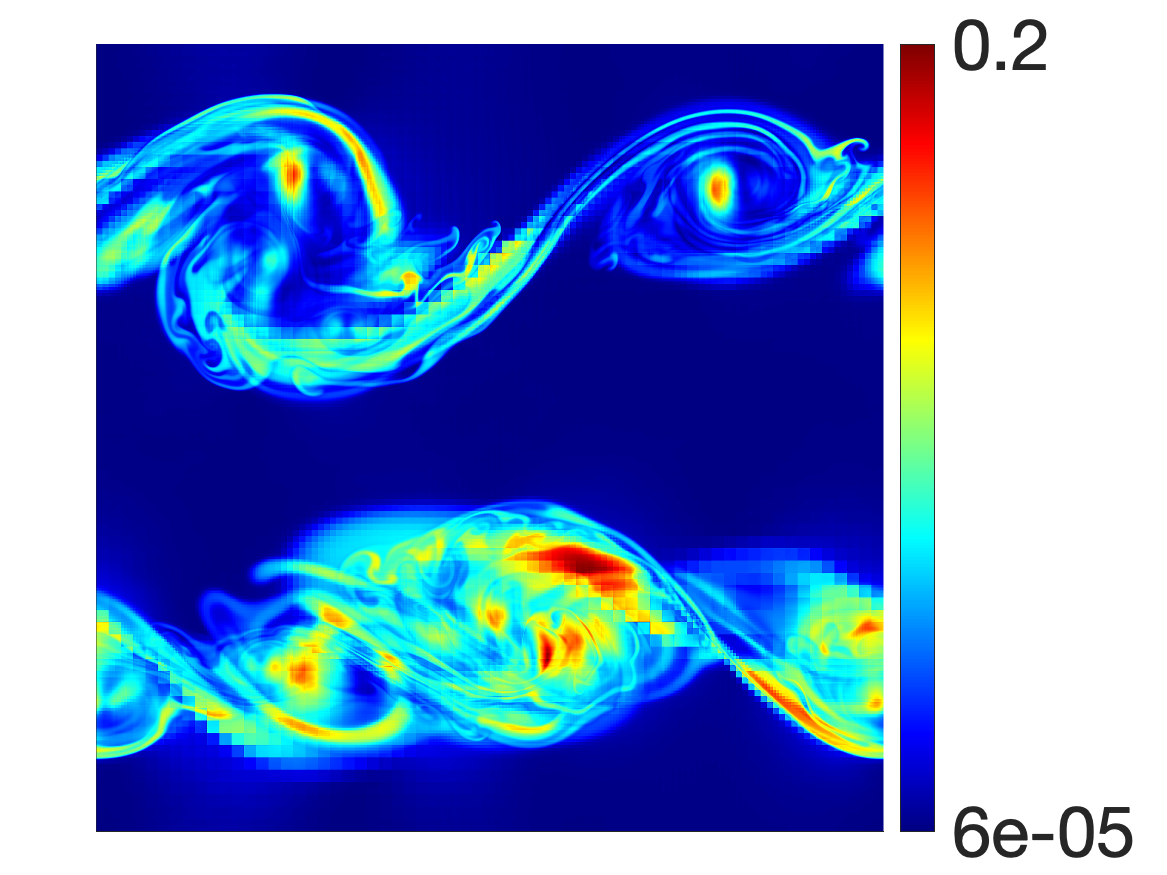}
\end{subfigure}
\\
\begin{subfigure}{0.245\linewidth} \centering
\includegraphics[width=1\textwidth]{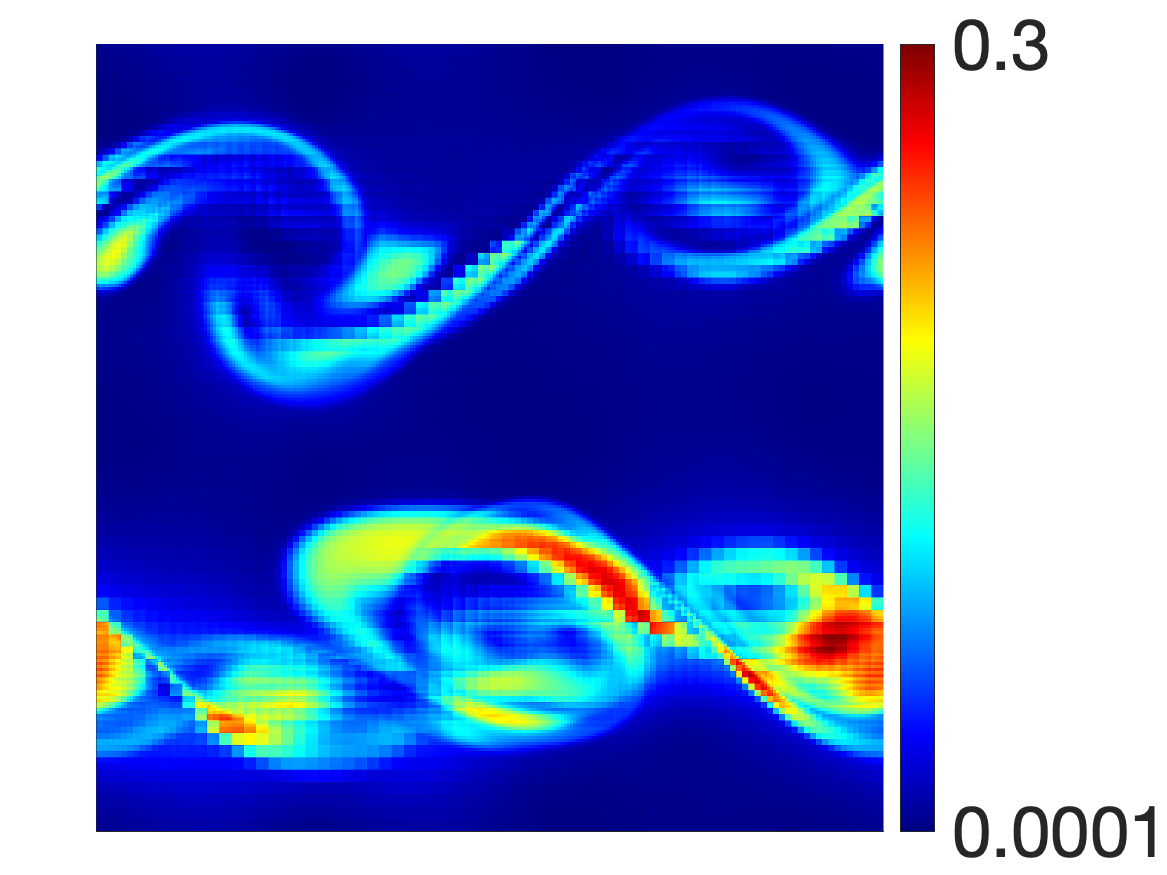}
\end{subfigure}
\begin{subfigure}{0.245\linewidth} \centering
\includegraphics[width=1\textwidth]{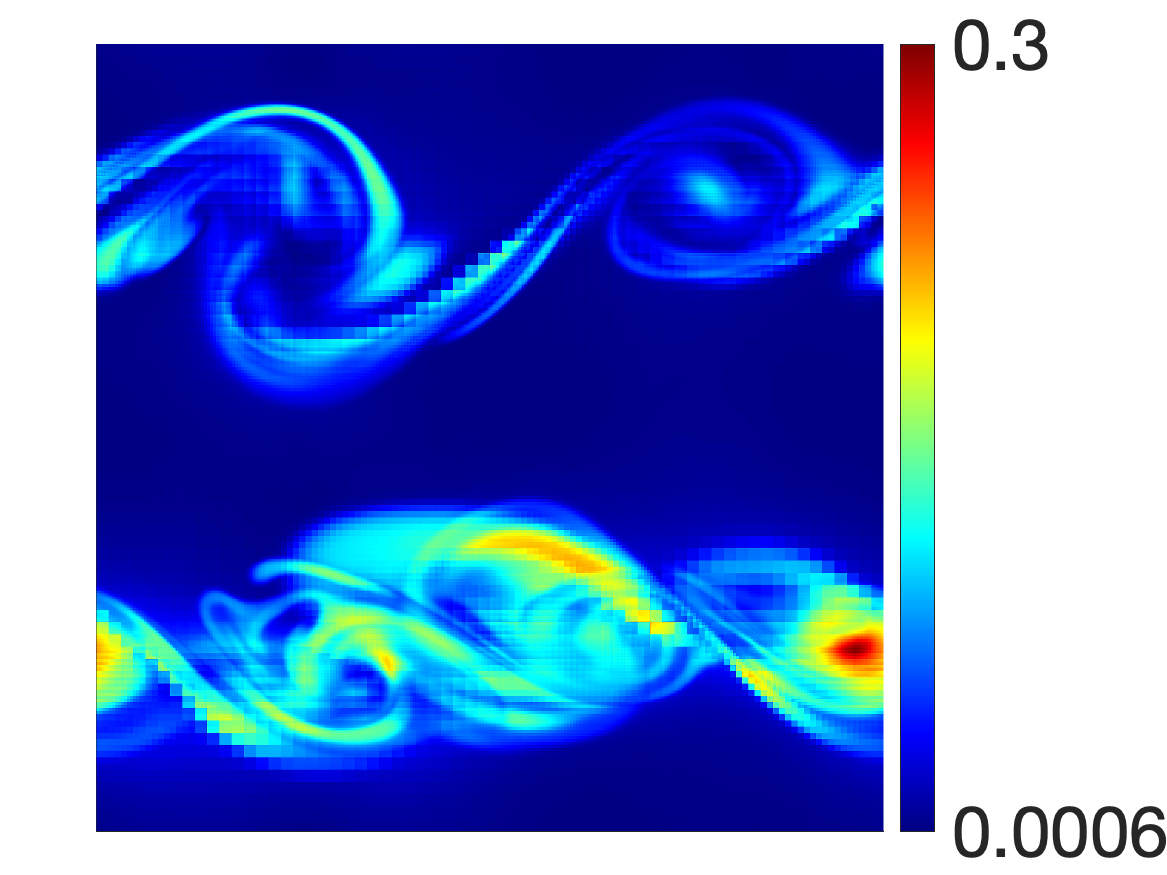}
\end{subfigure}
\begin{subfigure}{0.245\linewidth} \centering
\includegraphics[width=1\textwidth]{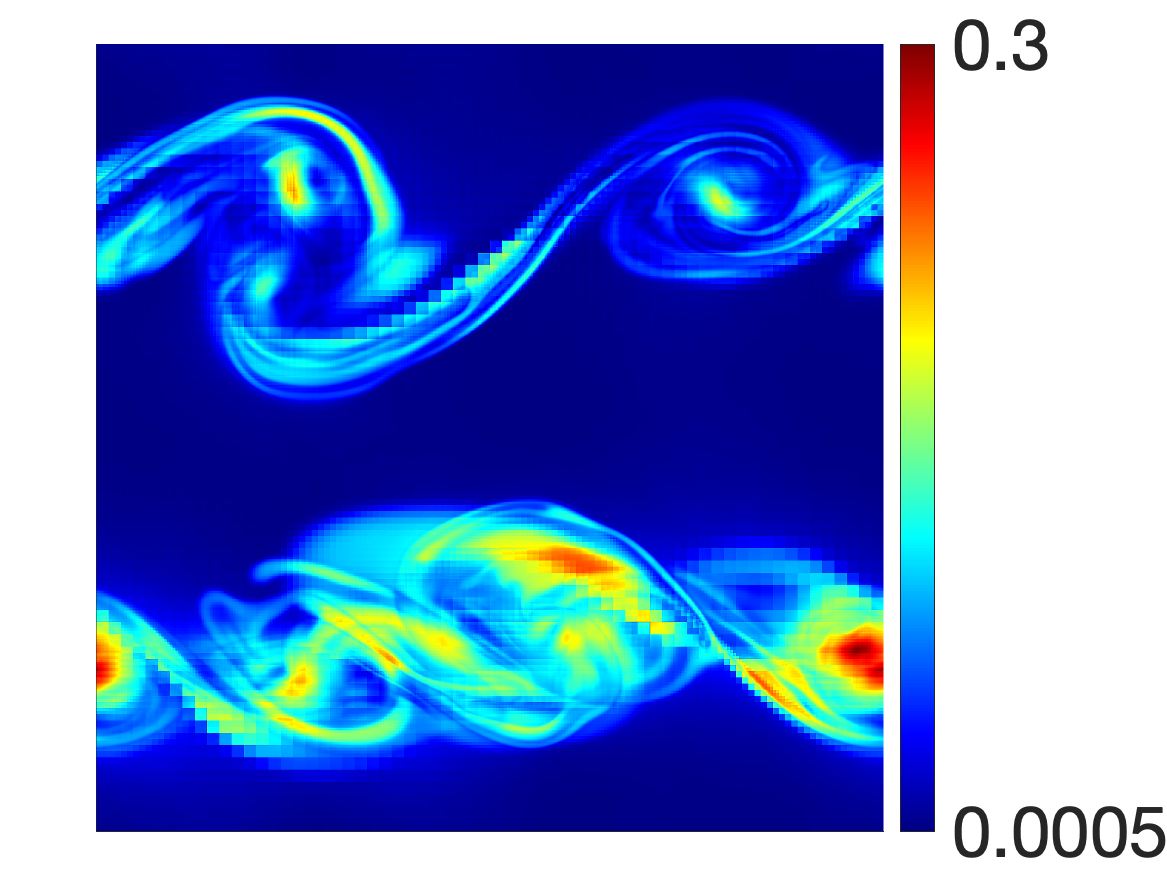}
\end{subfigure}
\begin{subfigure}{0.245\linewidth} \centering
\includegraphics[width=1\textwidth]{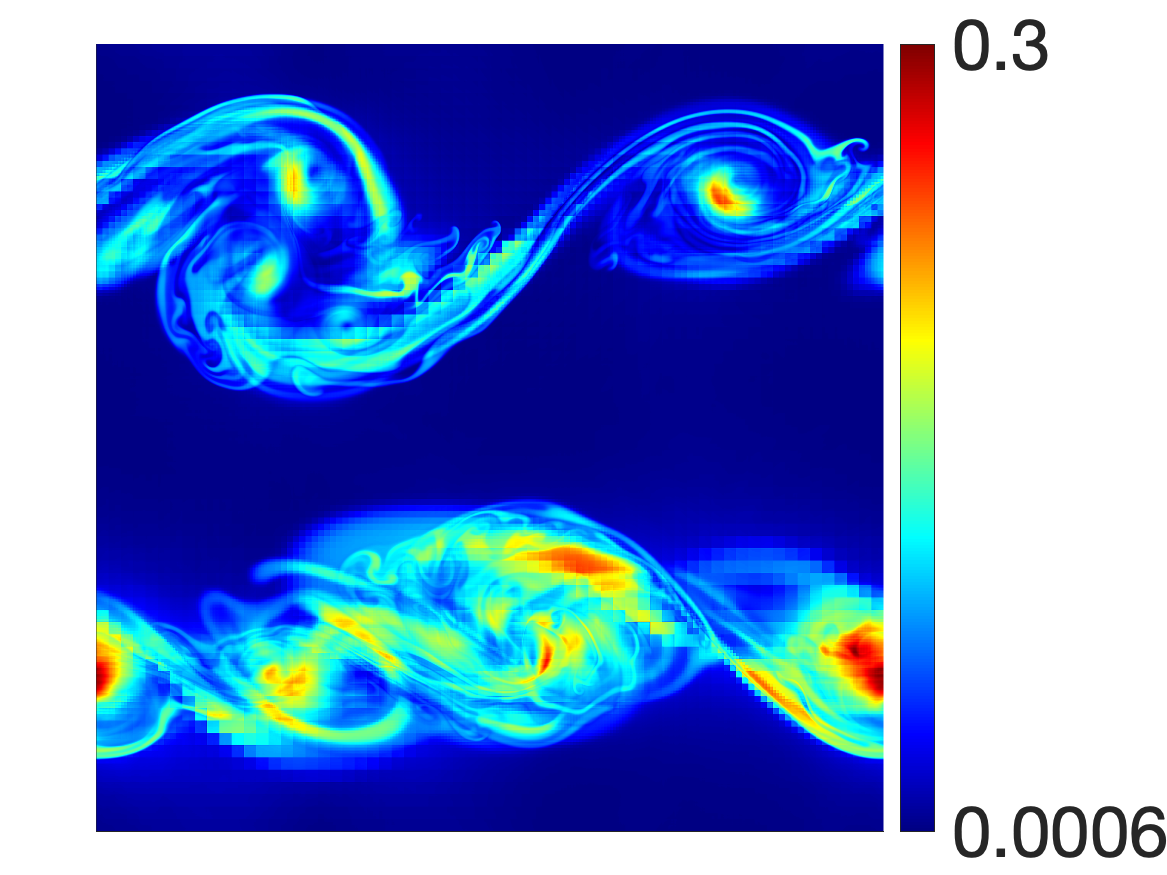}
\end{subfigure}
\\
\begin{subfigure}{0.245\linewidth} \centering
\includegraphics[width=1\textwidth]{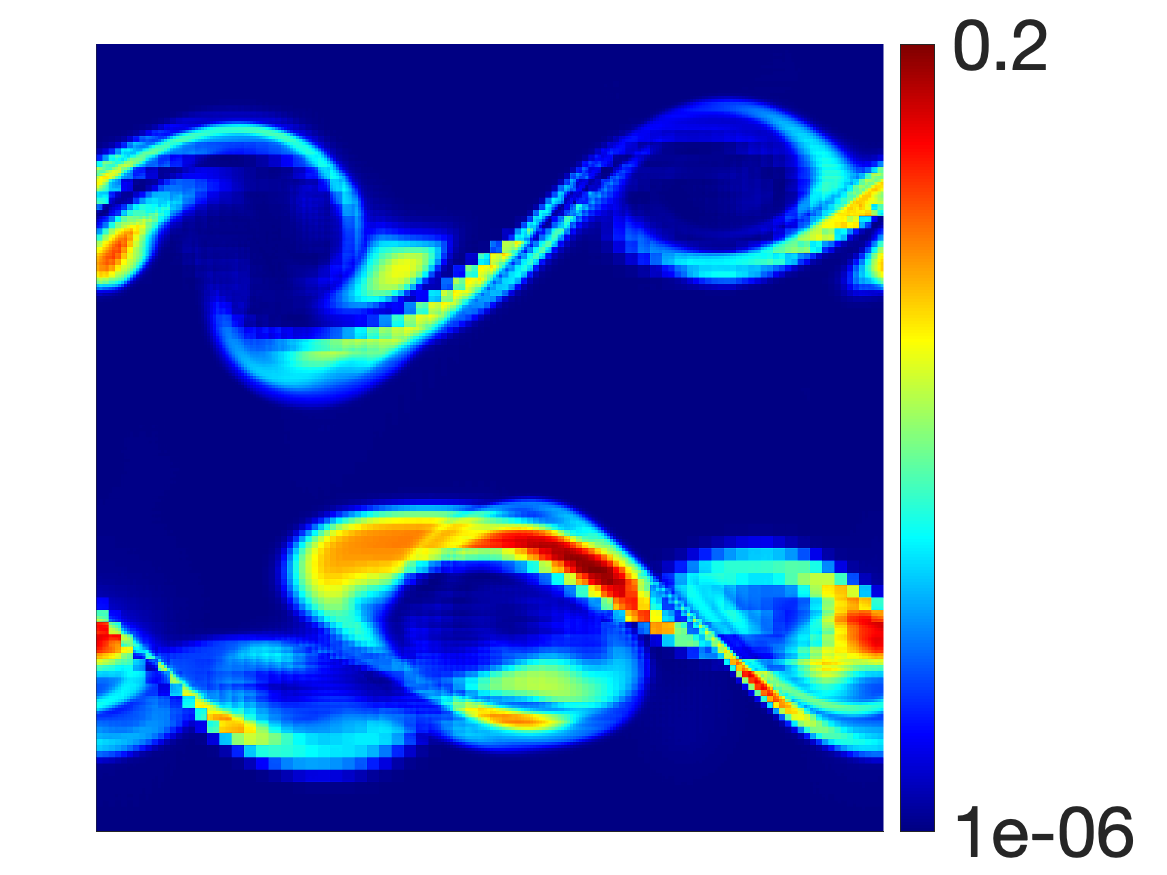}
\caption{$N=3$}
\end{subfigure}
\begin{subfigure}{0.245\linewidth} \centering
\includegraphics[width=1\textwidth]{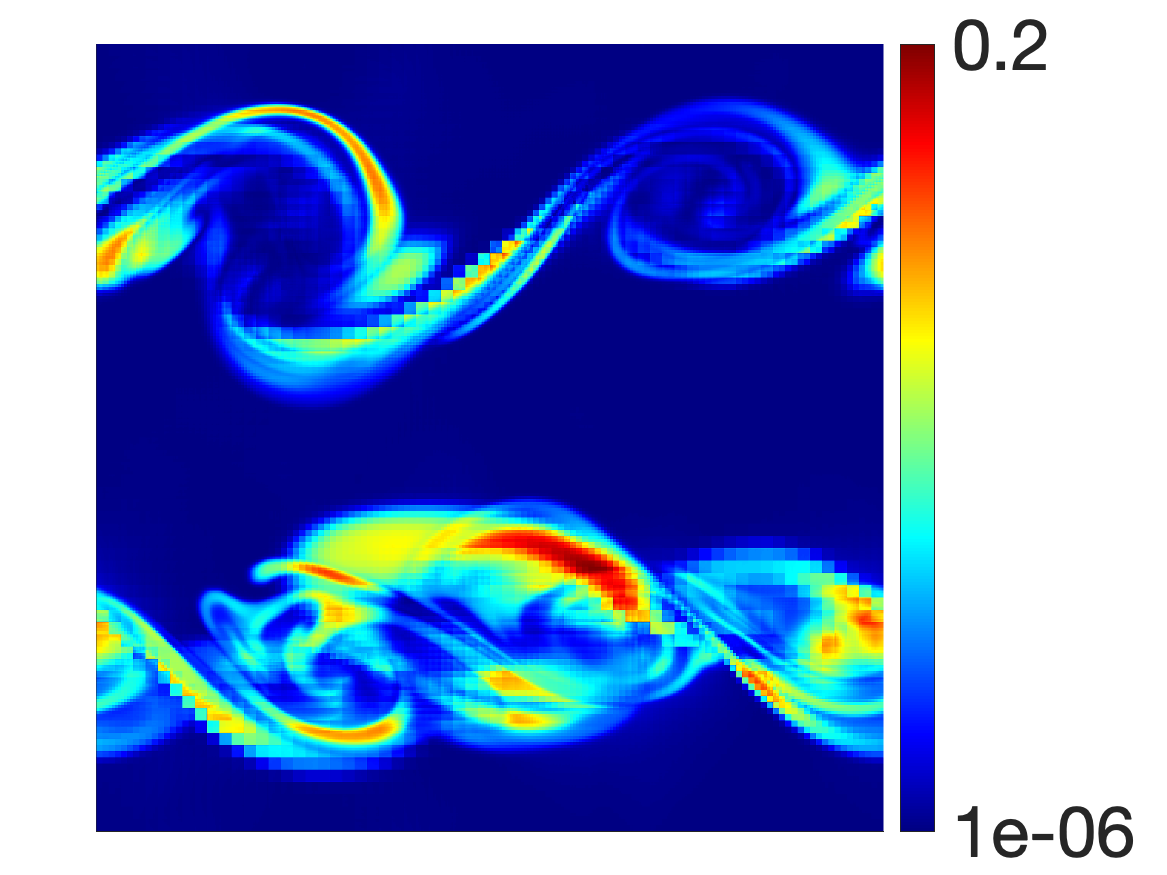}
\caption{$N=4$}
\end{subfigure}
\begin{subfigure}{0.245\linewidth} \centering
\includegraphics[width=1\textwidth]{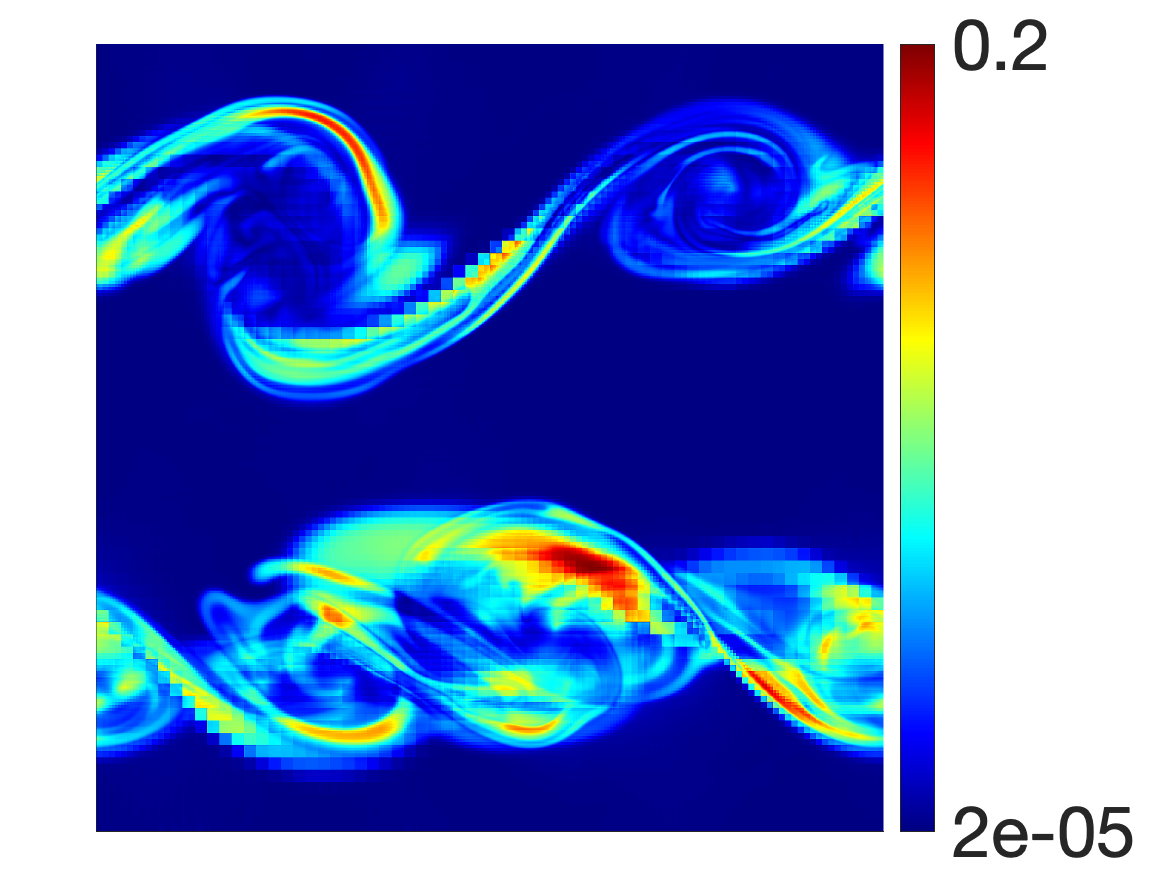}
\caption{ $N=5$}
\end{subfigure}
\begin{subfigure}{0.245\linewidth} \centering
\includegraphics[width=1\textwidth]{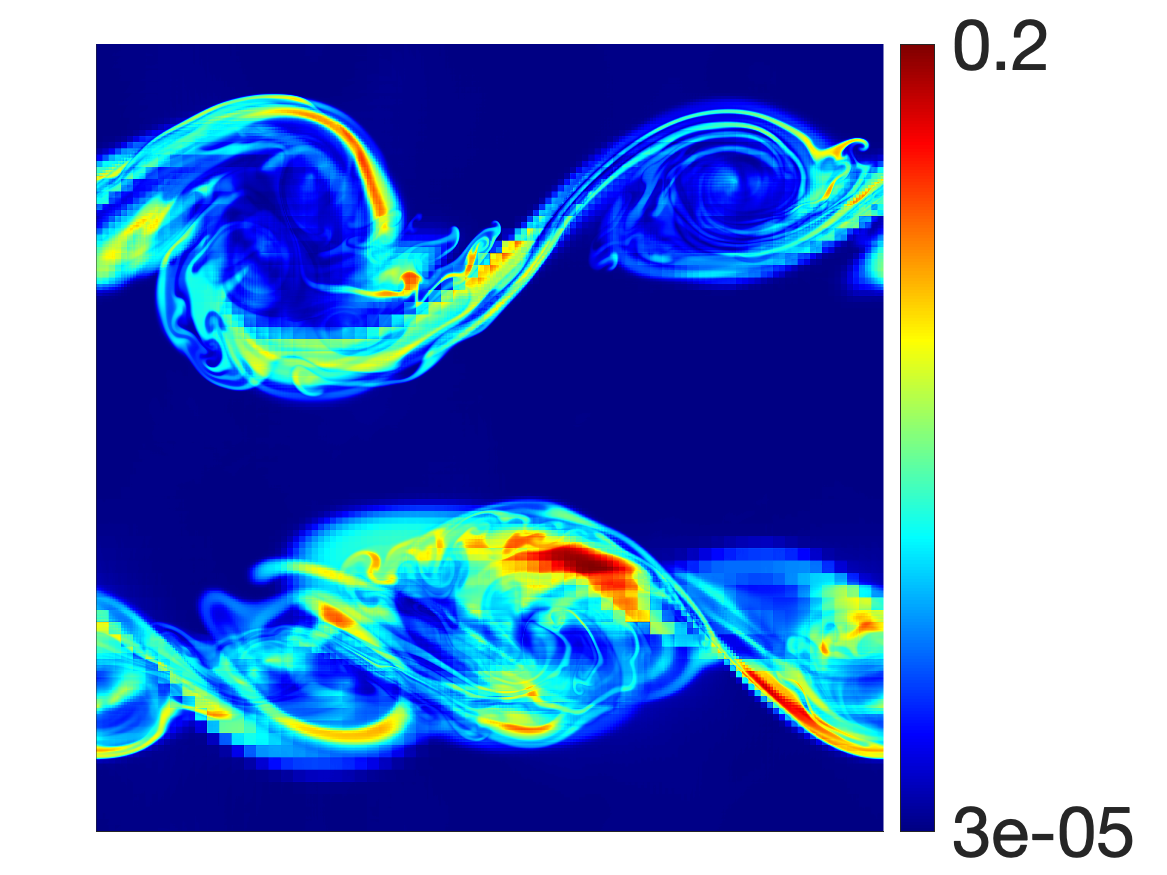}
\caption{$N=6$}
\end{subfigure}
\caption{Defects computed from the GRP solution.  From top to  bottom: $\mathfrak{E}_{N}$,  $\mathfrak{R}^{11}_N$, $\mathfrak{R}^{12}_N$, $\mathfrak{R}^{22}_N$,  $\lambda_1(\mathfrak{R}_N)$ $\lambda_2(\mathfrak{R}_N),$ $N=3, \dots, 6.$}
\label{fig_Rey_GRP}
\end{figure}

Next, in Fig.~\ref{fig-Rey-1} we present  the energy defects $\mathfrak{E}_N$ and the Reynolds stress defects  $\mathfrak{R}_N=( \mathfrak{R}^{11}_N,\mathfrak{R}^{12}_N;\mathfrak{R}^{12}_N,\mathfrak{R}^{22}_N)$ as well as their eigenvalues $\lambda_1(\mathfrak{R}_N), \lambda_2(\mathfrak{R}_N), \lambda_1(\mathfrak{R}_N) \geq \lambda_2(\mathfrak{R}_N)$, $N=3,4,5,6$, computed by the VFV method. Similar plots obtained from the GRP method are shown in Fig.~\ref{fig_Rey_GRP}.
As documented in Fig.~\ref{fig_err}  the Reynolds stress and energy defects do converge for $N \to \infty.$

For the VFV method we can observe that
{\sl i)}\,  the defects are mainly nonzero near the perturbed interface;
{\sl ii)}\, the Reynolds stress defects remain positive-definite;
{\sl iii)}\, the energy defects are positive, cf.~Fig.~\ref{fig-Rey-1}, where defects for different  $N$ are presented.
This perfectly confirms our theoretical results, cf.~Theorem~\ref{K}. Recall that the VFV method
has been proven to be weakly convergent to a DW solution with positive defects in \cite{FLM18_brenner}.
Fig.~\ref{fig_Rey_GRP} illustrates that the performance of the Reynolds stress and energy defects for the GRP method is similar, although the  convergence of the GRP method is not yet proved rigorously.


%

%
%

\vspace{-0.3cm}
\begin{acknowledgement}
This work was supported by  the Chinesisch-Deutschen Zentrum für Wissenschaftsförderung\chinese{(中德科学中心)} - Sino-German project number GZ1465. We thank Wang Yue (Beijing) for providing us with her simulation results of the GRP method.

M.L. gratefully acknowledges the support of Deutsche Forschungsgemeinschaft (DFG, German Research Foundation) - project number 233630050 - TRR 146 and  project number 525853336 -  SPP 2410 ``Hyperbolic Balance Laws: Complexity, Scales and Randomness''. She is grateful to the Gutenberg Research College and  Mainz Institute of Multiscale Modelling  for supporting her research.

The work of B.S. was supported by National Natural Science Foundation of China under grant No. 12201437.

Y.Y. gratefully acknowledges Nanjing University of Aeronautics and Astronautics  for supporting her research under the project number 90YAT23020.
\end{acknowledgement}

\end{document}